
\input eplain

\magnification\magstephalf

\def\lunghezza#1{\harrowlength#1pt
\varrowlength=0.618\harrowlength
\sarrowlength=1.18\harrowlength}

\def\diagram#1#2{\bgroup\lunghezza{#1} \commdiag{#2}\egroup}

\def\sqdiagram#1#2{\bgroup\harrowlength#1pt
\varrowlength=\harrowlength
\sarrowlength=0.42\harrowlength \commdiag{#2}\egroup}

\def\sequencelength{20pt}

\def\sequence#1{\bgroup\harrowlength\sequencelength
\varrowlength=0.618\harrowlength
\sarrowlength=0.42\harrowlength \commdiag{#1}\egroup}

\def\solleva#1#2{\raise#1\hbox{$#2$}}

\font\bigbold=cmbx12

\font\bigtype=cmr12
\font\biggtype=cmbx12 scaled \magstep 2
\font\Smcps = cmcsc10 scaled \magstep 1

\font\smallbulletseven=cmbsy10 scaled 1200
\font\smallbulletfive=cmbsy7 scaled 1200
\newfam\smallbulletfam
\textfont\smallbulletfam=\smallbulletseven
\scriptfont\smallbulletfam=\smallbulletseven
\scriptscriptfont\smallbulletfam=\smallbulletfive
\def\mini{{\fam\smallbulletfam \mathchar"7001}}

\font\gothicten=eufm10
\font\gothicseven=eufm7
\font\gothicfive=eufm5
\newfam\gothfam
\textfont\gothfam=\gothicten
\scriptfont\gothfam=\gothicseven
\scriptscriptfont\gothfam=\gothicfive
\def\goth{\fam\gothfam \gothicten}

\long\def\title#1{\bgroup
\biggtype \parindent=0pt\rightskip = 0in plus 10in
\leftskip=\rightskip \spaceskip=.3333em \xspaceskip=.5em
\parfillskip=0pt \baselineskip=21pt #1

\egroup\vskip .4cm}

\newcount\sezione

\newcount\enunciato

\enunciato=1

\newcount\lettera

\def\dotleaders{\leaders\hbox to1em{\hfil.\hfil}\hfil}

\def\iniziosezione{\vskip0pt plus.3\vsize\penalty-50\vskip0pt
plus-.3\vsize\bigskip\bigskip\vskip\parskip}

\def\beginsectionnonumber#1\par{\message{#1}
\iniziosezione
\leftline{\Smcps#1}\nobreak\medskip}

\def\beginsection #1[#2]{\message{#1}\advance\sezione by
1\enunciato=1
\iniziosezione
\leftline{\Smcps\the\sezione\enspace#1}
\definexref{#2}{\the\sezione}{section}
\definexref{item #2}{\leavevmode\hbox to 15pt{\the\sezione\hfil}
#1}{}
\xrdef{page #2}
\nobreak\medskip}

\outer\def\references{
\beginsection{References}[references]
\parindent=0pt \def\\{\smallbreak}}

\long\def\contents{\message{Table of contents}
\iniziosezione
\leftline{\Smcps Table of Contents}
\nobreak\medskip }

\def\indice#1{\line{\indent\ref{item #1}\rm\dotleaders\ref{page
#1}\hskip2cm}\smallskip}

\def\puntoletter{\ifcase\lettera \or a\or b\or c\or d\or e\or
f\or g\or h\or i\or j\or k\or l\or m\or n\or o\or p\or q\or r\or
s\or t\or u\or v\or w\or x\or y\or z\fi}

\def\proclaim #1. {\medbreak \def\ambiente{#1} \lettera=0
\noindent{\bf(\the\sezione{.}\the\enunciato)\enspace#1.}
\enspace\sl}

\def\call#1{\ifnum\lettera=0
\definexref{#1}{{\rm\the\sezione{.}\the\enunciato}}{\ambiente}
\definexref{bf #1}{{\bf\the\sezione{.}\the\enunciato}}{}
\else
\definexref{#1}{{\rm
\the\sezione{.}\the\enunciato{.}(\puntoletter)}}{\ambiente}
\definexref{part #1}{{\rm (\puntoletter)}}{}
\definexref{bf #1}{{\bf
\the\sezione{.}\the\enunciato{.}(\puntoletter)}}{}
\fi}

\def\endproofsymbol{\clubsuit}

\def\proclaimr #1. {\proclaim #1. \rm}

\def\punto
{\ifvmode\smallskip\leavevmode\else\relax\fi\advance\lettera by
1\hbox to 20pt{\rm (\puntoletter)\hfil}}

\def\endproclaim{\advance\enunciato by
1\rm\par\ifdim\lastskip<\medskipamount
\removelastskip\penalty55\medskip\fi}

\outer\def\proof{\smallbreak\noindent{\bf Proof.}\enspace}

\outer\def\proofof #1. {\smallbreak\noindent{\bf Proof of \ref{bf
#1}.}\enspace}

\outer\def\endproof{~\hfill~$\endproofsymbol$\par
\ifdim\lastskip<
\medskipamount\removelastskip\medbreak\fi}

\def\equaldef{\mathrel{\mathop=\limits^{\scriptscriptstyle\rm
def}}}

\def\into{\hookrightarrow}

\def\spec{\mathop{\rm Spec}}

\def\H{{\rm H}}

\def\cechH{\check{\rm H}}

\def\hom{\mathop{\rm Hom}\nolimits}

\def\curshom{\mathop{\mit{\cal H}om}\nolimits}

\def\gothext{{\goth Ext}}

\def\extcat{\gothext_{\cal O}({\cal F},{\cal G})}

\def\der{\mathop{\rm Der}\nolimits}

\def\m{{\goth \mathchar"716D}}

\def\a{{\goth \mathchar"7161}}

\def\cursext{\mathop{\mit{\cal E}\hskip-.5pt xt}\nolimits}

\def\ext{\mathop{\rm Ext}\nolimits}

\def\tor{\mathop{\rm Tor}\nolimits}

\def\pic{\mathop{\rm Pic}\nolimits}

\def\sym{\mathop{\rm Sym}\nolimits}

\def\loc{local complete intersection}

\def\im{\mathop{\rm im}\nolimits}

\def\cot #1{\m_{#1}/\m_{#1}^2}

\def\tang #1{(\m_{#1}/\m_{#1}^2)^\vee}

\def\ker{\mathop{\rm ker}\nolimits}

\def\cech{\v Cech}

\def\E{{\goth E}}

\def\obs{\mathop{\rm Obs}}

\def\indlim{\mathop{\vrule width0pt height7pt depth
4pt\smash{\lim\limits_{\raise 1.2pt\hbox {$\harrowlength13pt
\mapright$}}}}}

\def\projlim{\mathop{\vrule width0pt height7pt depth
4pt\smash{\lim\limits_{\raise 1.2pt\hbox {$\harrowlength13pt
\mapleft$}}}}}

\def \KS{Kodaira-Spencer}

\def\rest#1{\mathbin\mid\nobreak{}_{#1}}

\title{The deformation theory 

of local complete intersections}

\vskip .2cm

\centerline{\bigbold Fourth draft}

\vskip.5cm

\centerline{{\bigtype Angelo Vistoli}\footnote *{The first
version of this paper was written during a very pleasant stay at
the Mathematics Department of Harvard University. The author is
grateful for the hospitality.}}\vskip .3cm
\centerline {Dipartimento di Matematica}
\centerline {Universit\`a di Bologna}
\centerline{Piazza di Porta San Donato 5}
\centerline{40127 Bologna, Italy}
\smallbreak
\centerline{E-mail address: \tt vistoli@dm.unibo.it}

\vskip 1cm

This is an expository paper on the subject of the title. I tried
to make it as self-contained as I could, assuming only basic
facts in the theory of schemes and some commutative and
homological algebra, and to keep the formalism at a minimum.
There is very little pretense of originality on my part; the
results are of course very well known, and most of the ideas in
the proofs are also known. The possible exceptions are the
construction of the obstruction in
\ref{abstract}, which does not use simplicial techniques, and the
proof of the existence of versal deformations in \ref{versal},
in which the relations among the generators are obtained
directly from the obstructions.

The treatment is perhaps a little curt, and does not provide much
motivation for the results. Many of the details of the proofs
are left to the interested reader to fill in, and often I let
certain necessary compatibility conditions unstated. I believe
that not doing so would substantially increase the length of the
exposition without adding much to the reader's understanding;
however, a sufficient number of protests might change my mind
(so far, I got none.)

Comments and corrections are very welcome.

\beginsectionnonumber Acknowledgments

Joe Harris invited me to lecture on deformation theory, and then
encouraged me to write about it. His enthusiasm was infectious,
and I am grateful to him.

I have learned a large part of what I know about the subject from
Michael Artin. Very recently I followed a series of lectures
given by Hubert Flenner in Bologna, in July and October~1998; I
learned a lot from his lectures, and from discussion with him.

Lawrence Breen kindly answered my questions about the cohomology
of Picard stacks, and pointed out to me references [Retakh] and
[Ulbrich].

I learned a lot about deformation theory in a course that Hubert
Flenner taught in Bologna in July and October 1998, and from
discussions with him.

Finally, I thank all of those who helped me correct errors in
the previous versions of this paper, particularly Elham Izadi,
Ariane Mezard, Josko Plazonic and Kristina Rogale.

\contents

\indice{introduction}

\indice{embedded}

\indice{extensions}

\indice{abstract}

\indice{generalizations}

\indice{formal}

\indice{versal}

\indice{technical}

\indice{references}

\beginsection Introduction [introduction]

A very basic problem in algebraic geometry, which comes in all
sorts of guises, is to understand families of objects
(varieties, bundles, singularities, maps, \dots $\,$). This is
usually hard. The first step is to study the deformations of a
fixed object $X_0$ of the given type, that is, families of
objects depending on some parameters
$t_1,\dots,t_r$, such that for $t_1 = \cdots = t_r = 0$ we get
exactly $X_0$; we are interested in what happens near
$(0,\ldots,0)$. This is done in  three stages. 

First of all there is the problem of infinitesimal liftings,
which can be illustrated as follows. Suppose that we are given
an object $X_0$ over a field
$\kappa$, and a deformation $X$ of order $n-1$; we can think of
$X$ as obtained by perturbing the definition of $X_0$ by adding
a parameter $t$ with $t^n = 0$. This is very vague, of course;
in practice the definition usually involves defining
$X$ over the ring
$\kappa[t]/(t^n)$, or some other artinian ring, so that the
restriction of $X$ to $\kappa$ is $X_0$, and some other
``continuity'' condition (often involving flatness) is satisfied
(in practice it is very important to look also at deformations
on higher dimensional basis, that is, add more than one
parameter.) Then the problem is: can we lift $X$ to order $n$?
And if we can, how can we describe the liftings? Usually the
answer is in two parts: first there is a canonical element
$\omega$ of some vector space $V$ such that $\omega = 0$ if and
only if a lifting exists. Then there is some other vector space
$W$ such that if $\omega=0$ then
$W$ acts on the set of isomorphism classes of liftings making it
into a principal homogeneous space (with uncanny regularity $W$
is a cohomology group of a certain algebraic object, and $V$ the
cohomology group of the same object in one degree higher.)

The second stage is to look at formal deformations; a formal
deformation of $X_0$ can be roughly described as a lifting $X_1$
of $X_0$ to order 1, a lifting
$X_2$ of $X_1$ to order 2, and so, to all orders. Here the main
result is that, with very weak hypotheses, there is always a
formal deformation $V$ (involving several parameters,) called
``versal'', such that, very roughly, all other formal
deformations are obtained from $V$ (the actual definition is a
little technical.) This is defined over a quotient of a power
series algebra $R =
\kappa[[t_1,\ldots,t_r]]/I$ (or some other complete local ring,
when there is no base field); the $t_i$ are the parameters, and
$I$ the ideal of relations between the parameters. If the number
of parameters is as small as possible then saying that $I = 0$ is
equivalent to saying that the infinitesimal deformations are
unobstructed, that is, for any
$n$, given a deformation to order $n-1$ we can always lift it to
order
$n$. In practice knowing the ring
$R$ gives us a considerable degree of control on the
infinitesimal deformations.

The third stage is to pass from formal deformation to actual
deformations. In analytic geometry this amounts to passing from
formal solution of some equations to analytic solutions, and can
usually be done. In an algebraic context it is a much more
delicate question, sometimes called the problem of
algebraization; it can be solved for curves, but in general not
for surfaces.

Despite the fundamental importance of the subject I do not know
of any exposition that is both acceptably general and accessible
to the average algebraic geometer or number theorist. There is a
very thorough discussion of embedded deformations in [Koll\'ar],
and deformations of singularities are treated in a very readable
way in [Artin 1]; both are highly recommended. Going much beyond
these two references is [Illusie]; this is an excellent book,
and very carefully written, but uses substantial amounts of
simplicial machinery even to define the basic object, the
cotangent complex, and the exposition is at the kind of
topos-theoretic utmost level of generality that makes most
algebraic geometers' eyes glaze over in a fraction of the time
it takes to say ``simplicial object associated to a pair of
adjoint functors''.

Hopefully, this will change when a book by Buchweitz and
Flenner comes out in the next couple of years. It should treat
deformations of analytic spaces in full generality, using the
definition of the cotangent complex via Tate resolutions. These
are much easier to understand and work with then simplicial
resolutions; unfortunately at the moment this approach only works
in characteristic~0.

In this notes we study the infinitesimal and formal deformation
theory of a
\loc{} schemes. We limit ourselves to a few basic results that
can be obtained using sheaves of differentials, without
resorting to the cotangent complex, or to Tate resolutions, as
one has to do to go beyond this simple case. This is not likely
to keep most people happy for very long, but the task of giving a
reasonably self-contained explanation of the cotangent complex
is rather daunting to me.

The question of algebraization of formal deformations is only
touched upon very briefly (\ref{algdeformation}).

\ref{embedded} treats the infinitesimal liftings of \loc{}
subscheme of a given scheme.

Abstract liftings of generically smooth \loc{} schemes are
discussed in \ref{abstract}; here we use some standard, and some
less standard, facts on extensions of sheaves, proved in
\ref{extensions}. The construction of the obstruction seems to
be new, although I have not searched the literature long enough
to be sure.

\ref{generalizations} contains some generalizations, most
important to the case of \loc{} maps, and to the case of
deformations of pairs.

In \ref{formal} we define formal deformations of a \loc{}
generically smooth scheme over a field, define the \KS{} map and
the first obstruction map, and prove some basic properties.
Assuming the existence of the base field is absolutely not
necessary (a complete local ring would do very well) and is of
course a hindrance for applications to arithmetic, but it
simplifies to some extent the exposition.

In \ref{versal} we discuss versal deformations; in particular
prove the existence of versal deformations for generically
smooth \loc{} schemes over a field, using the results of
\ref{abstract}. This proof is different from the one in
[Schlessinger], and uses the obstruction theory. It has the
merit of illustrating how equations defining versal deformation
spaces arise from obstructions.

\ref{technical} contains the proof of a very important technical
lemma.

In \refs{formal} and \refn{versal} we become a little more
formal, and exploit the notion of homomorphism of deformations,
which allows to make the treatment somewhat less cluttered and
ultimately more transparent. Also, here rings are used, instead
of their spectra; I am aware that this might make many algebraic
geometers uncomfortable, but I feel that it is more natural in
this context.

The theory in these two sections can be extended to any of the
other cases considered in \ref{generalizations}, and beyond.
Indeed an axiomatic treatment would be possible (see for example
[Artin 2],) and perhaps will be added to a future version of
these notes.

\beginsection Liftings of embedded \loc s [embedded]

\proclaimr Notation. \call{notationring} The following notation
is used here and in
\ref{abstract}. Let
$A'$ be a noetherian local ring with maximal ideal $\m_{A'}$ and
residue field $\kappa = A'/\m_{A'}$, and $\a\subseteq A'$ an
ideal such that
$\m_{A'}\a = 0$; then $\a$ is a finite-dimensional vector space
over
$\kappa$. Set $A = A'/\a$. All schemes and morphisms will be
defined over $A'$.
\endproclaim

If $X$ is a scheme over a scheme $S$, and $T$ is a subscheme of
$S$, we denote by
$X\rest T$ the inverse image of  $T$ in $X$. If ``$X$'' denotes a
scheme over
$A$, where
$X$ is an arbitrary symbol, we will always set
$$X_0 = X\rest{\spec\kappa}.$$

Also, to avoid awkward terminology, if $X$ is a locally closed
subscheme of $M$, we will always talk about the sheaf of ideals
of
$X$ in $M$ meaning the sheaf of ideals of $X$ in some open
subset of
$M$ where $X$ is closed. More generally we will sometime omit to
say that $M$ has to be restricted to an open subset. This should
cause no confusion.

\proclaimr Hypotheses. \call{embhyp} Let $M'$ be a flat scheme of
finite type over
$A'$. Let $X$ be a \loc{} subscheme of $M = M'\rest{\spec A}$
(not necessarily closed); this means that locally the ideal of
$X$ in $M$ is generated by a regular sequence in ${\cal O}_M$.
Assume also that $X$ is flat over
$A$; then $X_0$ is still a \loc{} in $M_0$.

Call ${\cal N}_0$ the normal bundle to $X_0$ in $M_0$. If ${\cal
I}_0$ is the ideal of $X_0$ in $M_0$, then ${\cal C}_0 = {\cal
I}_0/{\cal I}_0^2$ is a locally free sheaf on
$X_0$, called the {\it conormal sheaf\/} to $X_0$ in $M_0$, and
${\cal N}_0$ is by definition its dual. \endproclaim

Let us check that indeed $X_0$ is still a \loc{} in $M_0$. This
is a particular case of the following lemma.

\proclaim Lemma. If $A\to B$ is a local homomorphism of local
noetherian rings, then
$X\times_{\spec A}\spec B$ is a \loc{} in $M\times_{\spec
A}\spec B$.
\endproclaim

\proof This is a local problem, so we may assume that $M$ is
affine,
$X$ is closed in
$M$, and  the ideal of $X$ in $M$ is generated by a regular
sequence {\bf f}. The Koszul complex
${\cal K}_\mini$ of {\bf f} is a resolution of ${\cal O}_X$ by
flat sheaves over
$A$; but ${\cal O}_X$ is by hypothesis flat over $A$, so ${\cal
K}_\mini\otimes_A B$ is exact in positive degree. This means
that the restriction of {\bf f} to $X\times_{\spec A}\spec B$ is
still a regular sequence, so $X\times_{\spec A}\spec B$ is a
\loc{} in $M\times_{\spec A}\spec B$.\endproof

\proclaim Definition.  A {\rm lifting} of $X$ to $M'$ is a
subscheme
$X'$ of
$M'$ which is flat over $A'$ and such that $X = X'\cap M$.
\endproclaim

\proclaim Theorem. \call{embdef} \punto \call{embdefa} Any
lifting of
$X$ is a \loc{} in $M'$.

\punto \call{embdefb} There is a canonical element
$$
\omega_{\rm emb} = \omega_{\rm emb}(X)\in
\a\otimes _\kappa\H^1(X_0,{\cal N}_0),
$$ called the\/ {\rm embedded obstruction} of $X$ in
$M$, such that $\omega_{\rm emb} = 0$ if and only if a lifting
exists.

\punto \call{embdefc} If a lifting exists, then there is a
canonical action of the group
$\a\otimes _\kappa\H^0(X_0,{\cal N}_0)$ on the set of liftings
making it into a principal homogeneous space.
\endproclaim

Let us begin the proof with a criterion for a given subscheme
$X'\subseteq M'$ with $X'\cap M = X$ to be flat over $A'$. Let
${\cal I}'$ be the ideal of
$X'$ in $M'$; then ${\cal I}'{\cal O}_M = {\cal I}$, so that
there is a natural surjective map ${\cal I}'/\a{\cal I}'\to {\cal
I}$.

\proclaim Lemma. \call{critflatness}The scheme $X'$ is flat over
$A'$ if and only if the natural surjective map ${\cal
I}'/\a{\cal I}'\to {\cal I}$ is an isomorphism.
\endproclaim

\proof This statement is local on $X'$, so we assume that all
schemes involved are affine and the embedding $X'\into M'$ is
closed. There is a short exact sequence
$$
\sequence{0&\mapright{\cal I}'&\mapright {\cal O}_{M'}&\mapright
{\cal O}_{X'}&\mapright&0};
$$ so if we tensor over $A'$ with $A$ we get an exact sequence
$$
\sequence{0 = \tor^{A'}_1({\cal O}_{M'},A)&\mapright
\tor^{A'}_1({\cal O}_{X'},A)&\mapright {\cal I}'/\a{\cal
I}'&\mapright {\cal O}_M&\mapright {\cal O}_{X}&\mapright 0}
$$ which shows that the map ${\cal I}'/\a{\cal I}'\to {\cal I}$
is an isomorphism if and only if $\tor^{A'}_1({\cal O}_{X'},A) =
0$. So if $X'$ is flat then the condition of the theorem is
verified; the converse statement is a particular case of
Grothendieck's local criterion of flatness (see for example
[Matsumura]), and can be proved very simply as follows. Let
$N$ be an arbitrary $A'$ module; we want to show that
$\tor^{A'}_1({\cal O}_{X'},N) = 0$, assuming this is true for $N
= A$. {}From the exact sequence of
$\tor$'s, it is enough to prove that $\tor^{A'}_1({\cal
O}_{X'},\a N) = \tor^{A'}_1({\cal O}_{X'},N/\a N) = 0$; so we may
assume that $\a N = 0$ (observe that $\a (\a N) = 0$.) In other
words, we assume that $N$ is an
$A$-module. Then the sequence
$$\sequence{0&\mapright {\cal I}'\otimes _{A'}N&\mapright {\cal
O}_{M'}\otimes _{A'}N&\mapright {\cal O}_{X'}\otimes
_{A'}N&\mapright 0}$$ is the same as the sequence
$$\sequence{0&\mapright {\cal I}\otimes_AN&\mapright {\cal
O}_M\otimes_AN&\mapright {\cal O}_X\otimes _AN&\mapright 0},$$
which is exact because $X$ is flat over $A$.
\endproof

Now we analyze the local situation; suppose that $M'$ is affine,
$X$ is closed in $M$, and the ideal
${\cal I}$ of $X$ in $M$ is generated by a regular sequence
$f_1,\dots,f_r$ in ${\cal O}_M$.

Let $X'$ be a lifting of $X$, ${\cal I}'$ the ideal of $X'$ in
$M'$. Choose liftings
$f'_1,\dots,f'_r$ of $f_1,\dots,f_r$ to ${\cal I}'$; then from
the equality ${\cal I}'/\a {\cal I}' = {\cal I}$ and the fact
that the ideal $\a$ is nilpotent we conclude that
$f'_1,\dots,f'_r$ generate ${\cal I}'$.

Let us check that  $f'_1,\dots,f'_r$ is a regular sequence in
${\cal O}_{M'}$; this will prove part
\ref{part embdefa}. Let ${\cal K}'_\mini$ be the Koszul complex
of 
$f'_1,\dots,f'_r$; then
${\cal K}_\mini = {\cal K}'_\mini \otimes _{A'} A$ is the Koszul
complex of $f_1,\dots,f_r$. We have a homology spectral sequence
$$E^2_{pq} = \tor^{A'}_p\bigl(\H_q({\cal
K}'_\mini),A\bigr)\Rightarrow\H_{p+q}({\cal K}_\mini) = \cases{0
&if
$p+q>0$\cr {\cal O}_X & if $p+q = 0$\cr}.$$ Notice that
$E^2_{p0} = 0$ for $p>0$, because
$\H_0({\cal K}'_\mini) = {\cal O}_{X'}$ is flat over $A'$.
{}From this, and the fact that the abutment is 0 in degree 1, we
get that
$\H_1({\cal K}'_\mini)\otimes_{A'}A = E^2_{01} = 0$. This
implies that
$\H_1({\cal K}'_\mini)$ is 0, and hence $E^2_{p1} = 0$ for all
$p$. Analogously one proves that
$\H_2({\cal K}'_\mini) = 0$, and by induction on $q$ that
$\H_q({\cal K}'_\mini) = 0$ for all
$q>0$. This proves that $f'_1,\dots,f'_r$ is a regular sequence.

Conversely, start from liftings $f'_1,\dots,f'_r$ of
$f_1,\dots,f_r$ to ${\cal O}_{M'}$, and define $X'$ via ${\cal
O}_{X'} = {\cal O}_{M'}/(f'_1,
\ldots, f'_r)$. Then we claim that $X'$ is a lifting of $X$ to
$M'$; for this we need to show that ${\cal I}'/\a{\cal I}' =
{\cal I}$ (\ref{critflatness}). Let
$\sum_i a'_if'_i$ be an element of
${\cal I}'$ whose image $\sum_ia_if_i$ in ${\cal I}$ is 0. Then
because $f_1,\dots,f_r$ is a regular sequence we can write
$(a_1,\ldots,a_r)\in {\cal O}_M^n$ as a linear combination of
standard relations of the form
$$ (0,\ldots,0,\hskip-6pt\underbrace{f_j}_{i^{\rm th}\ \rm
place}\hskip-6pt,0,\ldots,0,\hskip-6pt\underbrace{-f_i}_{j^{\rm
th}\
\rm place}\hskip-6pt,0,\ldots,0).
$$ These relations lift to relations
$$ (0,\ldots,0,\hskip-6pt\underbrace{f'_j}_{i^{\rm th}\ \rm
place}\hskip-6pt,0,\ldots,0,\hskip-6pt\underbrace{-f'_i}_{j^{\rm
th}\
\rm place}\hskip-6pt,0,\ldots,0).
$$ among the
$f'_i$. Then $(a'_1,\ldots,a'_r)$ can be written as a relation
among the $f'_i$, plus an element $(b'_1,\ldots,b'_r)\in
(\a{\cal O}_{M'})^n$, so that $\sum_i a'_if'_i = \sum_i
b'_if'_i\in
\a{\cal I}'$.

So we have proved that the liftings of $X$ are obtained locally
by lifting equations for
$X$. In particular we have proved the following.

\proclaim Lemma. \call{locnoobs} Assume that $X$ is affine and a
complete intersection in
$M$. Then $X$ has a lifting to $M'$.
\endproclaim

Here by a\/ {\it complete intersection} we mean that $X$ is
closed in
$M$, and its ideal is generated by a regular sequence.

Let $X'_1$ and $X'_2$ be two liftings of $X$ to $M'$; to these
we will associate a section
$\nu(X'_1,X'_2)$ of ${\cal N}_0$.

Call ${\cal I}'_1$ and ${\cal I}'_2$ the corresponding ideals in
${\cal O}_{M'}$. Take a local section $f$ of ${\cal I}$; then
$f$ can be lifted to sections
$f'_1$ and
$f'_2$ of ${\cal I}_1$ and ${\cal I}_2$ respectively. The
difference
$f'_1-f'_2$ is an element of
$\a{\cal O}_{M'} = \a\otimes_{A'}{\cal O}_{M'} =
\a\otimes_\kappa{\cal O}_{M_0}$. This element does not depend on
$f$ only; if we choose different liftings
$f'_1+g'_1$ and
$f'_2+g'_2$, with
$g'_i\in \a{\cal I}'_i$, then the difference $f'_1-f'_2$ will
change by an element
$g'_1-g'_2$ of $\a{\cal I}$; so the image of $f'_1-f'_2$ in
$\a\otimes_\kappa{\cal O}_{X_0}$ only depends on $f$. This
construction yields a function ${\cal I}\to
\a\otimes_\kappa{\cal O}_{X_0}$, which we denote by
$\nu(X'_1,X'_2)$. This function is
${\cal O}_M$-linear, and we think of it as a section of
$$
\eqalignno{\hom_{{\cal O}_M}({\cal I},\a\otimes_\kappa{\cal
O}_{X_0}) &= \hom_{{\cal O}_{M_0}}({\cal
I}_0,\a\otimes_\kappa{\cal O}_{X_0})  \cr &= \hom_{{\cal
O}_{X_0}}({\cal I}_0/{\cal I}_0^2,\a\otimes_\kappa{\cal
O}_{X_0}) \cr &= \H^0(X_0,\a\otimes_\kappa{\cal N}_0)\cr &=
\a\otimes_\kappa\H^0(X_0,{\cal N}_0).}
$$

This construction has the following properties.

\proclaim Proposition. \call{nucond} To each pair of liftings
$X'_1$,
$X'_2$ is associated a well defined element
$$\nu_{M'}(X'_1,X'_2) = \nu(X'_1,X'_2)\in
\a\otimes_\kappa\H^0(X_0,{\cal N}_0),$$ with the following
properties.

\punto \call{nuid} $\nu(X'_1,X'_2) = 0$ if and only if $X'_1 =
X'_2$.

\punto \call{nuadd} If $X'_1$, $X'_2$ and $X'_3$ are liftings,
then
$$\nu(X'_1,X'_3) = \nu(X'_1,X'_2)+\nu(X'_2,X'_3).$$

\punto \call{nuopp} $\nu(X'_2,X'_1) = -\nu(X'_1,X'_2)$.

\punto \call{nutrans} Given a lifting
$X'$ and a section $\nu\in \H^0(X,\a\otimes_A{\cal N}_0)$, there
is a lifting
$\widetilde X'$ with
$\nu(\widetilde X',X') = \nu$.

\punto \call{nuopen} If $Y$ is an open subscheme of $X$, $Y'_1$
and
$Y'_2$ are the restrictions of $X'_1$ and $X'_2$, then
$\nu(Y'_1,Y'_2)$ is the restriction of
$\nu(X'_1,X'_2)$.

\punto \call{nucomp} Let $\widetilde M'\to M'$ be a smooth
morphism of flat
$A'$-schemes of finite type, $X'_1\into
\widetilde M'$ and $X'_2\into\widetilde M'$ embeddings
compatible with the embeddings of
$X'_1$ and
$X'_2$ into
$M'$, inducing the same embedding of $X$ into $\widetilde M =
\widetilde M'\rest{\spec A}$. Let
$\widetilde {\cal I}_0$ be the ideal of $X_0$ in $\widetilde M_0
=
\widetilde M'\rest{\spec
\kappa}$. Then the homomorphism
$$\hom_{{\cal O}_{X_0}}(\widetilde {\cal I}_0/\widetilde {\cal
I}_0^2,\a\otimes_\kappa{\cal O}_{X_0})\longrightarrow
\hom_{{\cal O}_{X_0}}({\cal I}_0/{\cal
I}_0^2,\a\otimes_\kappa{\cal O}_{X_0})$$ induced by the natural
embedding of ${\cal I}_0/{\cal I}_0^2$ into
$\widetilde {\cal I}_0/\widetilde {\cal I}_0^2$ carries
$\nu_{\widetilde M'}(X'_1,X'_2)$ into
$\nu_{M'}(X'_1,X'_2)$.
\endproclaim

\proof We will only give a hint for part \ref{part nutrans}; the
remaining statements are straightforward and left to the reader.
The ideal $\widetilde {\cal I}'$ of $\widetilde X'$ in $M'$ can
be described as follows. A local section $\widetilde f'$ of
${\cal O}_{M'}$ is in $\widetilde {\cal I}'$ if and only if its
image $f$ in ${\cal O}_M$ lies in ${\cal I}$, and there exists a
local section
$f'$ of ${\cal I}'$ mapping to $f$ such that the image of
$\widetilde f'-f'$ in $\a\otimes_\kappa{\cal O}_{X_0}$ is
$\nu(f)$. One checks that $\widetilde {\cal I}'$ is an ideal in
${\cal O}_{M'}$, and that the subscheme $X'\subseteq M'$ with
ideal $\widetilde {\cal I}'$ is indeed a lifting of $X$ with
$\nu(\widetilde X',X') = \nu$.
\endproof

\ref{embdefc} follows immediately from \ref{nuid}, \ref{part
nuadd} and \ref{part nutrans}, and the following elementary fact.

\proclaim Lemma. \call{homogsp} Let $X$ be a set, $G$ a group.
Let there be given a function
$\phi\colon X\times X\to G$ with the following properties.

\punto $\phi(x_1,x_2) = 1$ if and only if $x_1 = x_2$.

\punto $\phi(x_1,x_2)\phi(x_2,x_3) = \phi(x_1,x_3)$ for all
$x_1$,
$x_2$ and $x_3$ in $X$.

\punto \call{homogsptrans} For each $g\in G$ and each $x\in X$
there exists $\widetilde x\in X$ such that
$\phi(\widetilde x,x) = g$.

\smallskip

Then the element $\widetilde x$ in \ref{part homogsptrans} is
unique, and $X$ has the structure of a principal left
homogeneous $G$-space, by defining $gx =
\widetilde x$ for all
$g\in G$ and $x\in X$.
\endproclaim

Let us prove part \ref{part embdefb} of the theorem. We may
assume that $X$ is closed in $M$. Choose a covering\/
${\cal U} =
\{U_\alpha\}$ of
$M$ by open affine subschemes, such that in each $U_\alpha$ the
subscheme $X_\alpha = X\cap U_\alpha$ is a complete
intersections, and call $U'_\alpha$ the corresponding open
subschemes of $M'$. By \ref{locnoobs}, each $X_\alpha$ has a
lifting
$X'_\alpha$ in
$M'$; there exists a global lifting if and only if after possibly
refining
${\cal U}$ we can choose the $X'_\alpha$ in such a way that
$X'_\alpha\cap(U'_\alpha\cap U'_\beta) =
X'_\beta\cap(U'_\alpha\cap U'_\beta)$ for all
$\alpha$ and $\beta$. To define the embedded obstruction
$\omega_{\rm emb}$, choose liftings
$X'_\alpha$ arbitrarily, and set
$$
\nu_{\alpha\beta} = \nu(X'_\alpha,X'_\beta).
$$ Of course this should have been written as
$$
\nu_{\alpha\beta} = \nu_{U'_\alpha\cap
U'_\beta}\bigl(X'_\alpha\cap(U'_\alpha\cap U'_\beta),
X'_\beta\cap(U'_\alpha\cap U'_\beta)\bigr);
$$ but now as in the future, we will commit a harmless and
convenient abuse of language by omitting to indicate the
restriction operators. Because of the cocycle relation of
\ref{nuadd} we see that
$\{\nu_{\alpha\beta}\}$ is a \v Cech 1-cocycle; a lifting exists
if and only if it is possible to choose local lifting so that
the associated cocycle is 0. If
$\widetilde X'_\alpha$ are different liftings, we set
$\nu_\alpha =
\nu(\widetilde X'_\alpha, X'_\alpha)$, $\widetilde
\nu_{\alpha\beta} =
\nu(\widetilde X'_\alpha,\widetilde X'_\beta)$. Again from
\ref{nuadd} and \ref{part nuopp} we get that
$$\widetilde \nu_{\alpha\beta} =
\nu_{\alpha\beta}+\nu_\alpha-\nu_\beta.$$ In other words,
cocycles associated to different local liftings are cobordant,
so the cohomology class
$\omega_{\rm emb}\in \check\H^1({\cal U},\a\otimes_A{\cal N}) =
\H^1(X,\a\otimes_A{\cal N})$ of
$\{\nu_{\alpha\beta}\}$ is independent of the local liftings.
Furthermore, if
$\widetilde \nu_{\alpha\beta}$ is a cocycle in $\omega_{\rm
emb}$, then there exists a 0-cochain $\{\nu_\alpha\}$ such that
$\widetilde \nu_{\alpha\beta} =
\nu_{\alpha\beta}+\nu_\alpha-\nu_\beta$. If we choose liftings
$\widetilde X_\alpha$ so that
$\nu_\alpha = \nu(\widetilde X'_\alpha,X'_\alpha)$
(\ref{nutrans}) we have $\widetilde
\nu_{\alpha\beta} = \nu(\widetilde X'_\beta,\widetilde
X'_\alpha)$. We can state this result as follows.

\proclaim Lemma. \call{descembcoc} The cocycles obtained from
different choices of local liftings are exactly the elements of
$\omega_{\rm emb}$. \endproclaim

So a lifting exists if and only if $\omega_{\rm emb} = 0$, as
desired. This concludes the proof of the theorem.\endproof

The proof of \ref{embdefb} can be described more conceptually as
follows. Consider the sheaf $\cal L$ of sets on
$X$, in which an open subset $U\subseteq X$ is sent to the set of
liftings of
$U$ to $M'$. Call $i\colon X_0 \into X$ the embedding. There is
an action of the the sheaf $i_*(\a\otimes_\kappa{\cal N}_0)$ on
$\cal L$;
\ref{embdefc} and
\ref{locnoobs} imply that $\cal L$ is a torsor for
$i_*(\a\otimes_\kappa{\cal N}_0)$. A lifting of $X$ to
$M'$ is global section of this torsor, so the obstruction to the
existence of a global section is the class $\omega_{\rm emb}$ of
$\cal L$ in $\H^1\bigl(X,i_*(\a\otimes_\kappa{\cal N}_0)\bigr) =
\a\otimes_\kappa
\H^1(X_0, {\cal N}_0)$.

The following property of the embedded obstruction $\omega_{\rm
emb}$ follows from its construction and from \ref{nucomp}.

\proclaim Lemma. \call{embobscomp}  Let $\pi\colon \widetilde
M'\to M'$ be smooth morphism of flat $A'$-schemes of finite
type, $j\colon X\into M'$ and
$\widetilde
\jmath\colon X\into \widetilde M'$ embeddings such that
$\pi\widetilde
\jmath = j$. Call
${\cal N}_0$ and
$\widetilde {\cal N}_0$ the normal bundles of $X_0$ in $M_0$ and
$\widetilde M_0$ respectively, and $h\colon H^1(X_0,\widetilde
{\cal N}_0)\to
\H^1(X_0,{\cal N}_0)$ the homomorphism obtained from the map
$\widetilde {\cal N}_0\to {\cal N}_0$ induced by $\pi$. Then $h$
carries the embedded obstruction of $X$ in
$\widetilde M'$ to the embedded obstruction of $X$ in $M'$.
\endproclaim

This construction has an obvious property of functoriality. Let
$B'$ be a local ring, ${\goth b}\subseteq B'$ be an ideal with
$\m_{B'}{\goth b} = 0$, $B = B'/{\goth b}$. Let $f\colon A'\to
B'$ a local homomorphism inducing an isomorphism of residue
fields, such that
$f(\a)\subseteq {\goth b}$. Set $f_*M' = M'\times_{\spec {A'}}
\spec {B'}$, $f_*M = M\times_{\spec {A}} \spec {B}$, $f_*X =
X\times_{\spec {A}} \spec {B}\subseteq f_*M$. If $X'$ is a
lifting of
$X$ in $M'$ set
$f'_*X' = X'\times_{\spec {A'}} \spec {B'}\subseteq f_*M'$; this
is a lifting of $f_*X$ in
$f_*M'$. Call $g = f\rest{\a}\colon \a\to {\goth b}$ the
restriction of $f$.

\proclaim Proposition. \punto \call{nufunc} Let $X'_1$ and
$X'_2$ be liftings of $X$ to $M'$. Then
$$\nu_{f_*M'}(f_*X'_1,f_*X'_2) = (g\otimes {\rm
id})\bigl(\nu_{M'}(X'_1,X'_2)\bigr)\in {\goth b}\otimes
\H^0(X_0,{\cal N}_0).$$

\punto \call{embobsfunc} $\omega_{\rm emb}(f_*X) = (g\otimes {\rm
id})\omega_{\rm emb}(X)\in {\goth b}\otimes \H^1(X_0,{\cal
N}_0)$.
\endproclaim

The case of local complete intersections
in projective spaces is particularly simple.

\proclaim Proposition. Let $X \subseteq {\bf P}^N_A$ a
complete intersection subscheme of codimension $r$, whose ideal
is generated by homogeneous polynomials $f_1, \ldots, f_r$. Then
any lifting of $X$ to ${\bf P}^N_{A'}$ is also a complete
intersection subscheme, with ideals generated by homogeneous
polynomials $f'_1, \ldots, f'_r$ which reduce to  $f_1, \ldots,
f_r$ modulo $\a$.
\endproclaim

\proof Fix liftings $f'_1, \ldots, f'_r$ of $f_1, \ldots,
f_r$ to homogeneous polynomials with coefficients in $A'$; these
define a lifting $X'$ of $X$ to ${\bf P}^N_{A'}$. If
we call $d_1, \ldots, d_r$ the degrees of $f_1, \ldots, f_r$,
then the normal bundle ${\cal N}_0$ decomposes as a direct sum
$\sum_{i=1}^r {\cal O}_{X_0}(d_i)$, so a section of $\a
\otimes_\kappa {\cal N}_0$ is given by a sequence $(g_1,
\ldots, g_r)$ of homogeneous polynomials of degrees $d_1,
\ldots, d_r$ with coefficients in $\a$. If $\widetilde X'$ is
the subscheme of ${\bf P}^N_{A'}$ whose ideal is generated by
$f'_1+g_1, \ldots, f'_r+g_r$, an analysis of the construction of
$\nu(\widetilde X',X')$ reveals that $\nu(\widetilde X',X')$ is
precisely $(g_1, \ldots, g_r)$. The result follows from
\ref{embdefc}.\endproof

Here is a typical application of \ref{embdef}.

\proclaim Corollary. Let $\pi\colon M\to S$ be a projective
morphism, where $S$ is a locally noetherian scheme. Let
$s_0\in S$ be a point, $M_0 = \pi^{-1}(s_0)$. Let $X_0\subseteq
M_0$ be a closed \loc{} subscheme with normal bundle ${\cal
N}_0$; assume that 
$\H^1(X_0,{\cal N}_0) = 0$. Then there is an \'etale neighborhood
$s_0\in U\to S$ of
$s_0$ and a closed subscheme $X\subseteq U\times_S M$ flat over
$U$ with $X\cap M_0 = X_0$.

Furthermore if $\H^0(X_0,{\cal N}_0) = 0$, then for any other
\'etale neighborhood $s_0\in U'\to S$, with a closed subscheme
$X'\subseteq U'\times_S M$ flat over $U'$ with
$X'\cap M_0 = X_0$, there exists a third \'etale neighborhood
$s_0\in U''\to S$ with morphisms of neighborhoods
$\phi\colon U''\to U$ and $\phi'\colon U''\to U'$ with
$\phi^{-1}(X) =
\phi'^{-1}(X')$.
\endproclaim

The basic example of this type of situation is a smooth
geometrically connected rational curve $X_0$ with self
intersection $-1$ on a smooth algebraic surface $M_0$.

\proof Assume that $\H^1(X_0,{\cal N}_0) = 0$. Let $H = {\rm
Hilb}(M/S)\to S$ be the relative Hilbert scheme, $\xi_0\in
H\bigl(\kappa(s_0)\bigr)$ the point corresponding to
$X_0\subseteq M_0$. Let $\spec B\into
\spec A$ be a closed embedding of spectra of local artinian
rings, with closed point $u_0\in
\spec B\into \spec A$, and suppose that it is given a commutative
diagram
$$
\diagram{30}{\spec B&\mapright^\beta&H\cr
\mapdown&&\mapdown\cr
\spec A&\mapright^\alpha&S\cr}
$$ such that $\beta(u_0) = \xi_0$, and hence $\alpha(u_0) =
s_0$. The morphism $\beta$ corresponds to a subscheme
$X_B\subseteq M\times_S \spec B$ flat over
$\spec B$ such that
$X_B\rest{u_0} = X_0\times_{s_0}u_0$. It follows easily from
\ref{embdefb}, by induction on the length of the kernel of the
homomorphism $A\to B$, that there exists a subscheme
$X_A\subseteq M\times_S \spec A$ flat over $\spec A$ such that
$X_A\rest{\spec B} = X_B$; furthermore from \ref{embdefc} we get
that if $\H^0(X_0,{\cal N}_0) = 0$ then $X_B$ is unique. This
means that there exists a morphism $\spec A\to H$ making the
diagram
$$\harrowlength=30pt
\diagram{30}{
\spec B&\mapright^\beta&H\cr
\mapdown&\arrow(3,2)&\mapdown\cr
\spec A&\mapright^\alpha&S\cr}
$$ commutative; moreover if $\H^0(X_0,{\cal N}_0) = 0$ this
morphism is unique. By Grothendieck's criteria this means that
$H$ is smooth at $\xi_0$, and if $\H^0(X_0,{\cal N}_0) = 0$ then
it is
\'etale. This implies that there exists an \'etale neighborhood
$s_0\in U\to S$ and a section $U\to H$ sending $s_0$ to $\xi_0$.
By taking for $X\subseteq U\times_S M$ the pullback of the
universal subscheme of $H\times_S M$ we have proved the first
statement.

For the second statement choose a Zariski neighborhood $\xi_0\in
H'\subseteq H$ which is
\'etale over $S$; by restricting $U$ and $U'$ we may assume that
the morphisms $U\to H$ induced by $U'\to H$ induced by the
subschemes
$X\subseteq U\times_S M$ and $X\subseteq U\times_S M$ have their
image in $H'$. Then we can take $U'' = U\times_{H'}U'$.\endproof

It is easy to give examples in which the \'etale neighborhood
$U$ can not be taken to be a Zariski neighborhood.

\vfill\break

\beginsection Extensions of sheaves [extensions]

In this section we will discuss briefly the theory of extensions,
which we will use to prove \ref{absdef}; I advise the reader to
skip it at first and then refer back to it as necessary. Working
directly with extensions, instead of elements of groups of
extensions, is critical in \ref{abstract}, because extensions
can be patched together, unlike classes in
$\ext^1$.

Let $X$ a topological space with a sheaf of commutative rings
${\cal O}$; in this section a {\it sheaf\/} will always be a
sheaf of ${\cal O}$ modules over $X$, and all homomorphisms will
be homomorphisms of sheaves of ${\cal O}$-modules. More
generally we could work with objects of a fixed abelian category.
I hope not to insult the reader by including some very standard
definitions.

\proclaim Definition. Let ${\cal F}$ and ${\cal G}$ be sheaves.
An\/ {\rm extension} $({\cal E},\iota,\kappa)$ of ${\cal F}$ by
${\cal G}$ is a sheaf ${\cal E}$ with two homomorphisms
$\iota\colon {\cal G}\to {\cal E}$ and
$\kappa\colon {\cal E}\to {\cal F}$ such that the sequence
$$
\sequence{0&\mapright& {\cal G}&\mapright^{\displaystyle\iota}&
{\cal E}&\mapright^{\displaystyle\kappa}& {\cal F}&\mapright &0}
$$ is exact.

If $({\cal E}_1,\iota_1,\kappa_1)$ and $({\cal
E}_2,\iota_2,\kappa_2)$ are extensions
$({\cal E},\iota,\kappa)$ of ${\cal F}$ by ${\cal G}$, a
homomorphism of extension
$$\phi\colon ({\cal E},\iota,\kappa)\to ({\cal
E}',\iota',\kappa')$$ is a homomorphism of sheaves $\phi\colon
{\cal E}\to {\cal E}'$ such that
$\phi\iota = \iota'$ and $\kappa'\phi = \kappa$.
\endproclaim

We will often talk about an extension ${\cal E}$, omitting the
homomorphisms $\iota$ and
$\kappa$ from the notation; if we need to refer to the them we
will call them $\iota_{\cal E}$ and $\kappa_{\cal E}$.

The content of the five lemma is that any homomorphism of
extensions is an isomorphism. With the obvious definition of
identities and composition of arrows, extensions form a category,
which is a groupoid, i.e., all arrows are invertible. We will
denote this category by
$\extcat$.

Call $\ext^1_{\cal O}({\cal F},{\cal G})$ the set of isomorphism
classes of extensions of $\cal F$ by $\cal G$. The link with the
usual definition of
$\ext^1_{\cal O}({\cal F},{\cal G})$ via injective resolutions
is as follows. Take an injective sheaf of
$\cal O$-modules
$\cal J$ containing
$\cal G$, and set ${\cal Q} = {\cal J}/{\cal G}$. Call
$\ext_{\cal O}({\cal F},{\cal G})$ the quotient of $\hom_{\cal
O}({\cal F},{\cal Q})$ by the subgroup of homomorphism
${\cal F}\to {\cal Q}$ which can be lifted to homomorphisms
${\cal F}\to {\cal J}$. Let $\cal E$ be an extension of $\cal F$
by $\cal G$. The embedding ${\cal G}\into {\cal J}$ can be
extended to a homomorphism ${\cal E}\to {\cal J}$, because
$\cal J$ is injective, which will induce a homomorphism ${\cal
F} = {\cal E}/{\cal G}\to {\cal Q}$. The image of this
homomorphism in $\ext_{\cal O}({\cal F},{\cal G})$ only depends
on the isomorphism class of $\cal E$, and the resulting map
$\ext^1_{\cal O}({\cal F},{\cal G})\to\ext_{\cal O}({\cal
F},{\cal G})$ is bijective.

This induces a structure of abelian group on $\ext^1_{\cal
O}({\cal F},{\cal G})$; this structure can be obtained directly
from operations on extensions, as follows.

The identity element corresponds to the split extension $({\cal
F}\oplus {\cal G},\iota,\kappa)$, where $\iota(y) = (0,y)$ and
$\kappa(x,y) = x$. We will denote the split extension by ${\bf
0}_{{\cal F},{\cal G}}$, or simply {\bf 0}.

\proclaim Definition. Let $({\cal E},\iota,\kappa)$ be an
extension of
${\cal F}$ by ${\cal G}$. The {\rm opposite\/} $-{\cal E}$ is
the extension $({\cal E},-\iota,\kappa)$.
\endproclaim

Notice that if $f\colon {\cal E}\to {\cal E}'$ is a homomorphism
of extensions, then the same sheaf homomorphism $f$ is also a
homomorphism from
$-{\cal E}$ to $-{\cal E}'$; we will denote $f$, thought of as a
homomorphism from $-{\cal E}$ to $-{\cal E}'$, by $\ominus f$.
If we assign to each extension
${\cal E}$ the extension $-{\cal E}$, and to each homomorphism
$f\colon {\cal E}\to {\cal E}'$ the homomorphism $\ominus
f\colon -{\cal E}\to -{\cal E}'$, we get a functor from
$\extcat$ to itself, whose square is the
identity\footnote{$^*$}{This is a very rare example of an
involution in a category, whose square is actually the identity,
as opposed to being canonically isomorphic to the identity.}.

Let $({\cal E}_1,\iota_1,\kappa_1),\dots,({\cal
E}_r,\iota_r,\kappa_r)$ be extensions of
${\cal F}$ by ${\cal G}$. Call
$${\cal A}\subseteq {\cal E}_1\oplus\cdots \oplus {\cal E}_r$$
the subsheaf whose sections over an open subset $U\subseteq X$
are of the form
$(e_1,\ldots,e_r)$, where
$e_i\in {\cal E}_i(U)$, $\kappa_1(e_1) = \cdots = \kappa_r(e_r)$.
Clearly ${\cal A}$ is a subsheaf of ${\cal O}$-modules of ${\cal
E}_1\oplus\cdots \oplus {\cal E}_r$. Call ${\cal B}\subseteq
{\cal A}$ the subsheaf of whose sections over
$U\subseteq X$ are of the form $\bigl(\iota_1(y_1),\ldots,
\iota_r(y_r)\bigr)$, where
$y_1,\dots,y_r$ are in ${\cal G}(U)$, and $\sum_{i+1}^r y_i = 0$.
Again ${\cal B}$ is a subsheaf of ${\cal O}$-modules of ${\cal
A}$. If $(e_1,\ldots,e_r)\in {\cal A}(U)$, we will denote the
image of $(e_1,\ldots,e_r)$ in ${\cal A}/{\cal B}$ by
$[e_1,\ldots,e_r]$. Notice that if $y\in{\cal G}(U)$, then
$$[\iota_1(y),0,\ldots,0] = \cdots = [0,\cdots,0,\iota_r(y)]\in
({\cal A}/{\cal B})(U).$$

\proclaim Definition. Let $({\cal
E}_1,\iota_1,\kappa_1),\dots,({\cal E}_r,\iota_r,\kappa_r)$ be
extensions of ${\cal F}$ by ${\cal G}$. The {\rm sum\/}
$$\sum_{i=1}^r {\cal E}_i = {\cal E}_1+\cdots+{\cal E}_r$$ is the
extension $({\cal A}/{\cal B},\iota,\kappa)$ of ${\cal F}$ by
${\cal G}$, where $\iota\colon {\cal G}\to {\cal A}/{\cal B}$ is
defined by $\iota(y) = [\iota_1(y),0,\ldots,0]$, and
$\kappa\colon {\cal A}/{\cal B}\to {\cal F}$ is defined by
$\kappa([e_1,\ldots,e_r]) =
\kappa_1(e_1)$.
\endproclaim

To avoid confusion it may be appropriate to observe that
$\gothext_{\cal O}({\cal F},{\cal G})$ is not an additive
category, and the sum defined here is neither a categorical sum
nor a product.

We leave it to the reader to check that the sum is indeed an
extension of ${\cal F}$ by
${\cal G}$. The sum ${\cal E}_1+(-{\cal E}_2)$ will be denoted
with
${\cal E}_1-{\cal E}_2$.

If $f_1\colon {\cal E}_1\to {\cal E}'_1,\dots,f_r\colon {\cal
E}_r\to {\cal E}'_r$, are homomorphisms of extensions, there is
an induced homomorphism
$$f_1+\cdots+f_r = \sum_{i=1}^r f_i\colon \sum_{i=1}^r {\cal
E}_i\to
\sum_{i=1}^r {\cal E}'_i$$ defined by the obvious rule
$$\biggl(\sum_{i=1}^r f_i\biggr)([e_1,\ldots,e_r]) =
[f_1(e_1),\ldots,f_r(e_r)].$$ This makes the direct sum a
functor from $\extcat^r$ to $\extcat$.

The following properties are straightforward to prove and left
to the reader.

\proclaim Proposition. \call{sumext} \punto\call{sumextass} If
${\cal E}_1,\dots,{\cal E}_r$, ${\cal E}_{r+1},\dots,{\cal
E}_{r+s}$ are extensions of ${\cal F}$ by ${\cal G}$, there is a
functorial isomorphism of extensions
$$\sum_{i=1}^r {\cal E}_i+\sum_{i=r+1}^{r+s} {\cal E}_i \simeq
\sum_{i=1}^{r+s} {\cal E}_i$$ which sends
$\bigl[[e_1,\ldots,e_r],[e_{r+1},\ldots,e_{r+s}]\bigr]$ into
$[e_1,\ldots,e_{r+s}]$.

\punto \call{sumextcomm}If ${\cal E}_1$ and ${\cal E}_2$ are
extensions, there is functorial isomorphism
$$\chi_{{\cal E}_1,{\cal E}_2}\colon {\cal E}_1\oplus {\cal
E}_2\simeq {\cal E}_2\oplus {\cal E}_1$$ which sends $[e_1,e_2]$
into $[e_2,e_1]$.

\punto \call{sumextminus}If ${\cal E}_1$ and ${\cal E}_2$ are
extensions, then
$$ -({\cal E}_1+{\cal E}_2) = (-{\cal E}_1)+(-{\cal E}_2).
$$

\punto \call{sumextid} For each extension ${\cal E}$ there are
functorial isomorphisms
$$\epsilon^{\rm r}_{\cal E}\colon {\cal E}+{\bf 0}_{{\cal
F},{\cal G}}\simeq {\cal E}$$ and
$$\epsilon^{\ell}_{\cal E}\colon {\bf 0}_{{\cal F},{\cal
G}}+{\cal E}\simeq {\cal E}$$ which send $[e,(x,y)]$ into
$e+\iota_{\cal E}(y)$, and $[(x,y),e]$ into $e+\iota_{\cal
E}(y)$, respectively.

\punto \call{sumextinv} If ${\cal E}$ is an extension, then
there is a functorial isomorphism of extensions
$$\delta_{\cal E}\colon {\cal E}-{\cal E}\simeq {\bf 0}_{{\cal
F}, {\cal G}},$$ such that a section $[e_1,e_2]$ of ${\cal
E}-{\cal E}$ is sent to
$(\kappa(e_1),y)$, if $y$ is the section of $\cal G$ such that
$\iota_{\cal E}(y) = e_1-e_2$.

\endproclaim

We will not distinguishing between $\sum_{i=1}^r {\cal
E}_i+\sum_{i=r+1}^{r+s} {\cal E}_i$ and $\sum_{i=1}^{r+s} {\cal
E}_i$, but we will use the isomorphism of
\ref{part sumextass} to identify them.

The functoriality statement in part \ref{part sumextid} should be
interpreted as saying that if
$f\colon {\cal E}\to {\cal E}'$ is a homomorphism of extensions,
then
$$\epsilon^{\rm r}_{{\cal E}',{\bf 0}}\circ (f+{\rm id}_{\bf 0})
= f\circ \epsilon^{\rm r}_{{\cal E},{\bf 0}}\colon {\cal E}+{\bf
0}_{{\cal F},{\cal G}}\longrightarrow {\cal E}',$$ and
analogously for $\epsilon^\ell$. For part
\ref{part sumextinv} it means that
$f\colon {\cal E}\to {\cal E}'$ is a homomorphism of extensions,
then
$$\delta_{\cal E}\circ(f+\ominus f) = \delta_{{\cal E}'}.$$

\proclaim  Definition. Let ${\cal E}$ be an extension of ${\rm
}$ by
${\cal G}$. A\/ {\rm splitting} of ${\cal E}$ is a sheaf
homomorphism $s\colon {\cal E}\to {\cal F}$ such that
$s\iota_{\cal E} = {\rm id}_{\cal F}$. \endproclaim

If we are given a splitting $s\colon {\cal E}\to {\cal F}$, then
the sheaf homomorphism $f_s\colon {\cal E}\to {\bf 0}_{{\cal F},
{\cal G}}$ defined by $f_s(e) = \bigl(s(e),\kappa_{\cal
E}(e)\bigr)$ is an isomorphism of extensions. Conversely, an
isomorphism of extensions ${\cal E}\to {\bf 0}_{{\cal F}, {\cal
G}}$ is of the form
$f_s$ for a unique splitting
$s\colon {\cal E}\to {\cal F}$, so we will identify splittings of
${\cal E}$ and isomorphisms ${\cal E}\to {\bf 0}_{{\cal F},
{\cal G}}$.

Given two splittings $s_1\colon {\cal E}_1\to {\bf 0}_{{\cal F},
{\cal G}}$ and
$s_2\colon {\cal E}_2\to {\bf 0}_{{\cal F}, {\cal G}}$, we can
define their sum
${\cal E}_1+{\cal E}_2\to {\cal F}$ by the formula
$(s_1+s_2)[e_1,e_2] = s_1(e_1)+s_2(e_2)$; it is readily checked
that $s_1+s_2$ is a well-defined splitting of
${\cal E}_1+{\cal E}_2$. In terms of isomorphisms of extensions,
we have that
$$f_{s_1+s_2} =
\epsilon^{\rm r}_{\bf 0}\circ (f_{s_1}+f_{s_2})\colon {\cal
E}_1+{\cal E}_2\to {\bf 0}_{{\cal F},{\cal G}}.$$

\proclaim Proposition. \call{extisom} Let ${\cal E}_1$ and ${\cal
E}_2$ be extensions. There is a canonical bijective
correspondence between splittings of ${\cal E}_1-{\cal E}_2$ and
isomorphisms
${\cal E}_1\simeq {\cal E}_2$. \endproclaim

\proof Given an isomorphism $f\colon {\cal E}_1\to {\cal E}_2$
we can associate to it the splitting
$$\delta_{{\cal E}_2}\circ(f+{\rm id}_{{\cal E}_2})\colon {\cal
E}_1-{\cal E}_2\to {\bf 0}_{{\cal F},{\cal G}}.$$ Conversely, to
each splitting $s\colon {\cal E}_1-{\cal E}_2\to {\bf 0}_{{\cal
F},{\cal G}}$ we can associate
$$\epsilon^\ell_{{\cal E}_2} \circ (s+{\rm id}_{{\cal E}_2})
\circ ({\rm id}_{{\cal E}_1} + \delta_{-{\cal E}_2}) \circ
(\epsilon^{\rm r}_{{\cal E}_1})^{-1}\colon {\cal
E}_1\longrightarrow {\cal E}_2.$$ It is easy to check that these
two operations are inverse to each other.\endproof

The groups $\ext^1_{\cal O}({\cal F},{\cal G})$ are functorial
both in
$\cal F$ and $\cal G$, and this functoriality already exists at
the level of extensions. Let us begin with the functoriality in
$\cal G$.

Let $g\colon {\cal G}\to {\cal G}'$ be a homomorphism of sheaves
of
$\cal O$-modules, and let $\cal E$ be an extension of $\cal F$
by $\cal G$. We define the pushforward $g_*{\cal E}$ of ${\cal
F}$ by ${\cal G}'$ as follows. As a sheaf, $g_*{\cal E}$ is the
direct sum ${\cal G}'\oplus {\cal E}$, divided by the image of
$\cal G$ under the homomorphism ${\cal G}\to {\cal G}'\oplus
{\cal E}$ which sends a local section $y$ to
$\bigl(g(y),-\iota_{\cal E}(y)\bigr)$. The homomorphism
$\iota_{g_*{\cal E}}\colon {\cal G}'\to g_*{\cal E}$ sends $y'$
into the class $[y',0]$, while
$\kappa_{g_*{\cal E}}\colon g_*{\cal E}\to {\cal F}$ sends
$[y',e]$ into $\kappa_{\cal E}(e)$. It is an easy exercise to
show that $g_*{\cal E}$ is an extension of $\cal F$ by ${\cal
G}'$.

Notice that there a homomorphism $\phi\colon {\cal E}\to
g_*{\cal E}$ which sends a local section $e$ into $[0,e]$. This
homomorphism fits into a commutative diagram
$$
\diagram{30}{0&\mapright& {\cal G}& \mapright^{\iota_{\cal
E}}&{\cal E}&\mapright^{\kappa_{\cal E}}&{\cal F}&\mapright&0\cr
&&\mapdown\lft{g}&&\mapdown\lft\phi&&\bivline\cr 0&\mapright&
{\cal G}'& \mapright^{\iota_{g_*{\cal E}}}&g_*{\cal
E}&\mapright^{\kappa_{g_*{\cal E}}}&{\cal F}&\mapright&0\cr}
$$ The extension $g_*{\cal E}$ is the ``only'' extension with
this property, in the following sense.

\proclaim Lemma. \call{homext} Let ${\cal E}'$ an extension of
$\cal F$ by ${\cal G}'$, and assume that there is a commutative
diagram
$$
\diagram{30}{0&\mapright& {\cal G}& \mapright^{\iota_{\cal
E}}&{\cal E}&\mapright^{\kappa_{\cal E}}&{\cal F}&\mapright&0\cr
&&\mapdown\lft g&&\mapdown\lft{\phi'}&&\bivline\cr 0&\mapright&
{\cal G}'& \mapright^{\iota_{{\cal E}'}}&{\cal
E}'&\mapright^{\kappa_{{\cal E}'}}&{\cal F}&\mapright&0.\cr}
$$

Then is a unique homomorphism of extensions $\psi\colon g_*{\cal
E}\simeq {\cal E}'$ such that
$\psi\phi = \phi'$. \endproclaim

\proof It is easy to see that if $\psi$ exists it must have the
form
$\psi([y',e]) =
\iota_{{\cal E}'}(y')+\sigma'(e)$. Conversely one checks that
this formula yields a well defined homomorphism $\psi$ such that
$\psi\sigma = \sigma'$.
\endproof

This construction has the following properties.

\proclaim Proposition. \call{push} Let $g,g_1,g_2\colon {\cal
G}\to {\cal G}'$ be homomorphisms of sheaves of $\cal
O$-modules, $\cal E$, ${\cal E}_1$ and ${\cal E}_2$ extensions
of $\cal F$ by $\cal G$.

\punto \call{pushzero} If $0\colon {\cal G}\to {\cal G}'$ is the
zero homomorphism, then $0_*{\cal E}$ is canonically isomorphic
to ${\bf 0}_{{\cal F},{\cal G}'}$.

\punto \call{pushextopp} $(-g)_*{\cal E} = -(g_*{\cal E})$.

\punto \call{pushexta} There is a canonical isomorphism of
extensions of $\cal F$ by ${\cal G}'$
$$g_*({\cal E}_1+{\cal E}_2)\simeq g_*{\cal E}_1+g_*{\cal E}_2.$$

\punto \call{pushextb} There is a canonical isomorphism of
extensions of $\cal F$ by ${\cal G}'$
$$(g_1+g_2)_*{\cal E}\simeq {g_1}_*{\cal E}+{g_2}_*{\cal E}.$$

\punto \call{pushextbound} The boundary homomorphism
$$\partial \colon \hom_{\cal O}({\cal G},{\cal
G}')\longrightarrow
\ext_{\cal O}^1({\cal F},{\cal G}')$$ coming from the sequence
$$\sequence{0&\mapright&{\cal G}&\mapright &{\cal E}& \mapright
&{\cal F}& \mapright&0}$$ sends $g\in \hom_{\cal O}({\cal
G},{\cal G}')$ into the isomorphism class of $g_*{\cal E}$.
\endproclaim

\proof Part \ref{part pushextopp} is straightforward.

For part \ref{part pushexta} we apply
\ref{homext} with
${\cal E}' = g_*{\cal E}_1+g_*{\cal E}_2$, and $\phi'\colon {\cal
E}_1+{\cal E}_2\to {\cal E}'$ defined by
$\phi'([e_1,e_2]) = \bigl[[0,e_1],[0,e_2]\bigr]$.

For part \ref{part pushextb} we take ${\cal E}' = {g_1}_*{\cal
E}+{g_2}_*{\cal E}$,
$\phi'([e_1,e_2]) = \bigl[[0,e_1],[0,e_2]\bigr]$.

Part \ref{part pushextbound} is standard.
\endproof

If $s\colon g_*{\cal E}\to {\cal G}'$ of $g_*{\cal E}$ is a
splitting of $g_*{\cal E}$, then the homomorphism $\sigma =
\phi\circ s$ has the property that
$\sigma\circ \iota_{\cal E} = g$. Conversely, given a
homomorphism $\sigma = \phi\circ s$ such that $\sigma\circ
\iota_{\cal E} = g$, we get a splitting $s\colon g_*{\cal E}\to
{\cal G}'$ by defining $s([y',e]) = y'+\sigma(e)$.

\proclaim Lemma. \call{pushsplit} There is a natural bijective
correspondence of splittings of $g_*{\cal E}$ with homomorphisms
$\sigma = \phi\circ s$ such that
$\sigma\circ
\iota_{\cal E} = g$.
\endproclaim

The functoriality in $\cal F$ is analogous. Let $\cal E$ be an
extension of $\cal F$ by $\cal G$, $f\colon {\cal F}'\to {\cal
F}$ be a homomorphism. We define the pullback $f^*{\cal E}$ of
${\cal F}'$ by $\cal G$ as the subsheaf of ${\cal F}'\oplus
{\cal E}$ whose sections are of type $(x',e)$ with $f(x') =
\iota_{\cal E}(e)$; it is a subsheaf of $\cal O$-modules of
${\cal F}'\oplus {\cal E}$. The homomorphism
$\iota_{f^*{\cal E}}$ from ${\cal G}$ to $f^*{\cal E}$ sends $y$
to
$(0,\iota_{\cal E}(y))$, and $\kappa_{f^*{\cal E}}$ from
$f^*{\cal E}$ to
${\cal F}'$ sends $(x',e)$ to $x'$. We leave it to the reader to
check that $f^*{\cal E}$ is indeed an extension of ${\cal F}'$
by $\cal G$.

The homomorphism $\phi\colon f^*{\cal E}\to {\cal E}$ which sends
$(x',e)$ into $e$ fits into a commutative diagram
$$
\diagram{30}{0&\mapright&{\cal G}&\mapright^{\iota_{f^*{\cal
E}}}& f^*{\cal E}&\mapright^{\kappa_{f^*{\cal E}}}&{\cal F}'&
\mapright& 0\cr
&&\bivline&&\mapdown\lft{\phi}&&\mapdown\lft{f}\cr
0&\mapright&{\cal G}&\mapright^{\iota_{\cal E}}&{\cal
E}&\mapright^{\kappa_{\cal E}}&{\cal F}&
\mapright& 0.\cr}
$$

The following results are dual to the results stated for
pushforwards. The proofs are left to the reader.

\proclaim Lemma. \call{homextpull} Let ${\cal E}'$ an extension
of
${\cal F}'$ by $\cal G$, and assume that there is a commutative
diagram
$$
\diagram{30}{0&\mapright& {\cal G}& \mapright^{\iota_{{\cal
E}'}}&{\cal E}'&\mapright^{\kappa_{{\cal E}'}}&{\cal
F}'&\mapright&0\cr &&\bivline&&\mapdown\lft{\phi'}&&\mapdown\lft
f\cr 0&\mapright& {\cal G}& \mapright^{\iota_{\cal E}}&{\cal
E}&\mapright^{\kappa_{\cal E}}&{\cal F}&\mapright&0.\cr}
$$

Then is a unique homomorphism of extensions $\psi\colon {\cal
E}'\simeq f^*{\cal E}'$ such that
$\phi\psi = \phi'$. \endproclaim

This construction has the following properties.

\proclaim Proposition. \call{pull} Let $f,f_1,f_2\colon {\cal
G}\to {\cal G}'$ be homomorphisms of sheaves of $\cal
O$-modules, $\cal E$, ${\cal E}_1$ and ${\cal E}_2$ extensions
of $\cal F$ by $\cal G$.

\punto \call{pullzero} If $0\colon {\cal F}'\to {\cal F}$ is the
zero homomorphism, then $0^*{\cal E}$ is canonically isomorphic
to ${\bf 0}_{{\cal F}',{\cal G}}$.

\punto \call{pullextopp} $(-f)^*{\cal E} = -(f^*{\cal E})$.

\punto \call{pullexta} There is a canonical isomorphism of
extensions of ${\cal F}'$ by ${\cal G}$
$$f^*({\cal E}_1+{\cal E}_2)\simeq f^*{\cal E}_1+f^*{\cal E}_2.$$

\punto \call{pullextb} There is a canonical isomorphism of
extensions of ${\cal F}'$ by ${\cal G}$
$$(f_1+f_2)^*{\cal E}\simeq {f_1}^*{\cal E}+{f_2}^*{\cal E}.$$
\endproclaim

The following concept arises naturally in our construction of the
obstruction in
\ref{absdefb}.

\proclaim Definition. \call{extcocycle} Let ${\cal U} =
\{X_\alpha\}$ be an open covering  of $X$. An {\rm extension
cocycle\/}
$$(\{{\cal E}_{\alpha\beta}\},\{F_{\alpha\beta\gamma}\})$$
 of $\cal F$ by
$\cal G$ on $\cal U$ is a collection of extensions $\{{\cal
E}_{\alpha\beta}\}$ of
${\cal F}\rest{X_\alpha\cap X_\beta}$ by ${\cal
G}\rest{X_\alpha\cap X_\beta}$, and isomorphisms
$F_{\alpha\beta\gamma}\colon {\cal E}_{\alpha\beta}+{\cal
E}_{\beta\gamma}\simeq {\cal E}_{\alpha\gamma}$ on $X_\alpha\cap
X_\beta\cap X_\gamma$, such that for any quadruple
$\alpha$, $\beta$, $\gamma$, $\delta$ we have
$$ F_{\alpha\gamma\delta}\circ\bigl(F_{\alpha\beta\gamma}+{\rm
id}_{{\cal E}_{\gamma\delta}}\bigr) = 
F_{\alpha\beta\delta}\circ\bigl({\rm id}_{{\cal
E}_{\alpha\beta}}+F_{\beta\gamma\delta}\bigr)
$$ as a homomorphism ${\cal E}_{\alpha\beta}+{\cal
E}_{\beta\gamma}+{\cal E}_{\gamma\delta}\longrightarrow {\cal
E}_{\alpha\delta}.$
\endproclaim

There is an obvious notion of isomorphism of extension cocycles.
An isomorphism
$$
\phi\colon (\{{\cal
E}_{\alpha\beta}\},\{F_{\alpha\beta\gamma}\})\simeq (\{{\cal
E}'_{\alpha\beta}\},\{F'_{\alpha\beta\gamma}\})
$$ consists of a collection of isomorphisms of extensions
$\phi_{\alpha\beta}\colon {\cal E}_{\alpha\beta}\simeq {\cal
E}'_{\alpha\beta}$ such that
$$\phi_{\alpha\gamma}\circ F_{\alpha\beta\gamma} =
F'_{\alpha\beta\gamma}\circ
(\phi_{\alpha\beta}+\phi_{\beta\phi})$$ for all $\alpha$,
$\beta$ and $\gamma$.

Furthermore, extension cocycles can be summed; we define
$$ (\{{\cal E}_{\alpha\beta}\},\{F_{\alpha\beta\gamma}\})+
(\{{\cal E}'_{\alpha\beta}\},\{F'_{\alpha\beta\gamma}\}) =
(\{{\cal E}_{\alpha\beta}+{\cal E}'_{\alpha\beta}\},
\{\widetilde F_{\alpha\beta\gamma}\})
$$ where
$$
\widetilde F_{\alpha\beta\gamma}\colon {\cal
E}_{\alpha\beta}+{\cal E}'_{\alpha\beta}- ({\cal
E}_{\beta\gamma}+ {\cal E}'_{\beta\gamma}) = {\cal
E}_{\alpha\beta}+{\cal E}'_{\alpha\beta}- {\cal E}_{\beta\gamma}-
{\cal E}'_{\beta\gamma}
$$ is
$$
\widetilde F_{\alpha\beta\gamma} = (F_{\alpha\beta\gamma}+
F'_{\alpha\beta\gamma})\circ ({\rm id}_{{\cal
E}_{\alpha\beta}}+\chi_{{\cal E}'_{\alpha\beta},-{\cal
E}_{\beta\gamma}}+{\rm id}_{{\cal E}_{\beta\gamma}})
$$ and leave it to the reader to check that this sum is still an
extension cocycle.

Consider the set of isomorphism classes of extension cocycles.
The operation of sum introduced above makes it into an abelian
group.

The zero is represented by the class of the trivial extension
cocycle, that is, the cocycle $(\{{\bf 0}\},\{\epsilon^{\rm
r}_{\bf 0}\})$ in which all the
${\cal E}_{\alpha\beta}$ are trivial, and the
$F_{\alpha\beta\gamma}$ are all
$\epsilon^{\rm r}_{\bf 0}\colon {\bf 0}+{\bf 0}\to {\bf 0}$. The
isomorphism
$$(\{{\bf 0}\},\{\epsilon^{\rm r}_{\bf 0}\})+(\{{\cal
E}_{\alpha\beta}\},\{F_{\alpha\beta\gamma}\})\simeq (\{{\cal
E}_{\alpha\beta}\},\{F_{\alpha\beta\gamma}\})$$ is given by
$$\epsilon^\ell_{{\cal E}_{\alpha\beta}}\colon {\bf 0}+{\cal
E}_{\alpha\beta}\simeq {\cal E}_{\alpha\beta}.$$

The inverse of the class of an extension cocycle
$(\{{\cal E}_{\alpha\beta}\},
\{F_{\alpha\beta\gamma}\})$ is $(\{-{\cal E}_{\alpha\beta}\},
\{\ominus F_{\alpha\beta\gamma}\})$, and the isomorphism
$$(\{{\cal E}_{\alpha\beta}\},
\{F_{\alpha\beta\gamma}\})+(\{-{\cal E}_{\alpha\beta}\},
\{\ominus F_{\alpha\beta\gamma}\})\simeq (\{{\bf
0}\},\{\epsilon^{\rm r}_{\bf 0}\})$$ is given by
$$\delta_{{\cal E}_{\alpha\beta}}\colon {\cal
E}_{\alpha\beta}-{\cal E}_{\alpha\beta}\simeq {\bf 0}.$$

Given a collection of extensions $\{{\cal E}_\alpha\}$ of ${\cal
F}\rest{X_\alpha}$ by
${\cal G}\rest{X_\alpha}$ on $X_\alpha$, we can define its {\it
boundary}
$$
\partial \{{\cal E}_\alpha\} = (\{{\cal E}_\alpha-{\cal
E}_\beta\}, F_{\alpha\beta\gamma})
$$ where
$$ F_{\alpha\beta\gamma} = (\epsilon^{\rm r}_{{\cal
E}_\alpha}+{\rm id}_{-{\cal E}_\gamma})
\circ ({\rm id}_{{\cal E}_\alpha}+\delta_{-{\cal
E}_\beta}+\delta_{-{\cal E}_\gamma})\colon {\cal E}_\alpha -
{\cal E}_\beta + {\cal E}_\beta -{\cal E}_\gamma\longrightarrow
{\cal E}_\alpha -{\cal E}\gamma.
$$

One check that $\partial \{{\cal E}_\alpha\}$ is an obstruction
cocycle, and that there is a canonical isomorphism of extension
cocycles
$$
\partial \{{\cal E}_\alpha+{\cal E}'_\alpha\} \simeq \partial
\{{\cal E}_\alpha\} + \partial
\{{\cal E}_\alpha\}
$$ given by
$$ {\rm id}_{{\cal E}_\alpha}+\chi_{{\cal E}'_\alpha,-{\cal
E}_\beta}+{\rm id}_{{\cal E}'_\beta}
\colon {\cal E}_\alpha+{\cal E}'_\alpha-({\cal E}_\beta + {\cal
E}'_\beta) =  {\cal E}_\alpha+{\cal E}'_\alpha-{\cal E}_\beta -
{\cal E}'_\beta\longrightarrow {\cal E}_\alpha- {\cal E}_\beta
+{\cal E}'_\alpha- {\cal E}'_\beta.
$$

We say that an extension cocycle is a {\it boundary\/} if it is
isomorphic to the boundary of a collection of extensions
$\{{\cal E}_\alpha\}$.

The isomorphism classes of boundaries form a subgroup of the
group of isomorphism classes of extension cocycles; the quotient
group will be called the {\it group of extension classes}, and
will be denoted by $\E_{\cal O}({\cal U};{\cal F},{\cal G})$.

Let ${\cal U}' = \{X'_{\alpha'}\}$ be refinement of\/ $\cal U$;
there is a function
$\rho\colon {\cal U}'\to {\cal U}$ such that
$X'_{\alpha'}\subseteq X_{\rho(\alpha')}$. An extension cocycle
$(\{{\cal E}_{\alpha\beta}\},
\{F_{\alpha\beta\gamma}\})$ on $\cal U$ induces an extension
cocycle
$(\{{\cal E}'_{\alpha'\beta'}\},\{F'_{\alpha'\beta'\gamma'}\})$,
where ${\cal E}'_{\alpha',\beta'}$ is the restriction of the
extension ${\cal E}_{\rho(\alpha')\rho(\beta')}$ to
$X'_{\alpha'}\cap X_\beta'$, and the $F_{\alpha\beta\gamma}$ are
also restricted. Boundaries are obviously brought to boundaries,
so we get a restriction map
$\E_{\cal O}({\cal U};{\cal F},{\cal G})\to \E_{\cal O}({\cal
U}';{\cal F},{\cal G})$; as we shall see during the proof of
\ref{extcocgrp}, this map only depends on $\cal U$ and
${\cal U}'$, and not on the map
$\rho$. Accepting this for the moment, there is a limit abelian
group
$$\E_{\cal O}({\cal F},{\cal G}) = \indlim_{\cal U}\E_{\cal
O}({\cal U};{\cal F},{\cal G}).$$

\proclaim Theorem. \call{extcocgrp} \punto There is canonical
group isomorphism of
$\E_{\cal O}({\cal F},{\cal G})$ with the kernel of the
localization map $\ext^2_{\cal O}(F{\cal },{\cal G})\to
\H^0\bigl(X,\cursext^2_{\cal O}(F{\cal },{\cal G})\bigr)$.

\punto The natural map $\E_{\cal O}({\cal U};{\cal F},{\cal
G})\to
\E_{\cal O}({\cal F},{\cal G})$ is injective. \endproclaim

The group $\E_{\cal O}({\cal F},{\cal G})$ can be interpreted as
the first cohomology group of the Picard stack
${\goth Ext}_{\cal O}({\cal F},{\cal G})$ of extensions of $\cal
F$ and $\cal G$ (see [Deligne]); the theorem does not seem to
appear in this form in the literature, but, as it was pointed
out to me by L. Breen, it was more or less known: see for
example [Retakh] and [Ulbrich].

\proof We will use the following notation. Let
${\cal E}$, ${\cal E}_1$, ${\cal E}_2$ be extensions of $\cal F$
by $\cal G$, and let $j\colon {\cal G}\to {\cal J}$ be a
homomorphism of sheaves.

If $\sigma_1\colon {\cal E}_1\to {\cal J}$ and $\sigma_2\colon
{\cal E}_2\to {\cal J}$ are homomorphisms with $\sigma_i\circ
\iota_{{\cal E}_i}= j$, we will take their sum
$\sigma_1+\sigma_2\colon {\cal E}_1+{\cal E}_2\to {\cal J}$ to
be the homomorphism defined by the formula
$(\sigma_1+\sigma_2)([e_1,e_2]) =
\sigma_1(e_1)+
\sigma_2(e_2)$. If we think of $\sigma_i$ as a splitting of
$j_*({\cal E}_i)$ (\ref{pushsplit}), then $\sigma_1+\sigma_2$
can be thought of as their sum as a splitting  of
$j_*({\cal E}_1+{\cal E}_2) = j_*({\cal E}_1)+j_*({\cal E}_2)$
(\ref{pushexta}).

On the other hand, let $\sigma_1,\sigma_2\colon {\cal E}\to
{\cal J}$ be homomorphisms with $\sigma_i\circ \iota_{\cal E}=
j$. Their difference
$\sigma_1-\sigma_2\colon {\cal E}\to {\cal J}$ is a homomorphism
with $(\sigma_1-\sigma_2)\circ \iota_{\cal E} = 0$, so there is
a unique homomorphism $\tau\colon {\cal E}\to {\cal J}$ such that
$\tau\circ \kappa_{\cal E} = \sigma_1-\sigma_2$. This $\tau$ we
will also denote by
$\sigma_1-\sigma_2$.

Finally, if $\sigma\colon {\cal E}\to {\cal J}$ is such that
$\sigma\circ\iota_{\cal E} = j$, and
$\tau\colon {\cal F}\to {\cal J}$ is a homomorphism, we will
write
$\sigma+\tau$ to mean $\sigma+\tau\circ \kappa_{\cal E}$.

Let $\cal J$ be an injective sheaf of $\cal O$-modules
containing $\cal G$, ${\cal Q} = {\cal J}/{\cal G}$. Call
$j\colon {\cal G}\to {\cal J}$ the inclusion, $\pi\colon {\cal
J}\to {\cal Q}$ the projection. Then the boundary operator
$$\partial \colon \ext^1_{\cal O}({\cal F},{\cal Q})\to
\ext^2_{\cal O}({\cal F},{\cal G})$$ is an isomorphism, and
induces an isomorphism of $\cursext^1_{\cal O}({\cal F},{\cal
Q})$ with $\cursext^2_{\cal O}({\cal F},{\cal G})$, such the
diagram
$$
\diagram{40}{\ext^1_{\cal O}({\cal F},{\cal
Q})&\mapright^{\displaystyle\sim}&\ext^2_{\cal O}({\cal F},{\cal
G})\cr
\mapdown&&\mapdown\cr
\H^0\bigl(X,\cursext^1_{\cal O}({\cal F},{\cal
Q})\bigr)&\mapright^{\displaystyle\sim}&\H^0\bigl(X,\cursext^2_{\cal
O}({\cal F},{\cal G})\bigr)\cr},
$$ where the rows are boundary maps and the columns are
localization maps, commutes. Hence the kernel of the left column
is isomorphic to the kernel of the right column. But from the
spectral sequence
$$ E^{pq}_2 = \H^p\bigl(X,\cursext^q_{\cal O}({\cal F},{\cal
G})\bigr)\Rightarrow
\ext^{p+q}_{\cal O}({\cal F},{\cal G})
$$ we get an exact sequence
$$\sequence{0&\mapright&\H^1\bigl(X,\curshom_{\cal O}({\cal
F},{\cal G})\bigr)& \mapright & \ext^1_{\cal O}({\cal F},{\cal
Q})&\mapright&\H^0\bigl(X,\cursext^1_{\cal O}({\cal F},{\cal
Q})\bigr),\cr}$$ and an isomorphism of
$\H^1\bigl(X,\curshom_{\cal O}({\cal F},{\cal G})\bigr)$ with the
kernel of the localization map from $\ext^2_{\cal O}(F{\cal
},{\cal G})$ to $\H^0\bigl(X,\cursext^2_{\cal O}(F{\cal },{\cal
G})\bigr)$. We will prove both statements, and the fact that the
restriction map from $\E_{\cal O}({\cal U};{\cal F},{\cal G})$
to $\E_{\cal O}({\cal U}';{\cal F},{\cal G})$ only depends on
${\cal U}$ and ${\cal U}'$, and not on the function ${\cal
U}'\to {\cal U}$, by proving the following.

\proclaim Lemma. \call{subextcoc}There is a canonical isomorphism
$$\E_{\cal O}({\cal U};{\cal F},{\cal G})\simeq
\cechH^1\bigl({\cal U},\curshom_{\cal O}({\cal F},{\cal
Q})\bigr)$$ which is compatible with restriction maps.
\endproclaim

If we remember that the restriction map
$$\cechH^1\bigl({\cal U},\curshom_{\cal O}({\cal F},{\cal
Q})\bigr)\longrightarrow
\cechH^1\bigl({\cal U}',\curshom_{\cal O}({\cal F},{\cal
Q})\bigr)$$ independent of the function ${\cal U}'\to {\cal U}$,
that the inductive limit of the
\v Cech cohomology groups is the ordinary cohomology, and that
the map
$$\cechH^1\bigl({\cal U},\curshom_{\cal O}({\cal F},{\cal
Q})\bigr)\longrightarrow
\H^1\bigl(X,\curshom_{\cal O}({\cal F},{\cal Q})\bigr)$$ is
always injective, we attain the proof of the theorem, together
with the statement about the restriction map from $\E_{\cal
O}({\cal U};{\cal F},{\cal G})$ to $\E_{\cal O}({\cal U}';{\cal
F},{\cal G})$ only depending on ${\cal U}$ and ${\cal U}'$.

\proofof subextcoc. Call $j\colon {\cal G}\into {\cal J}$ the
inclusion. Because
$\cal J$ is injective, for each
$\alpha$ and
$\beta$ we can find a homomorphism
$\sigma_{\alpha\beta}\colon {\cal E}_{\alpha\beta}\to {\cal J}$
such that
$\sigma_{\alpha\beta}\circ \iota_{{\cal E}_{\alpha\beta}} = j$.
Let us check that we can do it coherently, in the following
sense.

\proclaim Lemma. \call{mainextcoc} We can find homomorphisms
$\sigma_{\alpha\beta}\colon {\cal E}_{\alpha\beta}\to {\cal J}$
for each
$\alpha$, $\beta$, such that
$\sigma_{\alpha\beta}\circ \iota_{{\cal E}_{\alpha\beta}} = j$,
and such that
$$
\sigma_{\alpha\beta}+\sigma_{\beta\gamma} =
\sigma_{\alpha\gamma}\circ F_{\alpha\beta\gamma}\colon {\cal
E}_{\alpha\beta}+{\cal E}_{\beta\gamma}\longrightarrow {\cal J}
$$ for each triple $\alpha$, $\beta$ and $\gamma$.\endproclaim

We will call such a collection $\{\sigma_{\alpha\beta}\}$ a {\it
function\/} from the extension cocycle $(\{{\cal
E}_{\alpha\beta}\},\{F_{\alpha\beta\gamma}\})$ to $\cal J$.

\proof Choose homomorphisms $\sigma_{\alpha\beta}\colon {\cal
E}_{\alpha\beta}\to {\cal J}$ in such a way that
$\sigma_{\alpha\beta}\circ \iota_{{\cal E}_{\alpha\beta}} = j$
for all $\alpha$ and $\beta$, and consider the homomorphisms

$$
\tau_{\alpha\beta\gamma} =
(\sigma_{\alpha\beta}+\sigma_{\beta\gamma})-\sigma_{\alpha\gamma}\circ
F_{\alpha\beta\gamma}\colon {\cal F}\longrightarrow {\cal J}.
$$

\proclaim Lemma. \call{canfindcoc}The collection
$\{\tau_{\alpha\beta\gamma}\}$ is a \v Cech 2-cocycle in the
sheaf $\curshom_{\cal O}({\cal F},{\cal J})$, that is,
$$\tau_{\alpha\beta\gamma}+\tau_{\alpha\gamma\delta} =
\tau_{\alpha\beta\delta}+\tau_{\beta\gamma\delta}$$ for all
$\alpha$, $\beta$, $\gamma$ and $\delta$. \endproclaim

\proof Let $x$ be a local section of $\cal F$, and choose
sections
$e_{\alpha\beta}$,
$e_{\beta\gamma}$ and $e_{\gamma\delta}$ of ${\cal
E}_{\alpha\beta}$,
${\cal E}_{\beta\gamma}$ and ${\cal E}_{\gamma\delta}$ whose
image in $\cal F$ is $x$. Set
$e_{\alpha\gamma} =
F_{\alpha\beta\gamma}([e_{\alpha\beta},e_{\beta\gamma}])$,
$e_{\beta\delta} =
F_{\beta\gamma\delta}([e_{\beta\gamma},e_{\gamma\delta}])$,
$e_{\alpha\delta} =
F_{\alpha\beta\delta}([e_{\alpha\beta},e_{\beta\delta}])$.
Observe that the cocycle condition on the $F_{\alpha\beta\gamma}$
(\ref{extcocycle}) is tailor made to give us
$$ e_{\alpha\delta} =
F_{\alpha\beta\delta}([e_{\alpha\beta},e_{\beta\delta}]) =
F_{\alpha\gamma\delta}([e_{\alpha\gamma},e_{\gamma\delta}]).
$$ Then
$$
\eqalign{(\tau_{\alpha\beta\gamma}+\tau_{\alpha\gamma\delta}-
\tau_{\alpha\beta\delta}- \tau_{\beta\gamma\delta})(x)&=
\sigma_{\alpha\beta}(e_{\alpha\beta}) +
\sigma_{\beta\gamma}(e_{\beta\gamma}) -
\sigma_{\alpha\gamma}\bigl(F_{\alpha\beta\gamma}([e_{\alpha\beta},
e_{\beta\gamma}])\bigr)
\cr &\qquad +\sigma_{\alpha\gamma}(e_{\alpha\gamma}) +
\sigma_{\gamma\delta}(e_{\gamma\delta}) -
\sigma_{\alpha\delta}\bigl(F_{\alpha\gamma\delta}
([e_{\alpha\gamma},e_{\gamma\delta}])\bigr)\cr &\qquad
-\sigma_{\alpha\beta}(e_{\alpha\beta}) -
\sigma_{\beta\delta}(e_{\beta\delta}) +
\sigma_{\alpha\delta}\bigl(F_{\alpha\beta\delta}
([e_{\alpha\beta},e_{\beta\delta}])\bigr)\cr &\qquad
-\sigma_{\beta\gamma}(e_{\beta\gamma}) -
\sigma_{\gamma\delta}(e_{\gamma\delta}) +
\sigma_{\beta\delta}\bigl(F_{\beta\gamma\delta}
([e_{\beta\gamma},e_{\gamma\delta}])\bigr)\cr & =
\sigma_{\alpha\beta}(e_{\alpha\beta}) +
\sigma_{\beta\gamma}(e_{\beta\gamma}) -
\sigma_{\alpha\gamma}(e_{\alpha\gamma})\cr &\qquad
+\sigma_{\alpha\gamma}(e_{\alpha\gamma}) +
\sigma_{\gamma\delta}(e_{\gamma\delta}) -
\sigma_{\alpha\delta}(e_{\alpha\delta})\cr &\qquad
-\sigma_{\alpha\beta}(e_{\alpha\beta}) -
\sigma_{\beta\delta}(e_{\beta\delta}) +
\sigma_{\alpha\delta}(e_{\alpha\delta})\cr &\qquad
-\sigma_{\beta\gamma}(e_{\beta\gamma}) -
\sigma_{\gamma\delta}(e_{\gamma\delta}) +
\sigma_{\beta\delta}(e_{\beta\delta})\bigr)\cr &=0.\cr}
$$

This proves the lemma.
\endproof

Now observe that the sheaf $\curshom_{\cal O}({\cal F},{\cal
J})$ is flabby, because $\cal J$ is injective, so its second \v
Cech cohomology group
$\cechH^2\bigl({\cal U},\curshom_{\cal O}({\cal F},{\cal
J})\bigr)$ is 0. Hence we can find a 1-cochain
$\{\tau_{\alpha\beta}\}$ of
$\curshom_{\cal O}({\cal F},{\cal J})$ such that
$\tau_{\alpha\beta\gamma} =
\tau_{\alpha\beta}+\tau_{\beta\gamma}-\tau_{\alpha\gamma}$. If
we set
$\widetilde
\sigma_{\alpha\beta} = \sigma_{\alpha\beta}+\tau_{\alpha\beta}$
we see easily that the condition
$$\widetilde \sigma_{\alpha\beta}+\widetilde
\sigma_{\beta\gamma} =
\widetilde \sigma_{\alpha\gamma}\circ F_{\alpha\beta\gamma}\colon
{\cal E}_{\alpha\beta}+ {\cal E}_{\beta\gamma}\longrightarrow
{\cal J}$$ is satisfied. This proves \ref{canfindcoc}.
\endproof

Choose a function $\{\sigma_{\alpha\beta}\}$ from $(\{{\cal
E}_{\alpha\beta}\},\{F_{\alpha\beta\gamma}\})$ to $\cal J$. The
composition of the
$\sigma_{\alpha\beta}$ with the projection
$\pi\colon {\cal J}
\to {\cal Q}$ send
$\cal G$ to 0, and therefore induce homomorphisms
$\eta_{\alpha\beta}\colon {\cal F}\to{\cal Q}$ satisfying the
cocycle condition
$\eta_{\alpha\beta}+\eta_{\beta\gamma} = \eta_{\alpha\gamma}$.
So we have associated to the extension cocycle $(\{{\cal
E}_{\alpha\beta}\},\{F_{\alpha\beta\gamma}\})$ and the functions
$\{\sigma_{\alpha\beta}\}$ an element $\omega$ of the \v Cech
cohomology group
$\cechH^1\bigl({\cal U},\curshom_{\cal O}({\cal F},{\cal
Q})\bigr)$. Let us check that this element does not depend on
the function $\{\sigma_{\alpha\beta}\}$. Let
$\{\sigma'_{\alpha\beta}\}$ be another function, and call
$\omega'$ the element of $\cechH^1\bigl({\cal U},\curshom_{\cal
O}({\cal F},{\cal Q})\bigr)$ associated with
$(\{{\cal E}_{\alpha\beta}\},\{F_{\alpha\beta\gamma}\})$ and
$\{\sigma'_{\alpha\beta}\}$. Consider the cocycle
$\{\sigma_{\alpha\beta}-\sigma'_{\alpha\beta}\}$ in
$\curshom_{\cal O}({\cal F},{\cal J})$. The compositions
$\pi\circ(\sigma_{\alpha\beta}-\sigma'_{\alpha\beta})$ define a
cocycle in $\curshom_{\cal O}({\cal F},{\cal Q})$ whose
cohomology class is clearly
$\omega-\omega'$. But
$\cechH^1\bigl({\cal U},\curshom_{\cal O}({\cal F},{\cal
J})\bigr)$ is 0, because
$\curshom_{\cal O}({\cal F},{\cal J})$ is flabby, and so
$\omega-\omega' = 0$.

This construction gives a mapping from $\E_{\cal O}(U{\cal
};{\cal F},{\cal G})$ to $\cechH^1\bigl({\cal U},\curshom_{\cal
O}({\cal F},{\cal Q})\bigr)$. We leave to the reader to check
that it is a homomorphism. Let us prove that it is bijective.

Let $(\{{\cal E}_{\alpha\beta}\},\{F_{\alpha\beta\gamma}\})$ be
an extension cocycle whose associated cohomology element is 0.
Choose a function
$\{\sigma_{\alpha\beta}\}$ from
$(\{{\cal E}_{\alpha\beta}\},\{F_{\alpha\beta\gamma}\})$ to
$\cal J$, and call
$\{\eta_{\alpha\beta}\}$ the associated cocycle in
$\curshom_{\cal O}({\cal F},{\cal Q})$. Choose a 0-cochain
$\{\eta_\alpha\}$ such that $\eta_{\alpha\beta} =
\eta_\beta-\eta_\alpha$ for all $\alpha$ and $\beta$. We think
of $\cal J$ as an extension of
$\cal G$ by $\cal Q$, and set ${\cal E}_\alpha =
\eta_\alpha^*{\cal J}$. Observe that by definition of
$\eta_{\alpha\beta}$ the diagram
$$
\diagram{30}{0&\mapright &{\cal G}&\mapright &{\cal
E}_{\alpha\beta} &\mapright & {\cal F}&\mapright &0\cr
&&\bivline&&\mapdown\lft{\sigma_{\alpha\beta}}
&&\mapdown\lft{\eta_{\alpha\beta}}\cr 0&\mapright&{\cal G}&
\mapright& {\cal J}&\mapright& {\cal Q}&\mapright& 0\cr}
$$ commutes, so, by \ref{homextpull}, there is an induced
isomorphism of
${\cal E}_{\alpha\beta}$ with $\eta_{\alpha\beta}^*{\cal J} =
(\eta_\beta-\eta_\alpha)^*{\cal J}$. But by \ref{pullextopp} and
\ref{part pullextb} there is a canonical isomorphism of
$(\eta_\beta-\eta_\alpha)^*{\cal J}$ with $\eta_\beta^*{\cal
J}-\eta_\alpha^*{\cal J} = {\cal E}_\beta-{\cal E}_\alpha$. We
leave it to the reader to show that this collection of
isomorphisms is an isomorphism of the extension cocycle $(\{{\cal
E}_{\alpha\beta}\},\{F_{\alpha\beta\gamma}\})$ with the boundary
$\partial \{{\cal E}_\alpha\}$.

To prove surjectivity, let $\omega\in \cechH^1\bigl({\cal
U};{\cal F},{\cal Q}\bigr)$ be a class represented by a cocycle
$\{\eta_{\alpha\beta}\}$, and set
${\cal E}_{\alpha\beta} =
\eta_{\alpha\beta}^*{\cal J}$. The isomorphism
$$F_{\alpha\beta\gamma}\colon {\cal E}_{\alpha\beta}+{\cal
E}_{\beta\gamma} =
\eta_{\alpha\beta}^*{\cal J} + \eta_{\beta\gamma}^*{\cal J}\simeq
(\eta_{\alpha\beta}+\eta_{\beta\gamma})^*{\cal J} =
\eta_{\alpha\gamma}^*{\cal J} = {\cal E}_{\alpha\gamma}$$ is the
inverse of the one in \ref{pullextb}. One checks that the
$F_{\alpha\beta\gamma}$ satisfy the cocycle condition of
\ref{extcocycle}.

It is easy to see that the isomorphism between $\E_{\cal O}({\cal
U};{\cal F},{\cal G})$ and $\cechH^1\bigl({\cal
U},\curshom_{\cal O}({\cal F},{\cal Q})\bigr)$ that we have just
defined is compatible with refinement. This completes the proof
of \ref{mainextcoc}.

To finish the proof of \ref{extcocgrp} we only need to check 
that the resulting isomorphism between $\E_{\cal O}({\cal
F},{\cal G})$ and the kernel of the localization map
$\ext^2_{\cal O}(F{\cal },{\cal G})\to
\H^0\bigl(X,\cursext^2_{\cal O}(F{\cal },{\cal G})\bigr)$ is
independent of the choice of $\cal J$. This is straightforward
and left to the reader.
\endproof

In order to extend \ref{absdef} to maps, we need to generalize
what has been done above to complexes. There is no problem in
defining extensions of complexes, but in general extensions do
not represent enough classes in $\ext^1$. This problem, however,
will not arise in the case we are interested in, as we will see.

If ${\cal F}^\mini$ and ${\cal G}^\mini$ are complexes of
sheaves, we define an extension $({\cal E}^\mini,\iota,\kappa)$
of ${\cal F}^\mini$ by
${\cal G}^\mini$ is a complex of sheaves ${\cal E}^\mini$,
together with homomorphisms of complexes $\iota\colon {\cal
G}^\mini\to {\cal E}^\mini$ and $\kappa\colon {\cal E}^\mini\to
{\cal F}^\mini$ such that the sequence
$$\sequence{0&\mapright&{\cal
G}^n&\mapright^{\textstyle\iota^n}&{\cal E}^n &
\mapright^{\textstyle \kappa^n}&{\cal F}^n&\mapright&0}$$ is
exact for all integers $n$.

All of the theory that we have developed above before
\ref{extcocgrp} goes through without changes to this more
general case, when we substitute everywhere complexes to
sheaves, and homomorphisms of complexes to homomorphism of
complexes.

Let $\cal G$ be a sheaf, ${\cal F}^\mini$ a complex of sheaves,
${\cal J}^\mini$ an injective resolution of $\cal G$. Recall
that with ${\cal J}^\mini[n]$ we denote the complex with ${\cal
J}^\mini[n]^i = {\cal J}^{n+i}$, the differential being
$(-1)^n\partial_{{\cal J}^\mini}$. Then
$\ext^n_{\cal O}({\cal F}^\mini,{\cal G})$ is by definition the
group of homomorphisms of complexes ${\cal F}^\mini\to {\cal
J}^\mini$, modulo homotopy. Of course we could take $\cal G$
itself to be a complex bounded below, but we will not need this
more general case.

Think of $\cal G$ itself as a complex, which is zero everywhere
except in degree 0, and let
${\cal E}^\mini$ be an extension of ${\cal F}^\mini$ by $\cal
G$. To this we can associate an element of $\ext^1_{\cal
O}({\cal F}^\mini,{\cal G})$ as follows. The embedding of
complexes ${\cal G}\to {\cal J}^\mini$ can be extended to
homomorphism of {\bf Z}-graded sheaves $\phi\colon {\cal
E}^\mini\to {\cal J}^\mini$. To
$\phi$ we associate the homomorphism of complexes
$\partial_{{\cal J}^\mini}\phi-\phi\partial_{{\cal
E}^\mini}\colon {\cal E}^\mini\to {\cal J}^\mini[1]$. It sends
$\cal G$ to 0, and so induces a homomorphism
$\psi\colon {\cal F}^\mini\to {\cal J}^\mini[1]$. The homotopy
class of this complex does not depend on $\phi$. Thus we have
defined a function from the set of isomorphism classes of
extensions of ${\cal F}^\mini$ by $\cal G$ to $\ext^1_{\cal
O}({\cal F}^\mini, {\cal G})$. By standard arguments one shows
that this map is a group homomorphism; but in general it is not
injective nor surjective.

\proclaim Lemma. \call{extext1} Assume that ${\cal F}^i = 0$ for
$i>0$. Then every element of $\ext^1_{\cal O}({\cal F}^\mini,
{\cal G})$ is represented by a unique isomorphism class of
extensions.
\endproclaim

\proof Let ${\cal E}^\mini$ be an extension of ${\cal F}^\mini$
by
$\cal G$ such that the corresponding element in $\ext^1_{\cal
O}({\cal F}^\mini, {\cal G})$ is 0. Then I claim that the
embedding
${\cal G}\to {\cal J}^\mini$ extend to a homomorphism of
complexes
${\cal E}^\mini\to {\cal J}^\mini$. In fact, extend ${\cal G}\to
{\cal J}^\mini$ to a homomorphism of graded sheaves $\phi\colon
{\cal E}^\mini\to {\cal J}^\mini$, and call
$\psi\colon {\cal F}^\mini\to {\cal J}^\mini[1]$ the
homomorphism of complexes induced by
$\partial_{{\cal J}^\mini}\phi-\phi\partial_{{\cal E}^\mini}$, as
above. There is a homomorphism of graded sheaves $\lambda\colon
{\cal F}^\mini\to {\cal J}^\mini$ with
$\psi = \partial_{{\cal J}^\mini}\lambda-\lambda\partial_{{\cal
F}^\mini}$. Then $\phi' =
\phi-\lambda\kappa_{{\cal E}^\mini}$ is a homomorphism of
complexes extending the embedding
${\cal G}\into {\cal J}^\mini$.

The homomorphism $\phi'$ induces a homomorphism of sheaves
${\cal H}^0({\cal E}^\mini)\to {\cal H}^0({\cal J}^\mini) =
{\cal G}$. But
${\cal E}^i = 0$ for $i>0$, so this homomorphism ${\cal
H}^0({\cal E}^\mini)\to {\cal G}$ yields a homomorphism ${\cal
E}^\mini\to {\cal G}$ which splits the sequence.

Now start from an element of $\ext^1_{\cal O}({\cal F}^\mini,
{\cal G})$, represented by a homomorphism $f\colon {\cal
F}^\mini\to {\cal J}^\mini[1]$. Let ${\cal I}^\mini\subseteq
{\cal J}$ be the subcomplex with ${\cal I}^0 = {\cal J}^0$,
${\cal I}^1 =
\im({\cal J}^1\to {\cal J}^2) = \im({\cal J}^0\to {\cal J}^1)$,
and
${\cal I}^i = 0$ for $i>1$. Because of the condition on ${\cal
F}^\mini$ we see that the homomorphism  $f\colon {\cal
F}^\mini\to {\cal J}^\mini[1]$ factors through ${\cal
I}^\mini[1]$. There is an extension of ${\cal I}^\mini[1]$ by
$\cal G$ given by the diagram
$$
\def\grado#1{\hbox to 45pt{$\deg\ #1$\hfill}}
\diagram{30}{&&&\vdots&&\vdots&&\vdots\cr
&&&\mapdown&&\mapdown&&\mapdown\cr
\grado{-2}&0&\mapright&0&\mapright&0&\mapright&0&\mapright&0\cr
&&&\mapdown&&\mapdown&&\mapdown\cr
\grado{-1}&0&\mapright&0&\mapright&{\cal J}^0&\bihline&{\cal
I}^0&\mapright&0\cr
&&&\mapdown&&\bivline&&\mapdown\lft\partial\cr
\grado{0}&0&\mapright&{\cal G}&\mapright&{\cal
J}^0&\mapright^{\partial}&{\cal I}^1&\mapright&0\cr
&&&\mapdown&&\mapdown&&\mapdown\cr
\grado{1}&0&\mapright&0&\mapright&0&\mapright&0&\mapright&0\cr
&&&\mapdown&&\mapdown&&\mapdown\cr &&&\vdots&&\vdots&&\vdots\cr}
$$ One checks that the pullback of this extension to ${\cal
F}^\mini$ represents the given class of $\ext^1_{\cal O}({\cal
F}^\mini,{\cal G})$.
\endproof

With this at our disposal the whole theory in this section can be
extended to extensions of complexes with no terms in positive
degree by sheaves. The definition of extension cocycle extends
to this case also, and \ref{extcocgrp} remains true, and but the
proof has to be changed slightly. Here is the main point. If
${\cal A}^\mini$ and
${\cal B}^\mini$ are complexes of sheaves, we will use the
notation $\curshom_{\cal O}({\cal A}^\mini,{\cal B}^\mini)$ to
denote the sheaf of honest homomorphisms of complexes of
${\cal A}^\mini$ into
${\cal B}^\mini$, which is in general very different from
$\cursext^0_{\cal O}({\cal A}^\mini,{\cal B}^\mini)$, and
$\curshom^0_{\cal O}({\cal A}^\mini,{\cal B}^\mini)$ to denote
the sheaf of homomorphisms of graded sheaves of
${\cal A}^\mini$ into
${\cal B}^\mini$.

\proclaim Lemma. Let ${\cal F}^\mini$ be a complex of sheaves
bounded above, $\cal G$ a sheaf, ${\cal J}^\mini$ an injective
resolution of $\cal G$. Then there is a canonical isomorphism of
$\H^1
\bigl(X,\curshom_{\cal O}({\cal F}^\mini,{\cal
J}^\mini[n])\bigr)$ with the kernel of the localization map from
$\ext^{n+1}({\cal F}^\mini,{\cal G})$ to 
$\H^0(X,\cursext^{n+1}\bigl({\cal F}^\mini,{\cal G})\bigr)$.
\endproclaim

\proof Call $A$ the kernel of the localization map
$\ext^{n+1}({\cal F}^\mini,{\cal G})\to
\H^0(X,\cursext^{n+1}\bigl({\cal F}^\mini,{\cal G})\bigr)$; let
us define a homomorphism from $A$ to $\H^1
\bigl(X,\curshom_{\cal O}({\cal F}^\mini,{\cal
J}^\mini[n+1])\bigr)$ as follows. Let $\xi$ be an element of
$A$; then $\xi$ is represented by a homomorphism
$\phi\colon {\cal F}^\mini\to {\cal J}^\mini[n]$. There is an
open  covering $\{X_\alpha\}$ of $X$ and homomorphisms of graded
sheaves $\phi_\alpha\colon {\cal F}^\mini\rest{X_\alpha} \to
{\cal J}^\mini[n]\rest{X_\alpha}$ with $\partial_{{\cal
J}^\mini}\phi_\alpha -
\phi_\alpha\partial_{\cal F}^\mini = \phi\rest{X_\alpha}$. Define
$$
\phi_{\alpha\beta} =
\phi_\alpha-\phi_\beta\in \hom_{\cal O}({\cal
F}^\mini\rest{X_\alpha\cap X_\beta},{\cal
J}^\mini[n]\rest{X_\alpha\cap X_\beta});
$$ then $\{\phi_{\alpha\beta}\}$ is a \v Cech cocycle in
$\curshom_{\cal O}({\cal F}^\mini,{\cal J}^\mini)$. It is easy to
check that its class
$$[\{\phi_{\alpha\beta}\}]\in \H^1
\bigl(X,\curshom_{\cal O}({\cal F}^\mini,{\cal
J}^\mini[n])\bigr)$$ is independent of the choice of the
$\phi_\alpha$, and that in this way we obtain an injective group
homomorphism from $A$ into $\H^1
\bigl(X,\curshom_{\cal O}({\cal F}^\mini,{\cal
J}^\mini[n])\bigr)$.

Let us check surjectivity. Take a class in $\H^1
\bigl(X,\curshom_{\cal O}({\cal F}^\mini,{\cal
J}^\mini[n])\bigr)$, represented to some cocycle
$\phi_{\alpha\beta}$ relative to some open covering
$\{X_\alpha\}$ of $X$. The sheaf
$$\curshom^0_{\cal O}({\cal F}^\mini,{\cal J}^\mini[n]) =
\bigoplus_{i\in {\bf Z}}\curshom_{\cal O}({\cal F}^i,{\cal
J}^{i+n})$$ of homomorphisms of graded groups is a finite direct
sum of flabby sheaves, hence all of its higher \cech{}
cohomology groups are 0; so we can find homomorphisms of graded
sheaves
$\phi_\alpha\colon {\cal F}^\mini\to {\cal J}^\mini[n]$ on
$X_\alpha\cap X_\beta$ with
$\phi_{\alpha\beta} =
\phi_\beta - \phi_\alpha$. The homomorphisms of complexes
$\partial_{{\cal J}^\mini}\phi_\alpha -
\phi_\alpha\partial_{{\cal F}^\mini} \colon {\cal F}^\mini\to
{\cal J}^\mini[n+1]$ patch together to yield a homomorphism of
complexes $\phi\colon {\cal F}^\mini\to {\cal J}^\mini[n+1]$.
The class in $\ext^{n+1}({\cal G}^\mini,{\cal G})\to
\H^0(X,\cursext^{n+1}\bigl({\cal G}^\mini,{\cal G})\bigr)$ maps
to our class in $\H^1
\bigl(X,\curshom_{\cal O}({\cal F}^\mini,{\cal
J}^\mini[n])\bigr)$, and this concludes the proof.
\endproof

Once the have this, to prove the analogue of \ref{extcocgrp} it
is enough to prove the following.

\proclaim Lemma. There is a canonical isomorphism
$$\E_{\cal O}({\cal U};{\cal F}^\mini,{\cal G})\simeq
\cechH^1\bigl({\cal U},\curshom_{\cal O}({\cal F},{\cal
J}^\mini[1])\bigr)$$ which is compatible with restriction maps.
\endproclaim

Let $(\{{\cal E}^\mini_{\alpha\beta},F_{\alpha\beta\gamma}\})$
be an extension cocycle of ${\cal F}^\mini$ by
$\cal G$. Extend the homomorphism from ${\cal G}$ to
${\cal J}^\mini$ to homomorphisms of graded sheaves
$\sigma_{\alpha\beta}\colon {\cal E}^\mini_{\alpha\beta}\to
{\cal J}^\mini$ satisfying the compatibility condition of
\ref{mainextcoc}. The homomorphisms of complexes
$$
\partial_{{\cal
J}^\mini}\sigma_{\alpha\beta}-\sigma_{\alpha\beta}\partial_{{\cal
E}^\mini_{\alpha\beta}}\colon {\cal
E}^\mini_{\alpha\beta}\longrightarrow {\cal J}^\mini[1]
$$ sends $\cal G$ to 0, and induces homomorphisms
$\tau_{\alpha\beta}\colon {\cal F}^\mini\longrightarrow {\cal
J}^\mini[1]$ on $X_\alpha\cap X_\beta$ satisfying the cocycle
condition. These give the class in $\H^1\bigl(X,\curshom_{\cal
O}({\cal F},{\cal J}^\mini[1])\bigr)$ associated with $(\{{\cal
E}^\mini_{\alpha\beta},F_{\alpha\beta\gamma})\}$. The rest of the
proof of \ref{mainextcoc} goes through with obvious
modifications.
\endproof

\beginsection Abstract liftings of \loc s [abstract]

In this section we analyze abstract liftings. \ref{notationring}
is still in force.

\proclaimr Hypotheses. \call{abshyp} Let $X$ be a flat \loc{}
scheme of finite type over $A$. Assume also that $X_0$ is
generically smooth over
$\kappa$.\endproclaim

Recall that for $X$ to be a \loc{} means that if, locally on
$X$, we factor the structure morphism
$X\to \spec A$ as an embedding $X\into P$ followed by a smooth
morphism
$P\to \spec A$, then $X$ is a \loc{} in $P$; this condition is
independent of the factorization. Also, if $\kappa$ is perfect
then for $X_0$ to be generically smooth over $\kappa$ simply
means to be reduced.

\proclaim Definition. An\/ {\rm abstract lifting} of $X$ to
$A'$ is a flat scheme $X'$ over $A'$ with an embedding $X\into
X'$, which induces an isomorphism of $A$-schemes of $X$ with
$X'\rest{\spec A}$.
\endproclaim

\proclaim Definition. Let $X'_1$ and $X'_2$ be two abstract
liftings of $X$. An isomorphism of
$X'_1$ with $X'_2$ is an isomorphism $X'_1\simeq X'_2$ of
schemes over
$A'$ which induces the identity on $X$. \endproclaim

\proclaim Theorem. \call{absdef} Call $\Omega_{X_0/\kappa}$ the
sheaf of K\"ahler differentials of
$X_0$ relative to $\kappa$.

\punto \call{absdefa} Any abstract lifting of $X$ is a \loc{}
over
$A'$.

\punto \call{absdefd} If $X'$ is an abstract lifting of $X$ to
$A'$, then the group of automorphisms of  $X'$ is canonically
isomorphic to
$$
\a\otimes_\kappa\hom_{{\cal
O}_{X_0}}(\Omega_{{X_0}/\kappa},{\cal O}_{X_0}).
$$

\punto \call{absdefb} There is a canonical element
$$\omega_{\rm abs} = \omega_{\rm abs}(X)\in
\a\otimes_\kappa\ext^2_{{\cal
O}_{X_0}}(\Omega_{{X_0}/\kappa},{\cal O}_{X_0}),$$ called the\/
{\rm obstruction} of $X$, such that $\omega_{\rm abs} = 0$ if
and only if an abstract lifting of $X$ to $A'$ exists.

\punto \call{absdefc} If an abstract lifting exists, then there
is a canonical action of the group
$$\a\otimes_\kappa\ext^1_{{\cal
O}_{X_0}}(\Omega_{{X_0}/\kappa},{\cal O}_{X_0})$$ on the set of
isomorphism classes of abstract liftings making it into a
principal homogeneous space.
\endproclaim

The condition that $X_0$ be generically smooth is necessary for
the statement to hold. For example, consider that case that $A =
\kappa$, $A' = \kappa[x]/(x^2)$,
$X = X_0 = \spec\bigl(\kappa[t]/(t^2)\bigr)$. Then $X$ has many
non-isomorphic abstract lifting to
$A'$, given by $X' = \spec\bigl(\kappa[x,t]/(t^2-ax)\bigr)$, for
$a\in
\kappa$; but $\ext^1_{{\cal
O}_{X_0}}(\Omega_{{X_0}/\kappa},{\cal O}_{X_0})$ is 0.

Let us place ourselves again in the situation of \ref{embhyp}.
We need to understand the isomorphisms between liftings of $X$
to $M'$. Let $\Phi\colon X'_2\simeq X'_1$ be such an
isomorphism; $\Phi$ will yield an isomorphism of sheaves of
$A'$-algebras
$\phi\colon {\cal O}_{X'_1}\simeq {\cal O}_{X'_2}$, inducing the
identity on ${\cal O}_X = {\cal O}_{X'_1}/\a{\cal O}_{X'_1} =
{\cal O}_{X'_2}/\a{\cal O}_{X'_2}$. Consider the two projections
$\pi_i\colon {\cal O}_{M'}\to {\cal O}_{X'_i}$; the difference
$$D = \pi_2-\phi\pi_1\colon {\cal O}_{M'}\to {\cal O}_{X'_2}$$
will have its image inside $\a {\cal O}_{X'_2} = \a\otimes_{A'}
{\cal O}_{X'_2} =
\a\otimes_\kappa {\cal O}_{X_0}$, so we think of $D$ as a
function ${\cal O}_{M'}\to
\a\otimes_\kappa{\cal O}_{X_0}$. If $f$ and $g$ are local
sections of ${\cal O}_{M'}$, we have
$$\eqalign{D(fg) &= \pi_2(fg)-\phi\pi_1(fg) \cr &=
\pi_2(f)\pi_2(g)-
\pi_2(f)\phi\pi_1(g)+
\pi_2(f)\phi\pi_1(g) -
\phi\pi_1(f)\phi\pi_1(g)\cr &=\pi_2(f)D(g)+ D(f)\phi\pi_1(g).}$$
Notice that the two ${\cal O}_{M'}$-module structures on
$\a\otimes_\kappa{\cal O}_{X_0} = \a {\cal O}_{X'_2}$ induced by
$\pi_1$ and by $\phi\pi_2$ coincide, so we can think of $D$ is
derivation of ${\cal O}_{M'}$ into the ${\cal O}_{M'}$-module
$\a\otimes_\kappa{\cal O}_{X_0}$. Because $\phi$ is $A'$-linear
we have also that $D$ kills the elements of $A'$, so it is an
$A'$-derivation. Such a derivation will send all the elements of
the annihilator of ${\cal O}_{M_0}$ to 0, so we think of
$D$ as a
$\kappa$-derivation of ${\cal O}_{M_0}$ into
$\a\otimes_\kappa{\cal O}_{X_0}$, that is, as an element of
$\der_\kappa({\cal O}_{M_0},\a\otimes_\kappa{\cal O}_{X_0}) =
\hom_{{\cal
O}_{X_0}}(\Omega_{M_0/\kappa}\rest{X_0},\a\otimes_\kappa{\cal
O}_{X_0})$. The derivation
$D$ is not arbitrary.

{}From the homomorphism ${\cal I}_0/{\cal I}_0^2\to
\Omega_{M_0/\kappa}\rest{X_0}$ we get a restriction map
$$\hom_{{\cal
O}_{X_0}}(\Omega_{M_0/\kappa}\rest{X_0},\a\otimes_\kappa{\cal
O}_{X_0})\longrightarrow\hom_{{\cal O}_{X_0}}({\cal I}_0/{\cal
I}_0^2,\a\otimes_\kappa{\cal O}_{X_0}) = 
\H^0(X_0,\a\otimes_\kappa{\cal N}_0).$$

\proclaim Lemma. The restriction of $D$ to ${\cal I}_0/{\cal
I}_0^2$ is $\nu(X'_1,X'_2)$.
\endproclaim

\proof Let ${\cal I}'_1$ and ${\cal I}'_2$ be the ideals of
$X'_1$ and
$X'_2$ respectively. Take a local section $f$ of
$\cal I$, and lift it to two local sections $f'_1$ of ${\cal
I}'_1$ and
$f'_2$ of ${\cal I}'_2$. Then
$$ Df = \pi_2(f'_1)-\phi\pi_1(f'_1) = \pi_2(f'_1) =
\pi_2(f'_1-f'_2).
$$ But because $f'_1-f'_2$ is a local section of $\a{\cal
O}_{M'} =
\a\otimes_\kappa{\cal O}_{M_0}$, the homomorphism $\pi_2$ will
send it into its image in
$\a\otimes_\kappa{\cal O}_{X_0}$, which is by definition
$\nu(X'_1,X'_2)(f)$.
\endproof

We leave it as an exercise for the interested reader to check
that by assigning $D$ to $\Phi$ we get a bijective
correspondence between isomorphisms of liftings $X'_2\simeq
X'_1$ with elements of
$$\hom_{{\cal
O}_{X_0}}({\Omega_{M_0/\kappa}\rest{X_0},\a\otimes_\kappa{\cal
O}_{X_0}})$$ whose image in
$\H^0(X_0,\a\otimes_\kappa{\cal N}_0)$ is $\nu(X'_1,X'_2)$. This
bijective correspondence has the following properties.

\proclaim Proposition. \call{isomlif} There is a bijective
correspondence between isomorphisms of abstract liftings
$X'_2\simeq X'_1$ and elements of\/
$\hom_{{\cal
O}_{X_0}}({\Omega_{M_0/\kappa}\rest{X_0},\a\otimes_\kappa{\cal
O}_{X_0}})$ whose image in
$\H^0(X_0,\a\otimes{\cal N}_0)$ is $\nu(X'_1,X'_2)$, with the
following properties.

\punto \call{isomid} The identity ${\rm id}_{X'}\colon X'\simeq
X'$ corresponds to 0.

\punto \call{isomcoc} If $\Phi_1\colon X'_2\simeq X'_1$ and
$\Phi_2\colon X'_3\simeq X'_2$ are isomorphisms of abstract
liftings and $D_1$, $D_2$ are the corresponding elements of
$\hom_{{\cal
O}_{X_0}}({\Omega_{M_0/\kappa}\rest{X_0},\a\otimes_\kappa{\cal
O}_{X_0}})$, then the composition
$\Phi_1\circ \Phi_2$ corresponds to $D_1+D_2$.

\punto \call{isomopen} Let $Y$ be an open subset of $X$, $Y'_1$
and
$Y'_2$ the restrictions of
$X'_1$ and
$X'_2$ to $Y$, $\Phi\colon X'_1\simeq X'_2$ an isomorphism of
liftings, corresponding to
$D\in \hom_{{\cal
O}_{X_0}}({\Omega_{M_0/\kappa}\rest{X_0},\a\otimes_\kappa{\cal
O}_{X_0}})$. Then the element of
$\hom_{{\cal
O}_{Y_0}}({\Omega_{M_0/\kappa}\rest{Y_0},\a\otimes_\kappa{\cal
O}_{Y_0}})$ corresponding to the restriction
$\Phi\rest{Y'_1}\colon Y'_1\to Y'_2$ is the restriction of
$D$.
\endproclaim

The rest of the proof is entirely based on the constructions of
\ref{extensions}, to which the reader should refer for the
notation.

Let us put ourselves in the situation of \ref{abshyp}. Let $X'$
be an abstract lifting of $X$; if we set $M' = X' = X'_1 = X'_2$
in \ref{isomlif} we see that we have proved
\ref{absdefd}.

Let $X'_1$ and $X'_2$ be two abstract liftings of $X$; we want
to know when they are isomorphic. Assume that there exists a
smooth morphism $P'\to \spec{A'}$ and an embedding of
$X$ into $P = P'\rest{\spec A}$ which lifts to embeddings of
$X'_1$ and $X'_2$ into $P'$. Choose such liftings; we obtain an
${\cal O}_{X_0}$linear map
$\nu(X'_1,X'_2)\colon {\cal I}_0/{\cal I}_0^2\to
\a\otimes_\kappa{\cal O}_{X_0}$.

The following lemma is where we use the hypothesis that $X_0$ be
generically smooth.

\proclaim Lemma. \call{exseq} The usual sequence
$$\sequence{0&\mapright {\cal I}_0/{\cal I}_0^2&\mapright
\Omega_{P_0/\kappa}\rest{X_0} &\mapright
\Omega_{X_0/\kappa}&\mapright 0}$$ is exact.
\endproclaim

We will call this sequence the {\it fundamental exact
sequence\/} for the embedding of $X$ in $P$.

\proof This is standard, except for the injectivity of the first
arrow. But it is well known that this arrow is injective where
$X_0$ is smooth, so its kernel is concentrated on a nowhere
dense closed subset of $X_0$, because
$X_0$ is generically smooth. Since ${\cal I}_0/{\cal I}_0^2$ is
locally free on
$X_0$ and $X_0$ has no embedded point, being a \loc{} scheme
over a field, we see that the kernel must actually be 0.
\endproof

If $X_0$ is not generically smooth then the lemma does not hold
in general; the remedy is to substitute the complex ${\cal
I}_0/{\cal I}_0^2\to
\Omega_{P_0/\kappa}\rest{X_0}$ for the sheaf
$\Omega_{X_0/\kappa}$ in the statement of the theorem, and in the
remainder of the proof. This complex is defined up to a
canonical isomorphism in the derived category of coherent
sheaves on $X_0$, and it is the simplest nontrivial example of a
cotangent complex.

We define
$${\cal E}(X'_1,X'_2) \equaldef
\nu(X'_1,X'_2)_*\bigl(\Omega_{P_0/\kappa}\rest{X_0}\bigr)$$ as
an extension of
$\Omega_{X_0/\kappa}$ by $\a\otimes_\kappa{\cal O}_{X_0}$. We
need to check that this extension is independent of the choices
made, so choose two smooth morphisms $P'_1\to
\spec{A'}$ and $P'_2\to \spec{A'}$, embeddings
$X\into P_i$ and liftings $X'_1\into P'_i$ and $X'_2\into P'_i$.
These induce an embedding $X\into P_1\times _{\spec A}P_2$, and
liftings
$X'_1\into P'_1\times _{\spec {A'}}P'_2$ and $X'_2\into
P'_1\times _{\spec {A'}}P'_2$. Let ${\cal C}_1$,
${\cal C}_2$ and
${\cal C}_{12}$ be the conormal bundles of $X_0$ in $(P_1)_0$,
$(P_2)_0$, and
$(P'_1\times _{\spec {A'}}P'_2)_0 = (P_1)_0\times_\kappa (P_2)_0$
respectively. Denote by
$\nu_i\colon {\cal C}_i\to \a\otimes_\kappa{\cal O}_{X_0}$ the
corresponding sections of the normal bundles, and set
${\cal E}_i = {\nu_i}_*\Omega_{(P_i)_0/\kappa}$, ${\cal E}_{12} =
{\nu_{12}}_*\Omega_{(P_1)_0\times_\kappa (P_2)_0/\kappa}$. There
are homomorphisms
$\phi_i\colon {\cal C}_i\to {\cal C}_{12}$ fitting into
commutative diagrams with exact rows
$$
\diagram{30}{0&\mapright &{\cal
C}_i&\mapright&\Omega_{(P_i)_0/\kappa}\rest{X_0}&\mapright
&\Omega_{X_0/\kappa}&\mapright &0\cr
&&\mapdown\lft{\phi_i}&&\mapdown&&\bivline\cr 0&\mapright &{\cal
C}_{12}&\mapright&\Omega_{(P_1)_0\times_\kappa
(P_2)_0/\kappa}\rest{X_0}&\mapright
&\Omega_{X_0/\kappa}&\mapright &0\cr
&&\mapdown\lft{\nu_{12}}&&\mapdown&&\bivline\cr 0&\mapright
&\a\otimes_\kappa {\cal O}_{X_0}&\mapright&{\cal
E}_{12}&\mapright &\Omega_{X_0/\kappa}&\mapright &0\cr}
$$ Because of \ref{nucomp} we have that $\nu_{12}\circ\phi_i =
\nu_i\colon {\cal C}_i\to \a\otimes_\kappa{\cal O}_{X_0}$. By
\ref{homext} the diagram above induces an isomorphism of
extensions $\alpha_i\colon {\cal E}_i\simeq {\cal E}_{12}$; we
define the canonical isomorphism between ${\cal E}_2$ and ${\cal
E}_1$ to be
$\alpha_{12} =
\alpha_1\circ
\alpha_2^{-1}$. One checks that these isomorphisms satisfy the
cocycle condition, that is, if $P'_1$, $P'_2$ and $P'_3$ are
smooth over $A'$, and are given embeddings $X\into P_j$ and
liftings $X'_i\into P'_j$ for $i = 1$, 2 and
$j = 1$, 2,~3, then 
$\alpha_{13} =
\alpha_{12}\circ\alpha_{23}$, with the obvious notation.

{}From the construction above one sees clearly that if $Y$ is an
open subset of $X$, $Y'_1$ and $Y'_2$ are the restrictions of
$X'_1$ and $X'_2$ to $Y$, then the extension associated to the
abstract liftings $Y'_1$ and $Y'_2$ of $Y$ is simply the
restriction to $Y$ of the extension associated to $X'_1$ and
$X'_2$.

Such a morphism $P'\to \spec{A'}$ might not exist globally, but
certainly exists locally.

\proclaim Lemma. \call{locP'ex} Assume that $X$ is affine, and
$X\into {\bf A}^n_A$ is a closed embedding. If $X'$ is an
abstract lifting of $X$, then the embedding $X\into {\bf A}^n_A$
can be extended to an embedding $X'\into {\bf A}^n_{A'}$.
\endproclaim

\proof Set  $X =
\spec B$; the closed embedding
$X\into {\bf A}_A^n$ corresponds to a surjective homomorphism of
$A$-algebras
$A[x_1,\ldots,x_n]\to B$. Set
$X' = \spec B'$; by choosing liftings of the images of the $x_j$
to
$B'$ we obtain a homomorphism of $A'$-algebras
$A'[x_1,\ldots,x_n]\to B'$, which is easily shown to be
surjective, due once again to the fact that $\a$ is nilpotent.
\endproof

In the general case, let $\{X_\alpha\}$ be a covering of $X$
with open affine subschemes, and call $X'_{\alpha i}$ the
restriction of $X'_i$ to $X_\alpha$. For each
$X_\alpha$ choose a closed embedding $X_\alpha\into {\bf
A}^{n_\alpha}_A$, and extend to embeddings of $X'_{\alpha i}$
into ${\bf A}^{n_\alpha}_{A'}$.

For each $\alpha$ set ${\cal E}_\alpha = {\cal E}(X'_{\alpha
1},X'_{\alpha 2})$. Because of the above there are isomorphisms
of extensions of
${\cal E}_\alpha\rest{X_\alpha\cap X_\beta}$ with ${\cal
E}_\beta\rest{X_\alpha\cap X_\beta}$, satisfying the cocycle
condition. We use these to glue the
${\cal E}_\alpha$ together into an extension ${\cal
E}(X'_1,X'_2)$. It is clear that this is defined up to a
canonical isomorphism.

\proclaim Proposition. Given two abstract liftings $X'_1$ and
$X'_2$, there is an extension
${\cal E}(X'_1,X'_2)$ of $\Omega_{X_0/\kappa}$ by
$\a\otimes_\kappa{\cal O}_{X_0}$, well defined up to canonical
isomorphisms, with the following properties.

\punto \call{extid} If $X'$ is a lifting, ${\cal E}(X',X')$ is
canonically isomorphic to
${\bf 0}_{\Omega_{X_0/\kappa},\a\otimes_\kappa{\cal O}_{X_0}}$.

\punto \call{extopp} ${\cal E}(X'_2,X'_1)$ is $-{\cal
E}(X'_1,X'_2)$.

\punto\call{extcoc} Given three abstract liftings $X'_1$, $X'_2$
and
$X'_3$, there a canonical isomorphism
$$F_{X'_1,X'_2,X'_3}\colon {\cal E}(X'_1,X'_2)+{\cal
E}(X'_2,X'_3)\simeq  {\cal E}(X'_1,X'_3).$$

If $X'_1$, $X'_2$, $X'_3$ and $X'_4$ are abstract liftings, then
$$ F_{X'_1,X'_3,X'_4}\circ\bigl(F_{X'_1,X'_2,X'_3}+{\rm
id}_{{\cal E}(X'_3,X'_4)}\bigr) = 
F_{X'_1,X'_2,X'_4}\circ\bigl({\rm id}_{{\cal
E}(X'_1,X'_2)}+F_{X'_2,X'_3,X'_4}\bigr)
$$ as a homomorphism ${\cal E}(X'_1,X'_2)+{\cal
E}(X'_2,X'_3)+{\cal E}(X'_3,X'_4)\longrightarrow {\cal
E}(X'_1,X'_4)$.

\punto \call{extopen} If $Y$ is an open subscheme of $X$ and
$Y'_1$,
$Y'_2$ are the restrictions of $X'_1$ and $X'_2$ to $Y$, then
${\cal E}(Y'_1,Y'_2)$ is the restriction of
${\cal E}(X'_1,X'_2)$ to $Y$.

\punto \call{extconst} If there is given a smooth morphism $P'\to
\spec{A'}$ and an embedding of
$X$ into $P = P'\rest{\spec A}$ which lifts to embeddings of
$X'_1$ and $X'_2$ into $P'$, then
$${\cal E}(X'_1,X'_2) =
\nu(X'_1,X'_2)_*\bigl(\Omega_{P_0/\kappa}\rest{X_0}\bigr)$$
where $\Omega_{P_0/\kappa}\rest{X_0}$ is considered as an
extension of
$\Omega_{X_0/\kappa}$ by $\a\otimes_\kappa{\cal O}_{X_0}$ via
the fundamental exact sequence of the embedding of $X$ in $P$.

\punto\call{exttrans} Given an extension
${\cal E}$ of $\Omega_{X_0/\kappa}$ by $\a\otimes_\kappa{\cal
O}_{X_0}$ and an abstract lifting $X'$ of $X$, there is an
abstract lifting $\widetilde X'_1$ of
$X$ such that ${\cal E}(\widetilde  X',X')$ is isomorphic to
${\cal E}$.
\endproclaim

\proof For part \ref{part extid}, choose embeddings
$X_\alpha\subseteq {\bf A}_A^{n_\alpha}$, and liftings
$X'_\alpha\subseteq {\bf A}_{A'}^{n_\alpha}$; then
$\nu(X'_\alpha,X'_\alpha) = 0$ (\ref{nuid}) and so the statement
follows from \ref{pushzero}.

Part \ref{part extopp} follows from \ref{nuopp} and
\ref{pushextopp}.

The isomorphism in part \ref{part extcoc} is obtained from
\ref{nuadd} and \ref{pushextb}. We leave it the reader to check
the compatibility condition.

Part \ref{part extopen} is clear.

Part \ref{part extconst} follows from the construction.

The proof of part \ref{part exttrans} will be given after
\ref{splisom}.
\endproof

A splitting of ${\cal E}(X'_1,X'_2)$ corresponds to a collection
of splittings of each
${\cal E}_\alpha$, compatible with the gluing isomorphisms. But
by
\ref{pushsplit} these  splittings correspond to homomorphisms of
$\Omega_{{\bf A}^{n_\alpha}_\kappa/\kappa}$ into
$\a\otimes_\kappa {\cal O}_{X_0}$ whose restriction to the
conormal sheaf ${\cal C}_\alpha$ of $X_\alpha$ is
$\nu(X'_{\alpha 1},X'_{\alpha 2})$; because of \ref{isomlif} this
means exactly a compatible system of isomorphisms of liftings
between $X'_{\alpha 2}$ and
$X'_{\alpha 1}$, that is, an isomorphism of abstract liftings of
$X'_2$ with $X'_1$.

\proclaim Proposition. \call{splisom} There is a natural
bijective correspondence of isomorphisms of abstract liftings
$X'_2\simeq X'_1$ with splittings of
${\cal E}(X'_1,X'_2)$.

If $X'_1$, $X'_2$ and $X'_3$ are abstract liftings, $\phi_1\colon
X'_2\simeq X'_1$ and
$\phi_2\colon X'_3\simeq X'_2$ are isomorphisms corresponding to
splittings $s_1$ and
$s_2$ of ${\cal E}(X'_1,X'_2)$ and ${\cal E}(X'_2,X'_3)$, then
the composition
$\phi_1\phi_2$ corresponds to the splitting $(s_1+s_2) \circ
F_{X'_1,X'_2,X'_3}^{-1}$ of
${\cal E}(X'_1,X'_3)$.
\endproclaim

We leave it to the reader to unwind the various definitions and
prove the Proposition.

\proofof exttrans.  Let
$\{X_\alpha\}$ be a covering of $X$ by affine open subschemes,
and call $X'_{\alpha}$ the restriction of $X'$ to $X_\alpha$.
For each $\alpha$ we choose an embedding
$X'_\alpha\into {\bf A}^{n_\alpha}_{A'}$, and call ${\cal
C}_\alpha$ the conormal bundle of
$X_\alpha\cap X_0$ in ${\bf A}^{n_\alpha}_\kappa$. The
fundamental exact sequence
$$
\sequence{0&\mapright &{\cal C}_\alpha&\mapright &\Omega_{{\bf
A}^{n_\alpha}_\kappa/\kappa}\rest{X_\alpha\cap X_0}&\mapright
&\Omega_{{\cal O}_{X_\alpha\cap X_0}}&
\mapright&0}
$$ induces an exact sequence of abelian groups
$$
\lunghezza{30}\eqalign{\hom_{{\cal O}_{X_\alpha\cap X_0}}({\cal
C}_\alpha,\a\otimes_\kappa{\cal O}_{X_\alpha\cap
X_0})&\mapright^\partial\ext^1_{{\cal O}_{X_\alpha\cap
X_0}}(\Omega_{X_0/\kappa}\rest{X_\alpha\cap
X_0},\a\otimes_\kappa{\cal O}_{X_\alpha\cap X_0})\cr &\mapright
\ext^1_{{\cal O}_{X_\alpha\cap X_0}}(\Omega_{{\bf
A}^{n_\alpha}_\kappa/\kappa}\rest{X_\alpha\cap
X_0},\a\otimes_\kappa{\cal O}_{X_\alpha\cap X_0}) = 0}
$$ where the last group is 0 because ${X_\alpha\cap X_0}$ is
affine and
$\Omega_{{\bf A}^{n_\alpha}_\kappa/\kappa}\rest{X_\alpha\cap
X_0}$ is locally free. This means that
$\partial$ is surjective, so, by \ref{pushextbound}, for each
$\alpha$ we can find a homomorphism 
$f_\alpha\colon {\cal C}_\alpha\to \a\otimes_\kappa{\cal
O}_{X_\alpha\cap X_0}$ and an isomorphism
${f_\alpha}_*\Omega_{{\bf
A}^{n_\alpha}_\kappa/\kappa}\rest{X_\alpha\cap X_0}\simeq {\cal
E}\rest{X_\alpha\cap X_0}.$ But \ref{nutrans} implies that there
exists a lifting $\widetilde X'_\alpha$ of $X_\alpha$ in ${\bf
A}^{n_\alpha}_{A'}$ such that
$\nu(\widetilde X'_\alpha,X'_\alpha) = f_\alpha$, and from the
construction we get an isomorphism of extensions
$\phi_\alpha\colon {\cal E}(\widetilde X'_\alpha,X'_\alpha)\simeq
{\cal E}\rest{X_\alpha\cap X_0}$. Consider the composition
$$
\phi_\alpha^{-1}\phi_\beta\colon {\cal E}(\widetilde
X'_\beta,X'_\beta)\rest{X_\alpha\cap X_\beta \cap X_0}\simeq
{\cal E}(\widetilde X'_\alpha,X'_\alpha)\rest{X_\alpha\cap
X_\beta \cap X_0}
$$ By \ref{extisom} the isomorphisms
$\phi_\alpha^{-1}\phi_\beta$ correspond to splittings of
$$ {\cal E}(\widetilde X'_\beta,X'_\beta)-{\cal E}(\widetilde
X'_\alpha,X'_\alpha) = {\cal E}(\widetilde
X'_\beta,X'_\beta)+{\cal E}(X'_\alpha,\widetilde X'_\alpha)
\simeq {\cal E}(\widetilde X'_\beta,\widetilde X'_\alpha)
$$ By \ref{splisom} these splittings yield isomorphisms of
liftings
$X'_\beta\rest{X'_\alpha\cap X'_\beta}\simeq
X'_\alpha\rest{X'_\alpha\cap X'_\beta}$. One proves that they
satisfy the cocycle condition: therefore we can glue the various
$\widetilde X'_\alpha$ together to find the desired $\widetilde
X'$.
\endproof

If $X'_1$ and $X'_2$ are abstract liftings of $X$, call
$e(X'_1,X'_2)$ the class of ${\cal E}(X'_1,X'_2)$ in
$$
\ext^1_{{\cal
O}_{X_0}}(\Omega_{X_0/\kappa},\a\otimes_\kappa{\cal O}_{X_0}) =
\a\otimes_\kappa\ext^1_{{\cal O}_{X_0}}(\Omega_{X_0/\kappa},{\cal
O}_{X_0}).
$$ The various properties of the extension ${\cal E}(X'_1,X'_2)$
can be translated as follows.

\proclaim Proposition. \call{epsilon}Given two abstract liftings
$X'_1$ and $X'_2$, there is well defined element
$$e(X'_1,X'_2)\in \ext^1_{{\cal
O}_{X_0}}(\Omega_{X_0/\kappa},\a\otimes_\kappa{\cal O}_{X_0}) =
\a\otimes_\kappa\ext^1_{{\cal
O}_{X_0}}(\Omega_{X_0/\kappa},{\cal O}_{X_0})$$ with the
following properties.

\punto \call{epsilonisom}$e(X'_1,X'_2) = 0$ if and only if
$X'_1$ and
$X'_2$ are isomorphic.

\punto \call{epsilonadd}If $X'_1$, $X'_2$ and $X'_3$ are abstract
liftings of $X$, then
$e(X'_1, X'_2)+e(X'_2,X'_3) = e(X'_1,X'_3)$.

\punto \call{nuext} If there exists a smooth morphism $P'\to
\spec{A'}$ and an embedding of
$X$ into $P = P'\rest{\spec A}$ which lifts to embeddings of
$X'_1$ and $X'_2$ into $P'$, then
$$e(X'_1,X'_2) = \partial
\nu(X'_1,X'_2),$$ where
$$\partial\colon \H^0(X_0,{\cal N}_0) = \hom({\cal I}_0/{\cal
I}_0^2,\a\otimes_\kappa{\cal O}_{X_0})\to
\ext^1_{{\cal O}_{X_0}}(\Omega_{X_0/\kappa},\a\otimes_\kappa{\cal
O}_{X_0})$$ is the boundary homomorphism coming from the
fundamental exact sequence
$$\sequence{0&\mapright&{\cal I}_0/{\cal I}_0^2&\mapright
\Omega_{P_0/\kappa}\rest{X_0} &\mapright
\Omega_{X_0/\kappa}&\mapright 0}.$$
\endproclaim

\proof Part \ref{part epsilonisom} follows from \ref{splisom}.

Part \ref{epsilonadd} follows from \ref{extcoc}.

Part \ref{part nuext} follows from \ref{extconst} and
\ref{pushextbound}\endproof

Because of \ref{homogsp} we see that \ref{absdefc} follows.

We still have to prove \ref{absdefb}. First we will give an easy
proof under a simplifying hypothesis.

\proclaim Hypothesis. \call{simphyp} Assume that there exists a
smooth morphism
$P'\to
\spec{A'}$ and an embedding of $X$ into $P = P'\rest{\spec A}$.
\endproclaim

The only general case in which I know that \ref{simphyp} is true
is when
$X$ is quasiprojective over $\kappa$. However, quasiprojectivity
is not a very natural hypothesis; for example, if $X_0\subseteq
{\bf P}^3_{\bf C}$ is a smooth quartic surface with Picard
number 1 and
$X$ is a general lifting of $X_0$ to the ring of dual numbers,
then the Picard group of $X$ is trivial, so $X$ can not be
projective. This problem does not arise for curves, that is, if
$A$ is artinian and $X$ is one-dimensional then
$X$ is quasiprojective.

Assume that \ref{simphyp} holds, and choose such a factorization
$X\into P\to \spec A$. Call ${\cal I}_0$ the ideal of $X_0$ in
$P_0$, and consider the fundamental exact sequence
$$\sequence{0&\mapright {\cal I}_0/{\cal I}_0^2&\mapright
\Omega_{P_0/\kappa}\rest{X_0}
&\mapright\Omega_{X_0/\kappa}&\mapright 0}.$$ We define the
obstruction $\omega_{\rm abs}\in \ext^2_{{\cal
O}_{X_0}}(\Omega_{X_0/\kappa},\a\otimes_\kappa{\cal O}_{X_0}) =
\a\otimes_\kappa\ext^2_{{\cal
O}_{X_0}}(\Omega_{X_0/\kappa},{\cal O}_{X_0})$ to be the image
of the embedded obstruction
$\omega_{\rm emb}\in \H^1(X_0,\a\otimes_\kappa{\cal N}_0) =
\ext^1_{{\cal O}_{X_0}}({\cal I}_0/{\cal
I}_0^2,\a\otimes_\kappa{\cal O}_{X_0})$ by the boundary map
$$\partial \colon \ext^1_{{\cal O}_{X_0}}({\cal I}_0/{\cal
I}_0^2,\a\otimes_\kappa{\cal O}_{X_0})\longrightarrow
\ext^2_{{\cal O}_{X_0}}(\Omega_{X_0/\kappa},\a\otimes_\kappa{\cal
O}_{X_0})$$ coming from the fundamental exact sequence
$$\sequence{0&\mapright &{\cal I}_0/{\cal I}_0^2&\mapright
\Omega_{P_0/\kappa}\rest{X_0} &\mapright
\Omega_{X_0/\kappa}&\mapright 0}.$$

We need to show that $\omega_{\rm abs}$ does not depend on $P'$,
and that it is 0 if and only if an abstract lifting exists. Set
${\cal T}_0 = \curshom_{{\cal
O}_{X_0}}(\Omega_{P_0/\kappa}\rest{X_0},\a\otimes_\kappa{\cal
O}_{X_0})$, and observe that
$\omega_{\rm abs} = 0$ if and only if
$\omega_{\rm emb}$ is in the image of the  map
$$\rho\colon \H^1(X_0,{\cal T}_0) = \ext^1_{{\cal
O}_{X_0}}(\Omega_{P_0/\kappa}\rest{X_0},\a\otimes_\kappa{\cal
O}_{X_0})\to
\ext^1_{{\cal O}_{X_0}}({\cal I}_0/{\cal
I}_0^2,\a\otimes_\kappa{\cal O}_{X_0}) =
\a\otimes_\kappa\H^1(X_0,{\cal N}_0)$$ induced by the
fundamental exact sequence.

Notice first of all that if $X$ is affine then $\omega_{\rm
emb}\in
\H^1(X_0,{\cal N}_0) = 0$ vanishes, so a lifting $X'$ exists.
Also, the map
$$\partial\colon \hom({\cal I}_0/{\cal I}_0^2,
\a\otimes_\kappa{\cal O}_{X_0})\to
\ext^1_{{\cal O}_{X_0}}(\Omega_{X_0/\kappa},\a\otimes_\kappa{\cal
O}_{X_0})$$ of \ref{nuext} is surjective, so if $\widetilde X'$
is another abstract lifting of $X$ there is $\nu\in
\a\otimes_\kappa\H^0(X_0,{\cal N}_0)$ with $\partial \nu =
e(\widetilde X',X')$. But such a $\nu$ is of the form $\nu =
\nu(\widehat X',X')$ for a certain lifting $\widehat X'\subseteq
P'$ because of
\ref{nutrans}, and from \ref{nuext} we get that $e(\widetilde
X',X') = e(\widehat X',X')$.  Then
$$ e(\widehat X',\widetilde X') = e(\widehat X',X') -
e(\widetilde X',X') = 0,
$$ because of \ref{epsilon}. It follows that $\widetilde X'$ is
isomorphic to $\widehat X'$. So we have shown that if $X$ is
affine then an abstract lifting exists, and any abstract lifting
can be embedded in $P'$. The existence of an abstract lifting in
the affine case is consistent with the fact that
$\ext^2_{{\cal
O}_{X_0}}(\Omega_{X_0/\kappa},\a\otimes_\kappa{\cal O}_{X_0}) =
0$, because
$\Omega_{X_0/\kappa}$ has projective dimension 1, as a
consequence of
\ref{exseq}.

Now choose a covering $\{X_\alpha\}$ of $X$ by open affine
subsets. Then because of the discussion above an abstract
lifting exists if and only if there are liftings $X'_\alpha$ of
$X_\alpha$ in $P'$ and isomorphisms $\phi_{\alpha\beta}\colon
X'_\beta\rest{X_\alpha\cap X_\beta}\simeq
X'_\alpha\rest{X_\alpha\cap X_\beta}$ satisfying the cocycle
condition. But invoking \ref{isomlif} and \ref{descembcoc} we
see that this is true if and only if there exists a cocycle
$\{\nu_{\alpha\beta}\}$ in $\omega_{\rm emb}$, and a collection
$\{D_{\alpha\beta}\}$ of elements of $\H^0(X_\alpha\cap
X_\beta\cap X_0,{\cal T}_0)$ such that the restriction of
$D_{\alpha\beta}$ to $\H^0(X_\alpha\cap X_\beta\cap X_0, {\cal
N}_0)$ is $\nu_{\alpha\beta}$. But that this is the case if and
only if
$\omega_{\rm emb}$ is in the image of $\H^1(X_0,{\cal T}_0)$,
which is what we need.

Let us check that $\omega_{\rm abs}$ is independent of $P'$ and
of the embedding $j\colon X\into P$. Let
$\widetilde \jmath\colon X\into \widetilde P = \widetilde
P'\rest{\spec A}$ be another embedding. By the usual method of
considering the fiber product $P'\times _{\spec {A'}}\widetilde
P'$ we may assume that there exists a smooth morphism
$\pi'\colon \widetilde P'\to P'$ such that if we denote by $\pi$
the restriction of $\pi'$ to
$\spec A$ we have
$\pi\widetilde
\jmath = j$. If ${\cal I}_0$ is the ideal of $X_0$ in $P_0$ and
$\widetilde {\cal I}_0$ is the ideal of $X_0$ in
$\widetilde P_0$, we get a commutative diagram with exact rows
$$
\diagram{30}{0&\mapright&{\cal I}_0/{\cal I}_0^2&\mapright
&\Omega_{P_0/\kappa}\rest{X_0}&\mapright
&\Omega_{X_0/\kappa}&\mapright &0\cr
&&\mapdown&&\mapdown^{\pi^*}&&\bivline\cr 0&\mapright&\widetilde
{\cal I}_0/\widetilde {\cal I}_0^2&\mapright &\Omega_{\widetilde
P_0/\kappa}\rest{X_0}&\mapright &\Omega_{X_0/\kappa}&\mapright
&0\cr}
$$ inducing a commutative diagram
$$\diagram{40}{\ext^1({\cal I}_0/{\cal
I}_0^2,\a\otimes_\kappa{\cal O}_{X_0})&\mapright^{\partial }& 
\ext^2(\Omega_{X_0/\kappa},\a\otimes_\kappa{\cal O}_{X_0})\cr
\mapup&&\bivline\cr
\ext^1_{{\cal O}_X{_0}}(\widetilde {\cal I}_0/\widetilde {\cal
I}_0^2,\a\otimes_\kappa{\cal O}_{X_0})&\mapright^{\partial}& 
\ext^2_{{\cal O}_X{_0}}(\Omega_{X_0/\kappa},\a\otimes_\kappa{\cal
O}_{X_0}).\cr}$$ According to \ref{embobscomp} the first column
carries the embedded obstruction of $X$ in
$\widetilde P$ into the embedded obstruction of $X$ in $P$, and
so the two images in
$\ext^1_{{\cal
O}_{X_0}}(\Omega_{X_0/\kappa},\a\otimes_\kappa{\cal O}_{X_0})$
coincide.

To prove \ref{absdefb} in general, we need the machinery of
extension cocycles developed in
\ref{extensions}. Let $\{X_\alpha\}$ be a covering of $X$ with
affine subsets, and for each $\alpha$ let $X'_\alpha$ be an
abstract lifting of $X_\alpha$. For each triple
$\alpha$, $\beta$ and $\gamma$ consider the isomorphisms
$$ F_{X'_\alpha,X'_\beta,X'_\gamma}\colon {\cal
E}(X'_\alpha,X'_\beta)+ {\cal E}(X'_\beta,X'_\gamma)\simeq {\cal
E}(X'_\alpha,X'_\gamma)
$$ of \ref{extcoc}. Then $(\{{\cal
E}(X'_\alpha,X'_\beta)\},\{F_{X'_\alpha,X'_\beta,X'_\gamma}\})$
is an extension cocycle of  $\Omega_{X_0/\kappa}$ by
$\a\otimes_\kappa{\cal O}_{X_0}$, which we will denote simply by
$\{{\cal E}(X'_\alpha,X'_\beta)\}$.

Observe that an abstract lifting exists if and only if it is
possible to choose the
$X'_\alpha$ so that there are isomorphisms of abstract liftings
$X'_\beta\rest{X_\alpha \cap X_\beta} \simeq
X'_\alpha\rest{X_\alpha \cap X_\beta}$ satisfying the cocycle
condition. By
\ref{splisom} to give isomorphisms of abstract liftings
$X'_\beta\rest{X_\alpha \cap X_\beta}
\simeq X'_\alpha\rest{X_\alpha \cap X_\beta}$ is equivalent to
assigning a collection of splittings
${\cal E}(X'_\alpha,X'_\beta)\simeq {\bf 0}$; it is
straightforward to check that the cocycle condition corresponds
to asking that the splittings yield an isomorphism of the
extension cocycle $\{{\cal E}(X'_\alpha,X'_\beta)\}$ with the
trivial cocycle
$\{{\bf 0}\}$. In other words, an abstract lifting exists if and
only if it is possible to choose a collection of abstract
liftings
$\{X'_\alpha\}$ in such a way that the associated extension
cocycle is trivial.

If $\{\widetilde X'_\alpha\}$ is another collection of abstract
liftings, then by
\ref{extcoc} and
\ref{part extopp}, and
\ref{sumextcomm}, we get isomorphisms
$$ {\cal E}(\widetilde X'_\alpha,\widetilde X'_\beta) 
\simeq {\cal E}(\widetilde X'_\alpha,X'_\alpha) + {\cal
E}(X'_\alpha,X'_\beta) + {\cal E}(X'_\beta,\widetilde X'_\alpha)
\simeq {\cal E}(X'_\alpha,X'_\beta) + {\cal E}(\widetilde
X'_\alpha,X'_\alpha) - {\cal E}(\widetilde X'_\beta,X'_\beta).
$$ One checks that these give an isomorphism of the cocycle
$\{{\cal E}(\widetilde X'_\alpha,\widetilde X'_\beta)\}$ with
$\{{\cal E}(X'_\alpha,X'_\beta)\}+\partial\{{\cal E}(\widetilde
X'_\alpha,X'_\alpha)\}$. This means that the class
$$\omega = [\{{\cal E}(X'_\alpha,X'_\beta)\}]\in \E_{{\cal
O}_{X_0}}(\Omega_{X_0/\kappa},\a\otimes_\kappa{\cal O}_{X_0})$$
is independent of the liftings. What's more, if $\{{\cal
E}_{\alpha\beta}\}$ is a cocycle in $\omega$, then there exist
extensions ${\cal E}_\alpha$ such that $\{{\cal
E}_{\alpha\beta}\}$ is isomorphic to $\{{\cal
E}(X'_\alpha,X'_\beta)\}+\partial\{{\cal E}_\alpha\}$. If we
choose liftings
$\widetilde X'_\alpha$ with isomorphisms ${\cal E}_\alpha \simeq
{\cal E}(\widetilde X'_\alpha,X'_\alpha)$ (\ref{exttrans}) then
the cocycle $\{{\cal E}_{\alpha\beta}\}$ is isomorphic to
$\{{\cal E}(\widetilde X'_\alpha,\widetilde X'_\beta)\}$. So it
follows that $\omega = 0$ if and only if a lifting exists. By
\ref{extcocgrp} the proof is concluded.

This proof could be summarized as follows, for those whose taste
runs towards the abstract. {}From the proof of \ref{absdefc} we
get that the category of liftings is a principal bundle stack
over the commutative group stack ${\goth Ext}_{{\cal
O}_{X_0}}(\Omega_{X_0,\kappa},\a\otimes_\kappa {\cal O}_{X_0})$
of extensions. An abstract lifting is a section of this bundle
stack, so the obstruction to its existence is an element of the
first cohomology group of ${\goth Ext}_{{\cal
O}_{X_0}}(\Omega_{X_0,\kappa},\a\otimes_\kappa {\cal O}_{X_0})$.
But
\ref{extcocgrp} says that this group is contained in
$\ext^2_{{\cal O}_{X_0}}(\Omega_{X_0,\kappa},\a\otimes_\kappa
{\cal O}_{X_0})$

The class constructed here and the class constructed earlier
under
\ref{simphyp} coincide. The proof of this is omitted.

This construction has an obvious property of functoriality,
which will be exploited in \refs{formal} and \refn{versal}. Let
$B'$ be a local ring,
${\goth b}\subseteq B'$ be an ideal with $\m_{B'}{\goth b} = 0$,
$B = B'/{\goth b}$. Let
$f\colon A'\to B'$ a local homomorphism inducing an isomorphism
of residue fields, such that
$f(\a)\subseteq {\goth b}$.  Set $f_*X = X\times_{\spec A} \spec
B$; this scheme $f_*X$ is a flat \loc{} on $A$. If
$X'$ is an abstract lifting of $X$ we set $f_*X' =
X\times_{\spec A'}
\spec B'$; this is an abstract lifting of $f_*X$. Any
automorphism
$\phi$ of $X'$ as a lifting induces an automorphism $f_*\phi$ of
$f_*X'$.

Call $g = f\rest\a \colon \a\to {\goth b}$ the restriction of
$f$.

\proclaim Proposition. \punto \call{absautofunc} If $\phi$ is an
automorphism of an abstract lifting $X'$ corresponding to an
element $\xi\in
\a\otimes_\kappa\hom_{{\cal
O}_{X_0}}(\Omega_{{X_0}/\kappa},{\cal O}_{X_0})$, the element of
${\goth b}\otimes_\kappa\hom_{{\cal
O}_{X_0}}(\Omega_{{X_0}/\kappa},{\cal O}_{X_0})$ corresponding
to $f_*\phi$ is $(g\otimes{\rm id})(\xi)$.

\punto \call{extfunc} $e(f_*X'_1,f_*X'_2) = (g\otimes {\rm
id})\bigl(e(X'_1,X'_2)\bigr)\in {\goth b}\otimes \ext^1_{{\cal
O}_{X_0}}(\Omega_{X_0/\kappa},{\cal O}_{X_0})$.

\punto \call{absobsfunc} $\omega_{\rm abs}(f_*X') = (g\otimes
{\rm id})\omega_{\rm abs}(X)\in {\goth b}\otimes \ext^2_{{\cal
O}_{X_0}}(\Omega_{X_0/\kappa},{\cal O}_{X_0})$.
\endproclaim

\beginsection Generalizations [generalizations]

Here are two important generalizations.

Fix a scheme $M'$ of finite type over $A'$, a scheme $X$ of
finite type over $A$, such that $X_0$ is reduced and generically
smooth over $\kappa$, and a \loc{} morphism $f\colon X\to M =
M'\rest{\spec A}$ defined over $A$. A {\it lifting\/} of the
morphism $f$ consist of an abstract lifting $X'$ of $X$, and
morphism $f'\colon X'\to M'$ of $A'$-schemes whose restriction
to $X$ is $f$. There is an obvious notion of isomorphism of
liftings; if $f'_1\colon X'_1\to M'$ and $f'_2\colon X'_2\to M'$
are liftings, then an isomorphism of $f_1$ with $f_2$ is an
isomorphism of abstract liftings $\Phi\colon X'_1\simeq X'_2$
such that $f'_2\circ \Phi = f'_1$.

Let $f_0\colon X_0\to M_0$ be the restriction of $f$; $f_0$ is
again a
\loc{} morphism. We define the {\it complex of differentials\/}
$\Omega^\mini_{f_0}$ of the morphism $f_0$ to be the complex with
$\Omega^0_{f_0} =
\Omega_{X_0/\kappa}$, $\Omega^{-1}_{f_0} =
f^*\Omega_{M_0/\kappa}$,
$\Omega^i_{f_0} = 0$ for $i\neq 0$,~$-1$, the only nontrivial
differential
$f^*\Omega_{M_0/\kappa}\to
\Omega_{X_0/\kappa}$ being the pullback map.

\proclaim Theorem. \call{mapsdef}\punto \call{mapsdefa} Any
lifting of
$f$ is a \loc{} morphism.

\punto \call{mapsdefb} There is a canonical element $\omega\in
\a\otimes_\kappa\ext^2_{{\cal O}_{X_0}}(\Omega^\mini_{f_0},{\cal
O}_{X_0})$, called the\/ {\rm obstruction}, such that $\omega =
0$ if and only if a lifting exists.

\punto \call{mapsdefc} If a lifting exists, then there is a
canonical action of the group
$\a\otimes_\kappa\ext^1_{{\cal O}_{X_0}}(\Omega^\mini_{f_0},{\cal
O}_{X_0})$ on the set of isomorphism classes of liftings making
it into a principal homogeneous space.
\endproclaim

The proof of this result is essentially the same as the proof of
\ref{absdef} in
\ref{abstract}.

First one assumes that there is a smooth morphism $\pi'\colon
P'\to M'$ such that
$f\colon X\to M$ factors through $P = P'\rest{\spec A}$. Let
$X'_1$ and
$X'_2$ be two liftings of $X$ in $P'$; the resulting morphisms
$f'_1\colon X'_1\to M'$ and
$f'_2\colon X'_1\to M'$ are both liftings of $f$. Here is a
description of the isomorphisms between $f'_1$ and $f'_2$
analogous to \ref{isomlif}, with the same proof.

Call ${\cal I}_0$ the ideal of the embedding of $X_0$ in $P_0$,
and
$\pi_0\colon P_0\to M_0$ the restriction of $\pi$. The usual
homomorphism from
${\cal I}_0/{\cal I}_0^2$ to $\Omega_{P_0/M_0}\rest{X_0}$ yields
a restriction map from
$$\hom_{{\cal O}_{X_0}}({\Omega_{P_0/M_0},\a\otimes_\kappa{\cal
O}_{X_0}})\longrightarrow \H^0(X_0,\a\otimes{\cal N}_0).$$

\proclaim Proposition. There is a bijective correspondence
between isomorphisms of liftings
$f'_2\simeq f'_1$ and elements of\/
$\hom_{{\cal
O}_{X_0}}({\Omega_{P_0/M_0}\rest{X_0},\a\otimes_\kappa{\cal
O}_{X_0}})$ whose image in
$\H^0(X_0,\a\otimes{\cal N}_0)$ is $\nu(X'_1,X'_2)$, with the
following properties.

\punto \call{isomid} The identity ${\rm id}_{X'}\colon X'\simeq
X'$ corresponds to 0.

\punto \call{isomcoc} If $\Phi_1\colon X'_2\simeq X'_1$ and
$\Phi_2\colon X'_3\simeq X'_2$ are isomorphisms of abstract
liftings and $D_1$, $D_2$ are the corresponding elements of
$\hom_{{\cal
O}_{X_0}}({\Omega_{P_0/\kappa}\rest{X_0},\a\otimes_\kappa{\cal
O}_{X_0}})$, then the composition $\Phi_1\circ \Phi_2$
corresponds to $D_1+D_2$.
\endproclaim

The role of the fundamental exact sequence of \ref{exseq} is
played by the exact sequence
$$
\sequence{0&\mapright&{\cal I}_0/{\cal
I}_0^2&\mapright&\Omega^\mini_{\pi_0}\rest{X_0}&\mapright
&\Omega^\mini_{f_0}&\mapright&0}
$$ which looks as follows
$$
\def\grado#1{\hbox to 45pt{$\deg\ #1$\hfill}}
\diagram{30}{&&&\vdots&&\vdots&&\vdots\cr
&&&\mapdown&&\mapdown&&\mapdown\cr
\grado{-2}&0&\mapright&0&\mapright&0&\mapright&0&\mapright&0\cr
&&&\mapdown&&\mapdown&&\mapdown\cr
\grado{-1}&0&\mapright&0&\mapright&\pi_0^*\Omega_{M_0/\kappa}
\rest{X_0}&\bihline& f_0^*\Omega_{M_0/\kappa} &\mapright&0\cr
&&&\mapdown&&\mapdown&&\mapdown\lft\partial\cr
\grado{0}&0&\mapright&{\cal I}_0/{\cal I}_0^2 &\mapright&
\Omega_{P_0/\kappa}\rest{X_0} &\mapright
&\Omega_{X_0/\kappa}&\mapright&0\cr
&&&\mapdown&&\mapdown&&\mapdown\cr
\grado{1}&0&\mapright&0&\mapright&0&\mapright&0&\mapright&0\cr
&&&\mapdown&&\mapdown&&\mapdown\cr &&&\vdots&&\vdots&&\vdots\cr}
$$

If we observe that for any sheaf of ${\cal O}_{X_0}$-modules
$\cal G$ we have
$$\hom_{{\cal O}_{X_0}}(\Omega_{P_0/M_0}\rest{X_0},{\cal G}) =
\hom_{{\cal O}_{X_0}}(\Omega^\mini_{\pi_0}\rest{X_0},{\cal G}),$$
and we keep into account that, as observed in the last part of
\ref{extensions}, the theory of extensions of sheaves
generalizes to extensions of
$\Omega^\mini_{f_0}$ by sheaves, the proof of
\ref{absdef} goes through almost word for word.

Another situation that arises quite often in practice is when we
want to study deformations of a scheme inducing a fixed
deformation on a subscheme. For example, in the theory of
deformation of pointed curves the scheme is the curve itself,
while the subscheme is the union of the distinguished points.
Again we may look at embedded deformations or abstract
deformations. Here is the embedded setup.

Let $M'$ be a flat scheme of finite type over
$A'$, $Z'\subseteq M'$ a closed subscheme, also flat over $A'$.
Let
$X$ be a \loc{} subscheme of
$M = M'\rest{\spec A}$ containing $Z = Z'\rest{\spec A}$. Assume
also that $X$ is flat over
$A$. Call ${\cal N}_0$ the normal sheaf to $X_0$ in $M_0$,
${\cal J}_0$ the ideal of $Z_0$ in
$X_0$. A {\it lifting of $X$ to $M'$ relative to $Z'$} is a
lifting $X'\subseteq M'$ that contains $Z'$.

\proclaim Theorem. There is a canonical element
$$
\omega_{\rm emb}\in
\a\otimes_\kappa\H^1(X_0,{\cal N}_0\otimes_{{\cal O}_{X_0}}{\cal
J}_0)
$$ called the\/ {\rm embedded obstruction} of
$X$ in
$M$, such that $\omega_{\rm emb} = 0$ if and only if a lifting
of $X$ to $M'$ relative to
$Z'$ exists.

\punto If a lifting exists, then there is a canonical action of
the group
$$\a\otimes_\kappa\H^0(X_0,{\cal N}_0\otimes_{{\cal
O}_{X_0}}{\cal J}_0)$$ on the set of liftings making it into a
principal homogeneous space.
\endproclaim

The proof is essentially the same as the the proof of
\ref{embdef}. The key point is that if we have two relative
liftings $X'_1$ and $X'_2$, then the image of the homomorphism
$\nu(X'_1,X'_2)\colon{\cal I}_0/{\cal I}_0 \to \a\otimes_\kappa
{\cal O}_{X_0}$ is contained inside $\a\otimes_\kappa{\cal
J}_0$, so $\nu(X'_1,X'_2)$ can be considered as an element of
$$\hom_{{\cal O}_{X_0}}({\cal I}_0/{\cal
I}_0,\a\otimes_\kappa{\cal J}_0) =
\a\otimes_\kappa\H^0(X_0,{\cal N}_0\otimes_{{\cal O}_{X_0}}{\cal
J}_0).$$

Here is the setting in the abstract case. 

Let $Z'$ be a flat scheme over $A'$, and let set $Z =
Z'\rest{\spec A}$. Let $X$ be a flat \loc{} scheme of finite
type over $A$ with a closed embedding $Z\into X$. Assume also
that $X_0$ is generically smooth over $\kappa$.

An {\it abstract lifting of $X$ relative to $Z'$} is an abstract
lifting $X'$ of $X$ with a closed embedding $Z'\into X'$
extending the given embedding of $Z$ in
$X$.

An isomorphism of relative abstract liftings is an isomorphism of
abstract lifting inducing the identity on $Z'$.

Again, let ${\cal J}_0$ be the ideal of $Z_0$ in $X_0$.

\proclaim Theorem. \punto There is a canonical element
$\omega_{\rm abs}\in
\a\otimes_\kappa\ext^2_{{\cal
O}_{X_0}}(\Omega_{{X_0}/\kappa},{\cal J}_0)$, called the\/ {\rm
obstruction}, such that $\omega_{\rm abs} = 0$ if and only if an
abstract lifting exists.

\punto \call{relabsdefc} If an abstract lifting exists, then
there is a canonical action of the group
$$\a\otimes_\kappa\ext^1_{{\cal
O}_{X_0}}(\Omega_{{X_0}/\kappa},{\cal J}_0)$$ on the set of
isomorphism classes of abstract liftings making it into a
principal homogeneous space.
\endproclaim

The proof of \ref{absdef} goes through with minor variations,
substituting $\a\otimes_\kappa {\cal J}_0$ everywhere for
$\a\otimes_\kappa {\cal O}_{X_0}$. The only point that requires a
little care is the analogue of \ref{locP'ex}. Given two relative
abstract liftings $Z'\into X'_1$ and
$Z'\into X'_2$, to construct the relative extension we need to
compare them locally in a suitable ambient space. For this we
need to know that if
$X$ is affine, given a closed embedding $X\into {\bf A}^n_A$,
this can be lifted to embeddings
$X'_1\into {\bf A}^n_{A'}$ and $X'_2\into {\bf A}^n_{A'}$ which
agree on $Z'$. For this, set
$X = \spec B$, $X'_i = \spec B'_i$, $Z = \spec C$, $Z' = \spec
C'$. The embeddings $Z'\into X'_i$ and $X\into X'_i$, which
agree on $Z$, induce a homomorphism of $A'$-algebras $B'_i\to
B\times_C C'$, which is easily shown to be surjective. Then the
given embedding $X\into {\bf A}^n_A$ yields a surjective
homomorphism of $A$-algebras $A[x_1,\ldots,x_n]\to B$. Lift the
images of the $x_i$ to $B\times_C C'$, then lift them to $B'_1$
and $B'_2$. The resulting homomorphisms of $A'$-algebras
$A'[x_1,\ldots,x_n]\to B'_i$ are surjective, and induce the
required liftings
$X'_i\into {\bf A}^n_{A'}$.

Of course we can put together the two generalizations. This is
our final and most general setup. It is not in any essential
sense harder to treat than the case of \ref{abshyp}, just a
little more confusing.

Let $X$ be a \loc{} flat scheme of finite type over $A$ such that
$X_0$ is generically smooth over
$\kappa$, $M'$ and
$Z'$ flat schemes over
$A'$ with a morphism $j'\colon Z'\to M'$. Let it be given a
closed embedding of $Z = Z'\rest{\spec A}$ into $X =
X'\rest{\spec A}$, and morphism $f\colon X\to M = M'\rest{\spec
A}$ such that $f\rest Z = j'\rest Z$. A {\it lifting of $f$
relative to $j'$} is an abstract lifting $Z'\into X'$ relative
to $Z'$, together with a morphism $f'\colon X'\to M'$ such that
$f' \rest X = f$ and
$f'\rest{Z'} = j'$. Call
${\cal J}_0$ the ideal of $Z_0$ in $X_0$.

The notion of isomorphism of relative liftings of $f$ is the
obvious one.

\proclaim Theorem. \punto There is a canonical element $\omega\in
\a\otimes_\kappa\ext^2_{{\cal O}_{X_0}}(\Omega^\mini_{f_0},{\cal
J}_0)$, called the\/ {\rm obstruction}, such that $\omega = 0$
if and only if a lifting exists.

\punto  If a lifting exists, then there is a canonical action of
the group
$\a\otimes_\kappa\ext^1_{{\cal O}_{X_0}}(\Omega^\mini_{f_0},{\cal
J}_0)$ on the set of isomorphism classes of liftings making it
into a principal homogeneous space.
\endproclaim

Finally, all the results and the proofs generalize to the case of
algebraic spaces, by working with the \'etale topology instead
of the Zariski topology. Also, if we assume that
$A'$ is a finite {\bf C}-algebra they are valid for analytic
spaces; the proofs remain the same, if we substituting polydiscs
for affine spaces.

\beginsection Formal deformations [formal]

We fix a field $\kappa$. In this section and in the next a {\it
complete algebra\/} will be a complete noetherian local
$\kappa$-algebra with residue field $\kappa$. If $A$ is a
complete algebra, we will denote its maximal ideal by $\m_A$,
and set $A_n = A/\m_A^{n+1}$. An {\it artinian algebra\/} will
be a complete algebra which is artinian, or, equivalently,
finite over
$\kappa$. If $A$ is a complete algebra, the $\kappa$-vector
space $\cot A$ is called the {\it cotangent space} of
$A$, the dual space $\tang A$ its {\it tangent space}.

A homomorphism of complete algebras is homomorphisms of
$\kappa$-algebras; it is automatically local. A homomorphism of
complete algebras $f\colon A\to B$ induces a linear map
$f_*\colon
\cot A\to \cot B$; the dual map ${\rm d}f\colon \tang B\to \tang
A$ will be called the {\it differential\/} of $f$. If $f\colon
A\to B$ is a homomorphism of complete algebras then
$f(\m_A^{n+1})\subseteq \m_B^{n+1}$ for each $n\ge 0$; we denote
by
$f_n\colon A_n\to B_n$ the induced homomorphism.

\proclaim Definition. A\/ {\rm deformation} $(X,A)$ consists of a
complete algebra $A$, a scheme $X_n$ flat and of finite type
over $A_n$ for each $n\ge 0$, and a sequence of closed embeddings
$X_{n-1}\into X_n$ compatible with with the closed embeddings
$\spec{A_{n-1}}\into \spec{A_n}$, inducing an isomorphism of
$X_n\rest{\spec{A_{n-1}}}$ with
$X_{n-1}$.

We will say that $X$ is a\/ {\rm deformation of $X_0$ over $A$}. 

If $X$ and $\widetilde X$ are deformations of $X_0 = \widetilde
X_0$ over $A$, an {\rm isomorphism
$\phi\colon X\simeq \widetilde X$ of deformations over $A$} is a
sequence of isomorphisms
$\phi_n\colon X_n\simeq \widetilde X_n$ of schemes over $A_n$,
such that
$\phi_n\rest{X_{n-1}} =
\phi_{n-1}$, and $\phi_0\colon X_0\to X_0$ is the identity.
\endproclaim

The objects defined above should be properly called {\it formal
deformations}, but they are the only types of deformations we
will consider.

{}From now on we fix $X_0$; all deformations will be
deformations of the same $X_0$.

Deformations of $X_0$ over $A$ form a category, the arrows being
isomorphisms of deformations over $A$.

Let $f\colon A\to B$ be a homomorphism of complete algebras, $X$
a deformation of $X_0$ over
$A$. There is an induced deformation $f_*X$ on $B$, defined by
setting
$(f_*X)_n = X_n\times_{\spec{A_n}}
\spec{B_n}$; the embeddings
$(f_*X)_n\into (f_*X)_{n-1}$ are induced by the embeddings
$X_n\into X_{n-1}$. Also, if
$\phi\colon X\simeq \widetilde X$ is an isomorphism of
deformations of
$X_0$, then there is an induced isomorphism $f_*\phi\colon
f_*X\simeq f_*\widetilde X$, defined in the obvious way. This
makes $f_*$ into a functor from the category of deformations of
$X_0$ over $A$ to deformations of $X_0$ over $B$.

If $f\colon A\to B$ and $g\colon B\to C$ are homomorphism of
complete algebras and $X$ is a deformation over $A$, then there
is a canonical isomorphism
$(gf)_*X\simeq g_*f_*X$ of deformations over $C$. {}From now on
we identify $(gf)_*X$ with
$g_*f_*X$.

The {\it trivial deformation\/} of $X_0$ over a complete algebra
$A$ is $X_0^A = f_*X_0$, where $f\colon \kappa\to A$ is the
structure homomorphism. Concretely,
$\bigl(X_0^A\bigr)_n = X_0\times_\kappa{\spec{A_n}}$, and the
closed embeddings
$\bigl(X_0^A\bigr)_{n-1}\into
\bigl(X_0^A\bigr)_n$ are induced by the closed embeddings
$\spec{A_{n-1}}\into\spec{A_n}$

\proclaim Definition. Let $(X,A)$, $(Y,B)$ be deformations of
$X_0$. An homomorphism
$(\phi,f)\colon (X,A)\to (Y,B)$ of deformations consists of a
homomorphism of complete algebras $f\colon A\to B$ and an
isomorphism $\phi\colon f_*X\simeq Y$ of deformations over $B$.

A homomorphism $(\phi,f)$ is called\/ {\rm surjective} if $f$ is
surjective.

If $(\phi,f)\colon (X,A)\to (Y,B)$ and $(\psi,g)\colon (Y,B)\to
(Z,C)$ are homomorphism of deformations, the composition
$(\psi,g)\circ (\phi,f)\colon (X,A)\to (Z,C)$ is $(\psi\circ
g_*\phi, g\circ f)$.
\endproclaim

Extensions of $X_0$ are the objects of a category in which the
arrows are the homomorphisms. The object $(X_0,\kappa)$ is
terminal in this category, that is, given a deformation $(X,A)$
there is a unique homomorphism $(X,A)\to (X_0,\kappa)$.

The isomorphisms in the category of deformations are exactly the
homomorphisms $(\phi,f)$ where $f$ is an isomorphism.

If $A$ is artinian then
$A_n = A$ for $n\gg 0$, so $X_n$ is a flat scheme over $A$ for
$n\gg 0$. As one sees immediately, the category of deformations
of $X_0$ over $A$ is equivalent to the category of flat schemes
$X$ of finite type over $A$, together with closed embeddings of
$X_0\subseteq X$ inducing isomorphisms $X_0\simeq X\rest{\spec
\kappa}$; the arrows are isomorphisms of $A$-schemes inducing
the identity on $X_0$. {}From now on we will systematically
identify a deformation over an artinian algebra $A$ with the
corresponding scheme over $A$.

More generally, if $X$ is a flat scheme of finite type over $A$
with
$X\rest{\spec\kappa} = X_0$, then $X$ will induce a deformation
of $X_0$ over $A$ by setting
$X_n = X\rest{\spec{A_n}}$, the embeddings $X_{n-1}\into X_n$
being induced by the natural embeddings
$\spec A_{n-1}\into \spec A_n$. Such a deformation is called {\it
algebraic}.

Not every deformation is algebraic. For example, one can show
that if
$X_0\subseteq {\bf P}_{\bf C}^3$ is a smooth quartic surface,
then
$X_0$ has non-algebraic deformations over ${\bf C}[[t]]$.

This problem does not arises for projective curves. More
generally we have the following standard fact.

\proclaim Proposition. \call{algdeformation}Assume that $X_0$ is
projective, and
$\H^2(X_0,{\cal O}_{X_0}) = 0$. Then every deformation of $X_0$
is algebraic.
\endproclaim

\proof Let $X$ be a deformation of $X_0$. I claim that there for
each
$n>0$ the restriction map $\pic{X_n}\to \pic{X_{n-1}}$ is
surjective. In fact, there is an exact sequence
$$
\sequence{0& \mapright & (\m_A^n/\m_A^{n+1})\otimes_{\kappa}{\cal
O}_{X_0}&
\mapright^\alpha& {\cal O}_{X_n}^*& \mapright^\beta& {\cal
O}_{X_{n-1}}^*& \mapright 0}
$$ in which $\beta$ is the restriction map, and $\alpha$ is
defined by identifying
$$(\m_A^n/\m_A^{n+1})\otimes_{\kappa}{\cal O}_{X_0} =
(\m_A^n/\m_A^{n+1})\otimes_{A/\m_A^n}{\cal O}_{X_n}$$ with the
kernel of the restriction map
${\cal O}_{X_n}\to {\cal O}_{X_{n-1}}$, which can be done
because of the flatness of
$X_n$, then setting $\alpha(f) = 1+f$. The fact that
$\H^2\bigl(X_0,(\m_A^n/\m_A^{n+1})\otimes_{\kappa}{\cal
O}_{X_0}) = 0$ implies the surjectivity of the map on Picard
groups.

Let ${\cal L}_0$ be a very ample line bundle on $X_0$ such that
$\H^i(X_0,{\cal L}_0) = 0$ for $i>0$. For each $n>0$ we can
choose  a line bundle ${\cal L}_n$ on
$X_n$ such that ${\cal L}_n\rest{X_{n-1}}$ is isomorphic to
${\cal L}_{n-1}$; by semicontinuity we have that $\H^i(X_n,{\cal
L}_n) = 0$ for $i>0$. If
$\pi_n\colon X_n\to
\spec{A_n}$ is the structure morphism then ${\pi_n}_*{\cal L}_n$
satisfies base change, and if $N$ is the dimension of
$\H^0(X_0,{\cal L}_0)$ over
$\kappa$ the bundles
${\cal L}_n$ induce embeddings $X_n\into {\bf P}^{N-1}_{A_n}$,
such that $X_n\cap {\bf P}^{N-1}_{A_{n-1}} = X_{n-1}$. So the
system $\{X_n\}$ can be considered as a formal subscheme of
${\bf P}^{N-1}_A$, and the result follow from Grothendieck's
existence theorem.\endproof

I do not know whether this is still true if we do not assume that
$X_0$ is projective. Elkik proved that deformations of affine
schemes with isolated singularities are algebraic ([Elkik]).

{}From now on we will assume that $X_0$ is a generally smooth
\loc{} scheme on $\kappa$. We set
$$ {\rm T}^i(X_0) = \ext^i_{{\cal
O}_{X_0}}(\Omega_{X_0/\kappa},{\cal O}_{X_0}).
$$

If $X$ is a deformation, then $X_n$ is an abstract lifting of
$X_{n-1}$ to $A_n$, so a deformation of $X_0$ can be thought of
a sequence of abstract liftings; to these we can apply
\ref{absdef}.

\proclaim Proposition. \call{deft}\punto \call{deft0}Assume that
${\rm T}^0(X_0) = 0$. Then two deformations of
$X_0$ over the same algebra admit at most one isomorphism.

\punto \call{deft1}Assume that ${\rm T}^1(X_0) = 0$. Then any
deformation $(X,A)$ of
$X_0$ is isomorphic to the trivial deformation $X_0^A$.

\punto \call{deft2} Assume that ${\rm T}^2(X_0) = 0$. Then if
$(X,A)$ is a deformation of
$X_0$ and $f\colon B\to A$ is a surjective homomorphism of
algebras, there exists a deformation $Y$ of $X_0$ over $B$ and
an isomorphism $f_*Y\simeq X$.
\endproclaim

\proof For part \ref{part deft0}, notice that because of
\ref{absdefd} an isomorphism
$X_{n-1}\simeq \widetilde X_{n-1}$ over $A_n$ extends in at most
one way to an isomorphism
$X_n\simeq \widetilde X_n$.

For part \ref{part deft1}, let $\widetilde X = X_0^A$ let
$\phi_{n-1}\colon\widetilde X_{n-1}\simeq X_{n-1}$ be an
isomorphism inducing the identity on
$X_0$. Then $X_n$ is a lifting of
$X_{n-1}$ to
$A_n$, and we can also think of
$\widetilde X_n$ as a lifting, via the composition $\widetilde
X_{n-1}\simeq X_{n-1}\into X_n$. Then it follows from
\ref{absdefc} that the isomorphism
$\phi_{n-1}$ extends to an isomorphism $\phi_n\colon \widetilde
X_n\simeq X_n$.

Let us prove part \ref{part deft2}.  Call {\goth b} the kernel
of $f$.

Assume first that $A$ and $B$ are artinian, so a deformation $X$
on
$A$ is flat scheme over
$A$. By induction on the least integer $n$ such that
$\m_B^n{\goth b} = 0$ we can assume that $\m_B{\goth b} = 0$;
then the result follows from \ref{absdefb}.

In the general case we construct $Y_n$ by induction on $n$. Call
$\pi_n\colon B_n\to A_n = B/({\goth b}+\m_B^{n+1})$ the
projection. For
$n=0$ there is no problem, so suppose
$n>0$ and that we are given  a deformation $Y_{n-1}$ over
$B_{n-1} = A/(\a+\m_A^n)$ and a homomorphism
$(\phi_{n-1},\pi_{n-1})\colon Y_{n-1}\to X_{n-1}$. We are
looking for a deformation $(Y_n,B_n)$ and a commutative diagram
$$
\diagram{60}{(Y_n,B_n)&\mapright^{(\phi_n,\pi_n)}&(X_n,A_n)\cr
\mapdown&&\mapdown\cr
(Y_{n-1},B_{n-1})&\mapright^{(\phi_{n-1},\pi_{n-1})}&(X_{n-1},A_{n-1}).\cr
}
$$ {}From \ref{fiberproddef} we get a deformation $\widetilde Y$
over
$B/\bigl(\m_B^n\cap ({\goth b}+\m_B^{n+1}) \bigr)$ and a
commutative diagram of deformations
$$
\diagram{60}{\bigl(\widetilde Y,B/\bigl(\m_B^n\cap ({\goth
b}+\m_B^{n+1}) \bigr)\bigr)&
\mapright&(X_n,A_n)\cr
\mapdown&&\mapdown\cr (Y_{n-1},B_{n-1})&\mapright^{(\phi_{n-1},
\pi_{n-1})}&(X_{n-1},A_{n-1}).\cr }
$$ By the previous case there is a homomorphism of deformations
$(Y_n,A/\m_A^{n+1})\to\bigl(\widetilde Y,B/\bigl(\m_B^n\cap
({\goth b}+\m_B^{n+1})
\bigr)$; this is the deformations we were looking for.\endproof

It may happen that ${\rm T}^2(X_0)\neq 0$, but still the
conclusion of
\ref{deft2} holds. To clarify this we give a definition.

Let $(X,A)$ be a deformation of $X$ over an artinian algebra
$A$. A {\it small extension} of
$A$ is a surjective homomorphism of artinian algebras $A'\to A$
whose kernel $\a$ has length 1, and is therefore isomorphic to
$\kappa$. These data determine an element
$\omega_{\rm abs}\in {\rm T}^2(X_0)\simeq \a\otimes{\rm
T}^2(X_0)$, well defined up to multiplication by a nonzero
scalar.

\proclaim Definition. The\/ {\rm space of obstructions} $\obs
X_0$ of
$X_0$ is the subspace of ${\rm T}^2(X_0)$ generated by the
elements $\omega_{\rm abs}\in {\rm T}^2(X_0)$ for all
deformations $(X,A)$ and all small extensions $A'\to A$ as above.

If $\obs X_0 = 0$ we say that $X_0$ {\rm has unobstructed
deformations}.
\endproclaim

\proclaim Proposition. Let $A'$ be an artinian algebra,
$\a\subseteq A'$ an ideal with
$\m_{A'}\a = 0$, $A=A'/\a$. Let $X$ be a deformation of $X_0$
over
$A$; then the obstruction
$\omega_{\rm abs}(X)\in \a\otimes{\rm T}^2(X_0)$ is contained in
$\a\otimes\obs X_0$.
\endproclaim

\proof Choose a basis $v_1,\dots,v_n$ for $\a$ as a vector space
over
$\kappa$, and for each $i$ call $\a_i$ the quotient of $\a$ by
the subspace generated by
$v_1,\dots,v_{i-1}$,
$v_{i+1},\dots,v_n$. Set $A'_i = A/\a_i$; then $A'_i\to A$ is a
small extension, and by definition the obstruction to lifting
$X$ to $A'_i$ lies in
$(\a/\a_i)\otimes\obs X_0\subseteq (\a/\a_i)\otimes{\rm
T}^2(X_0)$. According to
\ref{absobsfunc} this obstruction is the image in
$\a/\a_i\otimes{\rm T}^2(X_0)$ of
$\omega_{\rm abs}\in \a\otimes{\rm T}^2(X_0)$, and by elementary
linear algebras this implies the thesis.\endproof

So we can improve on \ref{deft2}.

\proclaim Proposition. \call{deft2imp} Assume that $\obs X_0 =
0$. Then if $(X,A)$ is a deformation of
$X_0$ and $f\colon B\to A$ is a surjective homomorphism of
algebras, there exists a deformation $Y$ of $X_0$ over $B$ and
an isomorphism $f_*Y\simeq X$.
\endproclaim

The definition of unobstructed variety is relative to then base
field; I do not know any specific example, but it seems quite
plausible that there may be varieties over a field
$\kappa$ of positive characteristic which are unobstructed, but
cannot be lifted to some artinian ring $A$ with residue field
$\kappa$ (obviously $A$ can not be a $\kappa$-algebra.) On the
other hand one can prove that in the situation of
\ref{abstract} if the ring $A'$ is equicharacteristic then the
obstruction to lifting $X$ to $A'$ lives in $\a\otimes\obs X_0$.

It happens very frequently that ${\rm T}^2(X_0)\neq 0$ but
$\mathop{\rm Obs}{X_0} = 0$.

\proclaim Example. Let $X_0\subseteq {\bf P}^3_\kappa$ be a
smooth surface of degree $d\ge 6$. Then ${\rm T}^2(X_0)\neq 0$
but $\mathop{\rm Obs}{X_0} = 0$.
\endproclaim

\proof Denote by ${\cal T}_{X_0}$ the tangent bundle of $X_0$;
then
${\rm T}^i(X_0) =
\H^i(X_0,{\cal T}_{X_0})$. We have the sequence
$$
\harrowlength=\sequencelength \varrowlength=10pt
\commdiag{0&\mapright&{\cal T}_{X_0}& \mapright& {\cal T}_{{\bf
P}^3}\rest {X_0}& \mapright& {\cal N}_0& \mapright& 0\cr
&&&&&&{\bivline}\cr &&&&&&{\cal O}_{X_0}(d).\cr}
$$

Since $\H^1(X_0,{\cal N}_0) = 0$ the embedded deformations of
$X_0$ in
${\bf P}^3$ are unobstructed. Assume that we have proved that
$\H^1(X_0,{\cal T}_{{\bf P}^3}\rest {X_0}) = 0$. Then  the
boundary map
$$\partial\colon \H^0(X_0,{\cal N}_0)\to \H^1(X_0,{\cal
T}_{X_0})$$ is surjective: I claim that this implies that for
any deformation $X$ of $X_0$ over an artinian algebra
$A$ there is an embedding $X\into {\bf P}^3_A$ extending the
given embedding $X_0\into {\bf P}^3_\kappa$. Let $\ell$ be the
length of $A$; for $\ell = 1$ we have
$A = \kappa$, and the statement is vacuous. Suppose this true
when $A$ has length $\ell-1$, and take an ideal
$\a\subseteq A$ of length 1. The deformation $Y = X\times_{\spec
A}\spec{A/\a}$ can be embedded in ${\bf P}^3_{A/\a}$ by
inductive hypothesis; since
$\H^1(X_0,{\cal N}_0) = 0$ then
$Y$ has a lifting $Y'$ in ${\bf P}^3_A$. Choose another lifting
$X'$ of $Y$ in ${\bf P}^3_A$ such that $\partial \nu(X',Y') =
e(X,Y')$ (\ref{nutrans}). {}From
\ref{nuext} we see that
$e(X',Y') = e(X,Y')$, so that $e(X',X) = e(X',Y') - e(X,Y') = 0$
(\ref{epsilonadd}). Then
$X'$ and
$X$ are isomorphic (\ref{epsilonisom}), and $X$ can be embedded
in
${\bf P}^3_A$.

Since $\H^1(X_0,{\cal N}_0) = 0$ this proves that the
deformations of
$X_0$ are unobstructed, hence $\obs X_0 = 0$.

By Serre duality we have $\H^1(X_0,{\cal T}_{{\bf P}^3}\rest
{X_0})\simeq \H^1(X_0,
\Omega_{{\bf P}^3}(d-4)\rest {X_0})^\vee$. {}From the twisted
restricted Euler sequence
$$
\sequence{0& \mapright& \Omega_{{\bf P}^3}(d-4)\rest {X_0}&
\mapright {\cal O}_{X_0}(d-5)^4& \mapright& {\cal
O}_{X_0}(d-4)&\mapright& 0,\cr }
$$ the easily proved surjectivity of the induced map
$\H^0\bigl(X_0, {\cal O}_{X_0}(d-5)\bigr)^4\to
\H^0\bigl(X_0, {\cal O}_{X_0}(d-4)\bigr)$ for $d\ge 5$, and the
fact that
$\H^1\bigl({\cal O}_{X_0}(d-5)\bigr) = 0$ we see that 
$\H^1(X_0,{\cal T}_{{\bf P}^3}\rest {X_0}) = 0$, as claimed.

Finally, by Serre duality $\H^2(X_0,{\cal T}_{{\bf P}^3}\rest
{X_0})
\simeq \H^0(X_0,
\Omega_{{\bf P}^3}(d-4)\rest {X_0})^\vee$. Again from the Euler
sequence we see that is enough to prove that
$$ 4{d-2\choose 3} = \dim_\kappa\H^0\bigl(X_0, {\cal
O}_{X_0}(d-5)\bigr)^4 >
\dim_\kappa
\H^0\bigl(X_0, {\cal O}_{X_0}(d-4)\bigr) = {d-1\choose 3}
$$ for $d\ge 6$, and this is straightforward.\endproof

Let us go back to the general situation. Let $(X,A)$ be a
deformation of $X_0$. Then $X_1$ is a lifting of $X_0$ to $A_1 =
A/\m_A^2$, so it can be compared to the trivial lifting
$X_0^{A_1} = X_0\times_{\spec \kappa}\spec {A_1}$.

\proclaim Definition. The\/ {\rm \KS{} class} of $X$ is
$$ {\rm k}_X = e(X_1,X_0^{A_1})\in (\cot A)\otimes_\kappa {\rm
T}^1(X_0)\simeq \hom(\tang A, {\rm T}^1(X_0)).
$$

The associated linear map ${\rm K}_X\colon \tang A\to {\rm
T}^1(X_0)$ is called the\/ {\rm
\KS{} map}.
\endproclaim

\KS{} classes and maps have an important functorial property.

\proclaim Proposition. Let $f\colon A\to B$ a of complete
algebras,
$X$ a deformation of
$X_0$ on $A$. Then:

\punto \call{ksclassfunc} ${\rm k}_{f_*X} = (f_*\otimes{\rm
id})({\rm k}_X)\in \cot B\otimes_\kappa {\rm T}^1(X_0)$;

\punto \call{ksmapfunc} ${\rm K}_{f_*X} = {\rm K}_X\circ {\rm
d}f\colon \tang B\longrightarrow {\rm T}^1(X_0)$.
\endproclaim

\proof The two statements are obviously equivalent; part
\ref{part ksclassfunc} follows from
\ref{extfunc}.\endproof

An alternate and more traditional description of the \KS{} map
is as follows.

Consider the algebra of dual numbers $\kappa[\epsilon] =
\kappa[x]/(x^2)$; call $X_0[\epsilon] = X_0\times_{\spec
\kappa}\spec
\kappa[\epsilon]$ the trivial deformation. Deformations of $X_0$
on $\kappa[\epsilon]$ are abstract liftings of $X_0$ to
$\kappa[\epsilon]$. If $X$ is such a lifting, consider
$$ e(X,X_0[\epsilon])\in (\epsilon)\otimes_\kappa{\rm T}^1(X_0)
= {\rm T}^1(X_0).
$$ We get a map from the set of liftings of $X_0$ to
$\kappa[\epsilon]$ to ${\rm T}^1(X_0)$.

If $A$ is a complete algebra and $f\colon A\to
\kappa[\epsilon]$ is a homomorphism of complete algebras, the
associated homomorphism of
$\kappa$-vector spaces
$f_*\colon
\m_A/\m_A^2\to (\epsilon) = \kappa$ is an element of $\tang A$.
In this way we obtain a bijective correspondence of the set of
algebra homomorphisms from $A$ to $\kappa[\epsilon]$ with the
dual vector space $\tang A$.

\proclaim Proposition. \call{altKS} Let $X$ be a deformation of
$X_0$ over a complete algebra $A$. If $u\in \tang A$ and
$f\colon A\to
\kappa[\epsilon]$ is the corresponding homomorphism of algebras,
then ${\rm K}_X(u) = e(f_*X,X_0[\epsilon])$.
\endproclaim

{}From this description the linearity of the \KS{} map is not
obvious.

\proof Use \ref{ksmapfunc} to reduce to the case $A =
\kappa[\epsilon]$, where the statement holds by
definition.\endproof

The \KS{} class of a deformation $X$ is determined by its
first-order part $X_1$; conversely the \KS{} class determines
$X_1$ completely.

\proclaim Proposition. \call{def<->ks}Let $A$ be  an artinian
algebra with $\m_A^2 = 0$. The
\KS{} class gives a bijective correspondence between isomorphism
classes of deformations on $A$ and
$\m_A\otimes T\simeq \hom(\m_A^\vee,T)$.
\endproclaim

\proof This follows immediately from \ref{absdefc}.\endproof

If $W$ is a vector space, then
$\kappa\oplus W$ has a canonical ring structure given by
$(a,x)(b,y) = (ab,bx+ay)$. If
$w_1,\dots,w_r$ is a basis of $W$ then
$$\kappa\oplus W =
\kappa[[w_1,\ldots,w_r]]/\m_{\kappa[[w_1,\ldots,w_r]]}^2.$$ An
artinian algebra  $A$ with $\m_A^2 = 0$ is of the form
$\kappa\oplus \m_A$.

Assume that ${\rm T}^1(X_0)$ is finite-dimensional. Set
$T = {\rm T}^1(X_0)$, and consider the artinian algebra
$R_1 = \kappa\oplus T^\vee$. Let
$V_1$ be the deformation of $X_0$ over $R_1$ corresponding to the
identity in $\hom(T,T)
\simeq T^\vee\otimes T$ (this is a very important deformation,
and we'll meet it again in the construction of the minimal
versal deformation of $X_0$ in
\ref{versal}.) Now take the graded algebra $R_2 =
\kappa\oplus T^\vee\oplus \sym^2T^\vee$, where $\sym^2T^\vee$ is
the second symmetric power of
$T^\vee$. If $t_1,\dots,t_r$ is a basis for $T^\vee$ then
$$ R_2 =
\kappa[[t_1,\ldots,t_r]]/\m_{\kappa[[t_1,\ldots,t_r]]}^3.
$$ Obviously $R_2/\sym^2T^\vee = R_1$ and
$\m_{R_2}(\sym^2T^\vee) = 0$.

\proclaim Definition. \call{firstobs}Assume that ${\rm
T}^1(X_0)$ is finite-dimensional. The\/ {\rm first obstruction
map} of $X_0$ is the linear map
$$Q_{X_0}\colon \sym^2{\rm T}^1(X_0)\longrightarrow \obs X_0$$
which corresponds to the obstruction
$$q_{X_0}\in \sym^2T^\vee\otimes \obs X_0\simeq \hom(\sym^2{\rm
T}^1(X_0), \obs X_0)$$ to lifting $V_1$ from $R_1$ to $R_2$.
\endproclaim

This important map induces a vector-valued quadratic form $T\to
\obs X_0$  sending a vector $u\in T$ into $Q_{X_0}(u\cdot u)$.
This map has the following interpretation. Consider the algebra
$\kappa[\epsilon]$ as before, and choose a vector
$u\in T$. Call $X(u)$ the deformation on $\kappa[\epsilon]$ whose
\KS{} class is $\epsilon
\otimes u\in
\kappa\epsilon\otimes T$. Set $\kappa[\epsilon'] =
\kappa[t]/(t^3)$, with the projection $\kappa[\epsilon']\to
\kappa[\epsilon]$ sending
$\epsilon'$ to $\epsilon$.  The obstruction to lifting $X(u)$ to
$\kappa[\epsilon']$ lives in
$\kappa\epsilon'^2\otimes
\obs X_0 = \obs X_0$.

\proclaim Proposition. The obstruction to lifting $X(u)$ to
$\kappa[\epsilon']$ is
$Q_{X_0}(u\cdot u)$.
\endproclaim

\proof Think of $u$ as a linear map $u\colon T^\vee\to \kappa$;
there is a unique homomorphism of graded algebras $f\colon
R_2\to \kappa[\epsilon']$ whose restriction $T\to
\kappa\epsilon = \kappa$ is $u$. Call $\phi\colon \sym^2
T^\vee\to
\kappa\epsilon'^2 =
\kappa$ the restriction of $f$; as an element of
$\sym^2T^{\vee\vee} =
\sym^2T$ the map
$\phi$ is exactly $u\cdot u$. By
\ref{absobsfunc} the obstruction to lifting $X(u)$ to
$\kappa[\epsilon']$ is
$\phi(q_{X_0})$; by linear algebra this is equal to
$Q_{X_0}(u\cdot u)$.\endproof

\beginsection Versal deformations [versal]

\proclaim Definition. A\/ {\rm universal deformation} of $X_0$
is a deformation $(V,R)$ such that for any deformation $(X,A)$
there is a unique homomorphism of deformations
$(\phi,f)\colon (V,R)\to (X,A)$.
\endproclaim

So a universal deformation contains all the information about
deformations, and is therefore a very good think to have.
Unfortunately universal deformations do not always exist; if
${\rm T}^0(X_0)\neq 0$, then the trivial deformation
$X_0[\epsilon]$ has nontrivial automorphisms (\ref{absdefd}), so
there can not exist a unique isomorphism
$\phi\colon f_*V\simeq X_0[\epsilon]$, because any such $\phi$
can always be composed with a nontrivial automorphism of $X_0$.

One could feel that imposing the condition of unicity on $\phi$
is being too demanding, and only require unicity for $f$. This
is indeed reasonable, from a ``functor-theoretic'', as opposed
to ``stack-theoretic'' point of view, and corresponds to the
condition of prorepresentability of the functor of isomorphism
classes of deformations in [Schlessinger]. Let us call a
deformation satisfying this weaker condition a {\it weak
universal deformation}; these exist for many more schemes
$X_0$. For example, assume that ${\rm T}^0(X_0)\neq 0$ and ${\rm
T}^1(X_0) = 0$, (e.g., when
$X_0 = {\bf P}^n$ or $X_0 = {\bf A}^n$, $n>0$). According to
\ref{deft1} all deformations are isomorphic to the trivial
deformation, and therefore the trivial deformation $(X_0,\kappa)$
is a weak universal deformation, while a universal deformation
does not exists. On the other hand not all projective
$X_0$ have weak universal deformation; for example, one can
prove that the ``banana'' curve
$$
\Bigg(\hskip-7pt{\Bigg)}
$$ consisting of two copies of ${\bf P}^1$ glued together at two
pairs of rational points does not possess one.

One could think that to obtain an acceptable replacement for the
notion of universal deformation it is sufficient to drop the
condition of unicity on $f$ and simply require the existence of
$f$ and $\phi$; but this turns out to be too weak. For example
if $(V,R)$ satisfies this condition then any deformation $(W,S)$
such that there exists a homomorphism $(W,S)\to (V,R)$ also
satisfies it.

The correct general notion is the following.

\proclaim Definition. \call{versaldef} A {\rm versal
deformation} of
$X_0$ is a deformation
$(V,R)$ such that if $(\eta,p)\colon (X,A)\to (Y,B)$ is a
surjective homomorphism of deformations, and
$(\psi,g)\colon (V,R)\to (Y,B)$ is a homomorphism, then there
exists a homomorphism
$(\phi,f)\colon (V,R)\to (Y,B)$ with $(\eta,p)\circ(\phi,f) =
(\psi,g)$.
$$\varrowlength=30pt
\commdiag{&&(X,A)\cr
&\arrow(2,1)\lft{(\phi,f)}&\mapdown\rt{(\eta,p)}\cr
(V,R)&\mapright^{(\psi,g)}&(Y,B)\cr }
$$

A versal deformation $(V,R)$ is called\/ {\rm minimal} if the
\KS{} map ${\rm K}_V\colon
\tang R\to {\rm T}^1(X_0)$ is an isomorphism.
\endproclaim

A minimal versal deformation is often called {\it miniversal}.

Actually, the standard definition of a versal deformation only
assumes the lifting property above when $A$ and $B$ are
artinian. Let us prove that it is equivalent to the definition
above.

\proclaim Lemma. \call{weakcondvers} Assume that a deformation
$(V,R)$ has the property of
\ref{versaldef} in the case that $A$ and $B$ are artinian. Then
$(V,R)$ is versal.
\endproclaim

\proof Assume that we are in the situation of \ref{versaldef},
with
$A$ and $B$ not necessarily artinian. To give a homomorphism
$(\phi,f)$ as above is equivalent to giving a homomorphism
$$(\phi_n,f_n)\colon (V,R)\to (X_n,A_n)$$for each $n$ such that
the composition of $(\phi_n,f_n)$ with the projection
$(X_n,A_n)\to (X_{n-1},A_{n-1})$ is $(\phi_{n-1},f_{n-1})$, and
$$(\eta_n,p_n)\circ(\phi_n,f_n) = (\psi_n,g_n)$$ for all $n$.

Let us assume that we have constructed $(\phi_{n-1},f_{n-1})$,
and consider the commutative diagram
$$
\diagram{60}{(V,R)& \mapright^{(\phi_{n-1},f_{n-1})}&
(X_{n-1},A_{n-1})\cr
\mapdown\rt{(\psi_n,g_n)}&&\mapdown\rt{(\eta_{n-1},p_{n-1})}\cr
(Y_n,B_n)&\mapright&(Y_{n-1},B_{n-1}). }$$ The diagram above
induce a homomorphism $(V,R)\to (Z,C)$, where $(Z,C)$ is the
fiber product of $(Y_n,B_n)$ and
$(X_{n-1},A_{n-1})$ over $(Y_{n-1},B_{n-1})$
(\ref{fiberproddef}); this homomorphism can be lifted to a
homomorphism
$(\phi_n,f_n)\colon (V,R)\to (X_n,A_n)$ with the required
properties.\endproof

Here are several properties of versal deformations.

\proclaim Proposition. \punto \call{univ->vers}A universal
deformation is miniversal.

\punto \call{bigvers} If $(V,R)$ is a versal deformation and
$(X,A)$ is a deformation, there is a homomorphism $(V,R)\to
(X,A)$.

\punto \call{versKS} If $(V,R)$ is a versal deformation, the
\KS{} map
${\rm K}_V\colon \tang R\to {\rm T}^1(X_0)$ is surjective.

In particular, if $X_0$ has a versal deformation then ${\rm
T}^1(X_0)$ is finite dimensional.

\punto \call{univuniq} Two universal deformations are canonically
isomorphic.

\punto \call{versuniq} Two miniversal deformations are
non-canonically isomorphic.
\endproclaim

\proof Part \ref{part bigvers} follows immediately from the
definition, if we take $(Y,B) = (X_0,\kappa)$.

For part \ref{part versKS}, take an element $a\in {\rm
T}^1(X_0)$; according to
\ref{absdefb}, there is a deformation $X$ on $\kappa[\epsilon]$
with
$e(X,X_0[\epsilon]) = a$. Choose a homomorphism $(\phi,f)\colon
(V,R)\to (X,\kappa[\epsilon])$, and let $u = f_*\colon
\cot R\to \kappa$ the corresponding element of $\tang R$. But
then
${\rm K}_V(u) = a$ because of \ref{altKS}.

For part \ref{part univ->vers}, let $(V,R)$ be universal
deformation. To prove that it is versal, let $(\eta,p)\colon
(X,A)\to (Y,B)$ be a surjective homomorphism of extensions,
$(\psi,g)\colon (V,R)\to (Y,B)$ is a homomorphism. There exists a
unique homomorphism
$(\phi,f)\colon (V,R)\to (Y,B)$, and by unicity
$(\eta,p)\circ(\phi,f) = (\psi,g)$.

The \KS{} map ${\rm K}_V\colon \tang R\to {\rm T}^1(X_0)$ is
surjective because of part \ref{part versKS}. Let $u\in
\tang R$ with ${\rm K}_V(u) = 0$, and let $f\colon R\to
\kappa[\epsilon]$ be the corresponding homomorphism. Then
$e(f_*V,X_0[\epsilon]) = {\rm K}_V(u) = 0$, so we can choose an
isomorphism $\phi\colon f_*V\simeq X_0[\epsilon]$, and we get a
homomorphism of extensions
$(\phi,f)\colon (V,R)\to (X_0[\epsilon],\kappa[\epsilon])$. If
we denote by $k\colon R\to
\kappa[\epsilon]$ the trivial homomorphism, corresponding to
$0\in
\tang R$, then $k_*V = X_0[\epsilon]$, so $({\rm id},k)\colon
(V,R)\to (X_0[\epsilon],\kappa[\epsilon])$ is also a
homomorphism, hence $f = k$, and $u=0$. It follows that ${\rm
K}_V$ is injective, and this completes the proof.

Part \ref{part univuniq} is standard and easy.

For part \ref{part versuniq}, take two miniversal deformations
$(V,R)$ and $(W,S)$. {}From part \ref{part bigvers} we see that
there exists homomorphisms
$(\phi,f)\colon (V,R)\to (W,S)$ and $(\psi,g)\colon (W,S)\to
(V,R)$; but the functoriality of the \KS{} map (\ref{ksmapfunc})
and the minimality condition insure that $f_*\colon
\cot R\to \cot S$ and
$g_*\colon \cot S\to \cot R$ are isomorphisms. {}From
\ref{lemmaauto} below we see that
$fg$ and $gf$ are isomorphisms, so $f$ and $g$ are
isomorphisms.\endproof

\proclaim Lemma. \call{lemmaauto} Let $A$ be a complete algebra,
$f\colon A\to A$ a homomorphism. If $f_*\colon \cot A\to \cot A$
is surjective then $f$ is an automorphism.
\endproclaim

\proof The linear map $\m_A^n/\m_A^{n+1}\to
\m_A^n/\m_A^{n+1}$ induced by $f$ is surjective for all $n\ge
0$, so
$f_n\colon A_n\to A_n$ is also surjective, and therefore an
isomorphism. Hence
$$ f\colon A = \projlim A_n\longrightarrow \projlim A_n = A
$$ is also an isomorphism.\endproof

The following is a generalization of \ref{versuniq}, and can be
considered as a description of all versal deformations.

\proclaim Proposition. \call{charversal} Let $(V,R)$ be a
miniversal deformation,
${\bf t} = (t_1,\ldots,t_n)$ a sequence of indeterminates,
$j\colon R\to R[[{\bf t}]]$ the inclusion. The deformation
$(V,R)[[{\bf t}]] = (j_*V,R[[{\bf t}]])$ is versal.

Conversely, if $(W,S)$ is another versal deformation, and $n$ is
the dimension of the kernel of the \KS{} map ${\rm K}_W\colon
\tang S\to {\rm T}^1(X_0)$, then
$(W,S)$ is isomorphic to
$(V,R)[[{\bf t}]]$.
\endproclaim

\proof Let $(\eta,p)\colon (X,A)\to (Y,B)$ a surjective
homomorphism,
$(\psi,g)\colon (V,R)[[{\bf t}]]\to (Y,B)$ a homomorphism. Since
$(V,R)$ is versal there will be a  homomorphism $(\phi,f')\colon
(V,R)\to (X,A)$ such that $(\eta,p)\circ (\phi,f') =
(\psi,g)\circ({\rm id},j)\equaldef (\psi,gj)$.

The homomorphism $f'\colon R\to A$ can be lifted to a
homomorphism $f\colon R[[{\bf t}]]\to A$ such that $pf = g$; if
we choose $a_i\in
\m_A$ for each $i = 1,\dots,n$ so that $p(t_i) = a_i$, there is
then only a homomorphism
$f\colon R[[{\bf t}]]\to A$ such that $f(x) = f'(x)$ if $x\in
R$, and
$f(t_i) = a_i$ for each $i$. This homomorphism $f$ has the
desired property.

Since $f_*j_*V = f'_*V$ the pair $(\phi,f)$ gives a homomorphism
$(V,R)[[{\bf t}]]\to (X,A)$. {}From the fact that $(\eta,p)\circ
(\phi,f') = (\psi,gj)$ we get that $\psi =
\eta\circ p_*\phi$; this together with $pf = g$ implies that
$(\eta,p)\circ(\phi,f) = (\psi,g)$, as desired. This proves the
first part of the statement.

The diagram below illustrates the proof.
$$
\vgrid=22pt
\harrowlength=20pt
\varrowlength=28pt
\gridcommdiag{&&&&&&&&(X,A)\cr
&&&&\hskip-10pt{\sarrowlength=90pt\arrow(3,1)\lft{(\phi,f')}}
&&\hskip-10pt{\sarrowlength=40pt\arrow(3,2)\rt{(\phi,f)}}
&&\mapdown\rt{(\eta,p)}\cr (V,R)&&\hskip-5pt\mapright^{({\rm
id},j)}&&(V,R)[[{\bf t}]]&&\hskip8pt\mapright^{(\psi,g)}&&(Y,B)}
$$

Now take a versal deformation $(W,S)$. Because of \ref{bigvers}
there is a homomorphism
$(\phi,f')\colon (V,R)\to (W,S)$; from the functoriality of the
\KS{} class (\ref{ksmapfunc}) and the surjectivity of the \KS{}
map $K_W$ (\ref{versKS}) we see that the differential
${\rm d}f'\colon \tang S\to \tang R$ is surjective, so its
transpose map $f'_*\colon \cot R\to \cot S$ is injective. Choose
elements
$a_1,\dots,a_n$ in $\m_S$ whose class in $\cot S$ form a basis
for a complement of the image of $f'_*$, and consider the
homomorphism
$f\colon R[[{\bf t}]]\to S$ which sends $x\in R$ into $f'(x)$,
and $t_i$ into $a_i$. Let us prove that
$(\phi,f)\colon (j_*V,R[[{\bf t}]])\to (W,S)$ is an isomorphism.
But
$g_1\colon R[[{\bf t}]]_1\to S_1$ is an isomorphism, so consider
the isomorphism of deformations
$(\phi_1,g_1)^{-1}\colon (W_1,S_1)\to (j_*V_1,R[[{\bf t}]]_1)$:
this can lifted to a homomorphism $(\psi,g)\colon (W,S)\to
(j_*,R[[{\bf t}]])$, because 
$(W,S)$ is versal. The compositions $fg$ and
$gf$ induce the identity in $\cot S$ and $\cot {R[[{\bf t}]]}$
respectively, so from
\ref{lemmaauto} they are isomorphisms. It follows that $f$ is an
isomorphism, and this concludes the proof.\endproof

\proclaim Proposition. \call{charuniv} A miniversal deformation
of
$X_0$ is universal if and only if ${\rm T}^0(X_0) = 0$.
\endproclaim

For this we need a lemma.

Let $A$ and $B'$ be artinian algebras, ${\goth b}\subseteq B'$ an
ideal with $\m_{B'} = 0$. Set $B = B'/{\goth b}$, and call
$\pi\colon B'\to B$ the projection. Let $f'_1$ and
$f'_2$ be homomorphisms $A\to B'$ such that $\pi f'_1 = \pi f'_2
= f$. The difference
$f'_1-f'_2\colon A\to {\goth b}$ is a derivation of the algebra
$A$ into the
$\kappa$-module {\goth b} (this is a standard fact, which we have
already used in the proof of \ref{isomlif}). So
$(f'_1-f'_2)(\m_A^2) = 0$, and $f'_1-f'_2$ induces a
$\kappa$-linear map
$$\Delta(f'_1,f'_2)\colon \cot A\to {\goth b}.$$

Let $X$ be a deformation of $X_0$; then
${f'_1}_*X$ and ${f'_2}_*X$ are liftings of $f_*X$ to $B'$. We
want understand when
${f'_1}_*X$ and ${f'_2}_*X$ are isomorphic; for this we need a
formula for
$$e\bigl({f'_1}_*X,{f'_2}_*X\bigr)\in {\goth b}\otimes{\rm
T}^1(X_0).$$ This turns out to be determined by
$\Delta(f'_1,f'_2)$ and the \KS{} class
${\rm k}_X\in (\cot A)\otimes{\rm T}^1(X_0)$.

\proclaim Lemma.
\call{diffhomks}$e\bigl({f'_1}_*X,{f'_2}_*X\bigr) =
\bigl(\Delta(f'_1,f'_2)\otimes{\rm id}\bigr)({\rm k}_X)\in {\goth
b}\otimes{\rm T}^1(X_0)$
\endproclaim

\proof Set $V = \cot A$. Call $\rho\colon A\to A/\m_A = \kappa$
the projection, $D_A\colon A\to V$ the derivation which sends
$a\in A$ into the class of $a-\rho(a)\in
\m_A$. Give to
$A' = A\oplus V$ the obvious ring structure in which $0\oplus V$
becomes an ideal with square 0; the multiplication is defined by
$(a,x)(b,y) = (ab, bx+ay)$. We call $\pi\colon A'\to A$ the
projection. The algebra
$A'$ is local, and $V= 0\oplus V = \ker \pi\subseteq A'$ is an
ideal which is killed by the maximal ideal $\m_{A'} = \m_A\oplus
V$ of $A'$. There are two homomorphism of algebras $i\colon A\to
A'$ and $u\colon A\to A'$ defined respectively by $i(a) = (a,0)$
and $u(a) = \bigl(a,D_A(a)\bigr)$. Finally, consider the
homomorphism of algebras $F\colon A'\to B'$ defined by $F(a,x) =
f'_2(a)+\Delta(f'_1,f'_2)(x)$.

Obviously $F\circ i = f'_2\colon A\to B'$, while
$$\eqalign{(F\circ u)(a) &= F\bigl(a,D_A(a)\bigr)\cr& = f'_2(a)+
f'_2\bigl(a-\rho(a)\bigr)-f'_1\bigl(a-\rho(a)\bigr)\cr & =
f'_2(a)+ f'_1\bigl(a)-\rho(a)-f'_2\bigl(a)+\rho(a)\cr & =
f'_1(a),}$$ so that $F\circ u = f'_1$.

The deformations $u_*X$ and $i_*X$ on $A'$ are liftings of $X$;
since
$F\rest V =
\Delta(f'_1,f'_2)$ we get from
\ref{extfunc} that $e\bigl({f'_1}_*X,{f'_2}_*X\bigr) =
(\Delta(f'_1,f'_2)\otimes{\rm id})e\bigl(u_*X,i_*X\bigr)$, so it
is enough to prove that
$$e\bigl(u_*X,i_*X\bigr) = {\rm k}_X.$$

Consider now $A_1 = A/\m_A^2 = \kappa\oplus V$, call
$\sigma\colon A\to A_1$ the projection. The homomorphism
$h\colon A'\to A_1$ defined by $h(a,x) = \rho(a)+x$ has the
property that
$h\circ u = \sigma$, while $(h\circ i)(a) = \rho(a)\in A_1$. It
follows that $h_*u_*X = X_1$, while
$h_*i_*X = X_0[\epsilon]$; by \ref{ksclassfunc} applied to the
homomorphism $h$, which sends
$x\in V\subseteq A'$ into $x\in V\subseteq A_1$, we get that
$e(u_*X,i_*X) = e(X_1,X_0[\epsilon])\in V\otimes{\rm T}^1(X_0)$.
But
$e(X_1,X_0[\epsilon]) = {\rm k}_X$ by definition.\endproof

\proofof charuniv. Let $(V,R)$ be a deformation and
$(\phi,f),(\psi,g)\colon (V,R)\to (X,A)$ two homomorphism; we
need to show that $(\phi,f) = (\psi,g)$. If $f = g$ then
$\phi = \psi$ because of \ref{deft0}, so it is enough to prove
the following: if $f,g\colon R\to A$ are homomorphisms of
complete algebras and $f_*V$ is isomorphic to $g_*V$ as a
deformation over $A$, then $f=g$. Obviously $f_0 = g_0$; we
assume
$n\ge 1$ and prove that if $f_{n-1} = g_{n-1}\colon R_{n-1}\to
A_{n-1}$ then
$f_n = g_n\colon R_n\to A_n$. Then $f_n = g_n$ for all $n$, so
$f=g$, as claimed.

Consider $\Delta(f_n,g_n)\colon \cot R\to
\m_A^n/\m_A^{n+1}\subseteq A_n$. By hypothesis
${f_n}_*X_n$ and
${g_n}_*X_n$ are isomorphic as liftings of
${f_{n-1}}_*X_{n-1} = {g_{n-1}}_*X_{n-1}$, so from
\ref{diffhomks} we get
$$\bigl(\Delta(f_n,g_n)\otimes{\rm id}\bigr)({\rm k}_V) = 0\in
\m_A^n/\m_A^{n+1}\otimes{\rm T^1}(X_0).$$ This is equivalent to
saying that the adjoint map
$$\Delta(f_n,g_n)^\vee\colon
\tang A\to \tang R$$ composed with the \KS{} map ${\rm
K}_V\colon \tang R\to {\rm T}^1(X_0)$ is 0. But ${\rm K}_V$ is
an isomorphism, so $\Delta(f_n,g_n) = 0$, and
$f_n = g_n$, as claimed.\endproof

We have seen that if $X_0$ has a versal deformation then ${\rm
T}^1(X_0)$ is finite-dimensional (\ref{versKS}). The main result
of this section is that the converse holds.

\proclaim Theorem. \call{mainversal} If ${\rm T}^1(X_0)$ is
finite-dimensional then $X_0$ has a miniversal deformation
$(V,R)$.

If $r = \dim_\kappa{\rm T}^1(X_0)$ then $R$ is  of the form
$\kappa[[t_1,\ldots,t_r]]/I$ with $I\subseteq
\m_{\kappa[[t_1,\ldots,t_r]]}^2$, and the minimal number of
generators of $I$ is the dimension of the space $\obs X_0$.
\endproclaim

So if ${\rm T}^1(X_0)$ is finite-dimensional then the
obstruction space $\obs X_0$ is finite-dimensional, even thought
${\rm T}^2(X_0)$ might not be.

The condition that ${\rm T}^1(X_0)$ be finite-dimensional is
satisfied for example when $X$ is proper, or affine with
isolated singularities.

Before going to the proof let us draw two consequences.

\proclaim Corollary. Let $(V,R)$ be a miniversal deformation of
$X_0$. Then $R$ is a power series algebra if and only if $X_0$
is unobstructed.
\endproclaim

The following is a consequence of the theorem and of
\ref{charuniv}.

\proclaim Corollary. The scheme $X_0$ has a universal
deformation if and only if ${\rm T}^0(X_0) = 0$ and ${\rm
T}^1(X_0)$ is finite-dimensional.
\endproclaim

\proofof mainversal. Let us start with a definition.

\proclaim Definition. The\/ {\rm order} of an artinian algebra
$A$ is the least $n$ such that
$\m_A^{n+1} = 0$.

A deformation $(V,R)$ is $n$-versal if the condition of
\ref{versaldef} is satisfied when
$A$ and $B$ are artinian of order at most $n$.
\endproclaim

Obviously every deformation is $0$-versal. A deformation $(V,R)$
is
$n$-versal of and only if $(V_n,R_n)$ is $n$-versal.

The content of \ref{weakcondvers} is that a deformation is versal
if and only if it is $n$-versal for all $n$.

Set $T = {\rm T}^1(X_0)$, and let $t_1,\dots,t_r$ be a basis for
the dual space
$T^\vee$. Assume also that $\obs X_0$ is finite-dimensional, with
basis
$\omega_1,\dots,\omega_\ell$. We do not know a priori that $\obs
X_0$ is finite-dimensional; if we do not assume this the proof
goes through with minor changes in notation. Set
$$
\Lambda = \prod_{i = 1}^\infty \sym^iT^\vee =
\kappa[[t_1,\ldots,t_r]].
$$

We will construct $R$ as a quotient of $\Lambda$ by an ideal
generated by $\ell$ generators
$f_1,\dots,f_\ell$ in $\Lambda$, with no terms of degree less
than 2. If we denote by
$f_i^{(n)}$ the part of $f_i$ consisting of terms of degree at
most
$n$, then $f_i^{(1)} = 0$, and the $f_i^{(n)}$ are constructed
by an inductive procedure in such a way that there is an
$n$-versal deformation $V_n$ over $R_n =
\Lambda\big/\bigl(\bigl(f_1^{(n)},\ldots,f_\ell^{(n)}\bigr)+
\m_\Lambda^{n+1}\bigr)$; then the $f_i$ are defined as the limit
of the
$f_i^{(n)}$.

We start from the deformation $(V_1,R_1)$ constructed after
\ref{def<->ks}; here
$R_1 = \kappa\oplus T^\vee = \Lambda/\m_\Lambda^2$, and $V_1$ is
a deformation having as the identity
$T\to T$ as its \KS{} map.

Let us prove that $(V_1,R_1)$ is 1-versal. Let $(X,A)$ and
$(Y,B)$ be deformations over artinian algebras of order at most
1 with a surjective homomorphism $(\eta,p)\colon (X,A)\to
(Y,B)$, and
$$(\psi,g)\colon (V_1,R_1)\to (Y,B)$$ a homomorphism. The \KS{}
class
$k_X\in \m_A\otimes T\simeq \hom(T^\vee,\m_A)$ induces a linear
map
$T^\vee\to 
\m_A$, and an algebra homomorphism $f\colon R_1\to A$; by the
functoriality of \KS{} classes (\ref{ksclassfunc}) and
\ref{def<->ks} we see that the composition $pf$ is equal to
$g$, and $p_*V_1$ is isomorphic to $X$. Fix an isomorphism
$\phi'\colon p_*V_1\simeq X$; if the composition
$$\psi' = \eta\circ p_*\phi'\colon g_*V_1 = p_*f_*V_1\to Y$$ is
equal to $\psi$ then $(\psi,g) = (\eta,p)\circ (\phi',f)$, and we
are done. In general this is not true.

Set $\beta = \eta^{-1}\circ \psi\circ p_*\phi'^{-1}\colon
p_*X\simeq p_*X$; then $\beta$ will correspond to an element of
$\overline \beta\in \m_B\otimes {\rm T}^0(X_0)$ (\ref{absdefd}).
Because of the surjectivity of $p$ the element $\overline \beta$
can be lifted to an element
$\overline
\alpha\in \m_A\otimes {\rm T}^0(X_0)$ corresponding to an
automorphism
$\alpha$ of $X$, and, because of
\ref{absautofunc}, we have
$p_*\alpha = \beta$. Set $\phi = \alpha\circ \phi'$; then
$$
\eta\circ p_*\phi  = \eta\circ p_*\alpha \circ p_*\phi' =
\eta\circ
\beta\circ p_*\phi' =
\eta\circ\eta^{-1}\circ\gamma\circ p_*\phi'^{-1}\circ p_*\phi =
\psi,
$$ so $(\phi,f)$ gives us the lifting.

Next consider the obstruction $q_{X_0}\in \sym^2T^\vee\otimes
\obs X_0 =
\m_\Lambda^2/\m_\Lambda^3\otimes \obs X_0$, as in \ref{firstobs},
which we write as
$$q_{X_0} = \sum_{i=1}^\ell f_i^{(2)}\otimes \omega_i,$$ and
think of $f_i^{(2)}$ as homogeneous polynomials of order 2. Then
we set
$$R_2 = \bigl(\kappa\oplus T^\vee\oplus
\sym^2T^\vee\bigr)\big/\bigl(f_1^{(2)},\ldots,f_\ell^{(2)}\bigr)
=
\Lambda\big/\bigl((f_1^{(2)},\ldots,f_\ell^{(2)})+
\m_\Lambda^3\bigr);$$ because of \ref{absobsfunc} the
obstruction to lifting $V_1$ to $R_2$ vanishes, so we choose an
arbitrary lifting
$V_2$. This deformation $(V_2,R_2)$ is 2-versal, as we shall see.

In general we proceed by induction. Let us assume that we have
lifted
$V_1$ to a deformation $(V_{n-1},R_{n-1})$, with $R_{n-1}$ an
algebra of the form
$$R_{n-1} = \bigg(\bigoplus_{k=0}^{n-1} \sym^k T^\vee\bigg)
\bigg/\bigl(f_1^{(n-1)},\ldots,f_\ell^{(n-1)}\bigr)=
\Lambda\big/\bigl(\bigl(f_1^{(n-1)}, \ldots,f_\ell^{(n-1)}\bigr)
+\m_\Lambda^n\bigr) =\Lambda/I_{n-1}$$ where the $f_1^{(n-1)}$
are polynomials in $t_1,\dots,t_r$ of degree at most $n-1$ with
no terms of degree less than 2. Suppose also that we know that
$R_{n-1}$ is $(n-1)$-versal. Then we look for homogeneous
polynomials
$g_1^{(n)},\dots,g_\ell^{(n)}$ of degree $n$ so that if we set
$f_i^{(n)} = f_i^{(n-1)}+g_i^{(n)}$, the obstruction to lifting
$V_{n-1}$ to
$$ R_n = \bigg(\bigoplus_{k=0}^n\sym^k T^\vee\bigg)
\bigg/\bigl(f_1^{(n)},\ldots,f_\ell^{(n)}\bigr)=
\Lambda\big/\bigl(\bigl(f_1^{(n)},\ldots,f_\ell^{(n)}\bigr)+
\m_\Lambda^{n+1}\bigr) =
\Lambda/I_n
$$ vanishes. 

To find the $g_i^{(n)}$ set
$$
\widetilde R_n
 = \Lambda/\m_\Lambda I_{n-1};
$$ the projection $\pi\colon \widetilde R_n\to R_{n-1}$ has
kernel
$I_{n-1}/\m_\Lambda I_{n-1}$. The obstruction $\omega\in
(I_{n-1}/\m_\Lambda I_{n-1})\otimes \obs X_0$ to lifting
$V_{n-1}$ to
$\widetilde R_n$ can be written as
$$
\omega = \sum_{i=1}^\ell u_i\otimes \omega_i,
$$ where the $u_i\in I_{n-1}$ are sums of a homogeneous
polynomial of degree $n$ and a linear combination of the
$f_i^{(n-1)}$ with coefficients in $\kappa$. Set
$$ I_n = \m_\Lambda I_{n-1} + (u_1,\ldots,u_\ell) =
\m_\Lambda^{n+1} +
\m_\Lambda f_1^{(n-1)}+\cdots + \m_\Lambda f_\ell^{(n-1)}+
(u_1,\dots,u_\ell)\subseteq \Lambda
$$ and consider the algebra
$$ R_n = \Lambda/I_n;
$$ with its projection $R_n\to R_{n-1}$. By \ref{absobsfunc} the
obstruction to extending
$V_{n-1}$ to $R_n$ vanishes. Let us check that $I_n$ is
generated by
$\m_\Lambda^{n+1}$ and by $\ell$ polynomials of the form
$f_i^{(n-1)}+g_i^{(n)}$, where
$g_i^{(n)}$ is homogeneous of degree $n$.

\proclaim Lemma. We have
$$ I_{n-1} = \m_\Lambda^n+ (u_1,\ldots,u_\ell)
$$ and
$$ I_n = \m_\Lambda^{n+1}+(u_1,\ldots,u_\ell).
$$
\endproclaim

\proof Clearly $\m_\Lambda^n+ (u_1,\ldots,u_\ell)\subseteq
I_{n-1}$ and
$\m_\Lambda^{n+1}+(u_1,\ldots,u_\ell)\subseteq I_n$, so we only
need to show the reverse inclusions.

Consider the natural algebra homomorphism $\pi\colon R'_{n-1}
\equaldef
\Lambda/\bigl(\m_\Lambda^n+I_n\bigr)\to R_{n-1}$. The algebra
$R'_{n-1}$ is a quotient of $R_n$, so there is a lifting of
$V_{n-1}$ to $R'_{n-1}$. By hypothesis
$(V_{n-1},R_{n-1})$ is $(n-1)$-versal; hence the projection
$\pi$ is split by a algebra homomorphism $\rho\colon R_{n-1}\to
R'_{n-1}$. But
$R'_{n-1}$ has order at most
$n-1$, and the differential ${\rm d}\pi$ induces an isomorphism
of tangent spaces, whose inverse is ${\rm d}\rho$; it follows
that
$\rho$ is surjective, and $\pi$ is an isomorphism. So
$$
\m_\Lambda^n+\m_\Lambda f_1^{(n-1)}+ \cdots+ \m_\Lambda
f_\ell^{(n-1)}+ (u_1,\ldots,u_\ell)
\equaldef
\m_\Lambda^n+I_n = I_{n-1} \equaldef
\m_\Lambda^n+\bigl(f_1^{(n-1)},\ldots,f_\ell^{(n-1)}\bigr)
$$ and we can write
$$ f_i^{(n-1)} = h_i+ \sum_{j=1}^\ell\alpha_{ij}f_j^{(n-1)}
+\sum_{j=1}^\ell
\lambda_{ij}u_j
$$ where $h_i\in \m_\Lambda^n$, $\alpha_{ij}\in \m_\Lambda$ and
$\lambda_{ij}\in \Lambda$. The matrix ${\rm
Id}_{\ell\times\ell}-(\alpha_{ij})\in {\rm GL}_n(\Lambda)$ is
invertible, and  the system of equations above implies that
$$
\bigl({\rm Id}_{\ell\times\ell}-(\alpha_{ij})\bigr)
\pmatrix{f_1^{(n-1)}\cr \vdots\cr f_\ell^{(n-1)}}\in
\bigl(\m_\Lambda^n+(u_1,\ldots,u_\ell)\bigr)^{\oplus \ell}
$$ and therefore $f_i^{(n-1)}\in
\m_\Lambda^n+(u_1,\ldots,u_\ell)$ for all $i$. This proves the
inclusion $I_{n-1}\subseteq \m_\Lambda^n+ (u_1,\ldots,u_\ell)$.
{}From this we also get
$$
\m_\Lambda^{n+1}+\m_\Lambda f_1^{(n-1)}+\cdots + \m_\Lambda
f_\ell^{(n-1)} \subseteq
\m_\Lambda^{n+1}+(u_1,\ldots,u_\ell)
$$ and this proves the inclusion $I_n\subseteq
\m_\Lambda^{n+1}+(u_1,\ldots,u_\ell)$.\endproof

So we have made progress; we have managed to express $R_n$ as the
quotient of
$\Lambda/\m_\Lambda^{n+1}$ by an ideal generated by $\ell$
polynomial
$u_1,\ldots,u_\ell$. This is not quite what what we want, as the
parts of degree less than
$n$ of these polynomials are not necessarily the $f_i^{(n-1)}$.
To fix this we need an elementary lemma.

\proclaim Lemma. Let $M$ be a finite module over a local ring
$\Lambda$, $x_1,\ldots,x_\ell$ and $y_1,\ldots,y_\ell$ two
sequences of generators. Then there exists an invertible
$\ell\times\ell$ matrix $(\lambda_{ij})\in {\rm
GL}_\ell(\Lambda)$ such that $y_i =
\sum_{j=1}^\ell \lambda_{ij}x_i$ for all $i$.
\endproclaim

\proof If the sequences are minimal then any matrix
$(\lambda_{ij})$ such that $y_i =
\sum_{j=1}^\ell \lambda_{ij}x_i$ for all $i$ is invertible, so
we are done. In general we permute the $x_i$ and the $y_i$ so
that $x_1,\ldots,x_p$ and
$y_1,\ldots,y_p$ are minimal sequences of generators, and write
$y_i = \sum_{j = 1}^p\lambda_{ij}x_i$ for all
$1\le i\le p$, $y_i = x_i+\sum_{j=1}^p \lambda_{ij}x_j$ for
$p+1\le i\le \ell$. The resulting matrix
$$
\pmatrix{\lambda_{1\,1}&\ldots&\lambda_{1\,p}&0&0&\ldots&0\cr
\vdots&\ddots&\vdots&\vdots&\vdots&\ddots&\vdots \cr
\lambda_{p\,1}&\ldots&\lambda_{p\,p}&0&0&\ldots&0\cr
\lambda_{p+1\,1}&\ldots&\lambda_{p+1\,p}&1&0&\ldots&0\cr
\lambda_{p+2\,1}&\ldots&\lambda_{p+2\,p}&0&1&\ldots&0\cr
\vdots&\ddots&\vdots&\vdots&\vdots&\ddots&\vdots \cr
\lambda_{\ell\,1}&\ldots&\lambda_{\ell\,p}&0&0&\ldots&1\cr }
$$ is invertible.\endproof

Now we apply the lemma to the ideal
$$ {\textstyle\m_\Lambda^n+(u_1,\ldots,u_\ell)\over \textstyle
\m_\Lambda^n} = {\textstyle\m_\Lambda^n+\bigl(f_1^{(n-1)},
\ldots,f_\ell^{(n-1)}\bigr)\over \textstyle
\m_\Lambda^n}\subseteq \Lambda_n
$$ and conclude that there is an invertible matrix
$(\lambda_{ij})\in {\rm GL}_{\ell\times\ell}(\Lambda)$ with
$f_i^{(n-1)}\equiv \sum_{j=1}^\ell
\lambda_{ij}u_j
\pmod{\m_l^n}$ for all $i$; therefore there are homogeneous
polynomials $g_i^{(n)}\in
\Lambda$ of degree
$n$ such that
$$f_i^{(n-1)}+g_i^{(n)}\equiv \sum_{j=1}^\ell \lambda_{ij}u_j
\pmod{\m_l^{n+1}}$$ for all $i$, so that the $f_i^{(n)} =
f_i^{(n-1)}+g_i^{(n)}$ generate
$I_n$, as we wanted.

Lift $V_{n-1}$ to a deformation $V_n$ over $R_n$. We want to
prove that $(V_n,R_n)$ is $n$-versal.

We know that $(V_n,R_n)$ is 
$(n-1)$-versal, because $R_{n-1} = R_n/\m_{R_n}^n$, and
$(V_{n-1},R_{n-1})$ is
$(n-1)$-versal. Let $(\eta,p)\colon (X,A)\to (Y,B)$ a surjective
homomorphism of deformations with $A$ and $B$ algebras of order
at most $n$,
$(\psi,g)\colon (V_n,R_n)\to (Y,B)$ a homomorphism; set $\a =
\ker p$.

First of all assume that $\a\subseteq \m_A^n$. In this case the
induced homomorphism
$A_{n-1}\to B_{n-1}$ is an isomorphism. Lift $g\colon R_n =
\Lambda/I_n\to B$ to a homomorphism $F'\colon \Lambda\to A$ by
lifting the images of the
$t_i$ in $B$ to $A$; clearly $F'(\m_\Lambda I_{n-1})\subseteq
\m_A(\m_A^n +\a) = 0$, so
$F'$ induces a homomorphism $F\colon \widetilde R_n =
\Lambda/\m_\Lambda I_{n-1}\to A$. By functoriality the
obstruction
$\sum_{i=1}^\ell u_i\otimes \omega_i\in I_{n-1}/\m_\Lambda
I_{n-1}\otimes {\rm T}^1(X_0)$ to lifting $V_{n-1}$ from
$R_{n-1}$ to $\widetilde R_n$ maps to 0 in
$\m_A^n\otimes {\rm T}^2(X_0)$, so
$F'(u_i) = 0$ for all $i = 1,\ldots,\ell$, and $F$ induces a
homomorphism $f'\colon R_n\to A$ which is a lifting of $g$. The
isomorphisms $\psi\colon g_*V_n = p_*f'_*V_n\simeq Y$ and
$\eta\colon p_*X\simeq Y$ make $f'_*V_n$ and $X$ into liftings
of $Y$; if there were an isomorphism of liftings $\phi\colon
f'_*V_n\simeq X$ then
$(\eta,p)\circ(\phi,f) = (\psi,g)$. In general there is no such
isomorphism, so consider the element
$e(X,f'_*V_n)\in
\a\otimes T$, and the corresponding linear map $u\colon
T^\vee\to \a$. Now we apply
\ref{diffhomks}: if $f\colon R_n\to A$ is a homomorphism such
that $\Delta(f,f') = u$, then
$e(f_*V_n,f'_*V_n) =  (u\otimes {\rm id})(k_{V_n})$. But
$k_{V_n}\in T^\vee\otimes T\simeq
\hom(T,T)$ is the element corresponding to the identity, so
$(u\otimes {\rm id})(k_{V_n}) = e(X,f'_*V_n)$, and
$e(f_*V_n,X) = e(f_*V_n,f'_*V_n)-e(X,f'_*V_n) = 0$. Therefore
$f_*V_n$ and $X$ are isomorphic as liftings, and the conclusion
follows.

In the general case one can factor $p\colon A\to B$ as the
projection
$\pi\colon A\to A/(\m_A^n\cap
\a)$ followed by a homomorphism $A/(\m_A^n\cap \a)\to B$; if we
can lift $(\psi,g)\colon (V,R)\to (Y,B)$ to a homomorphism
$(V,R)\to (\pi_*X,A/(\m_A^{n+1}\cap
\a))$ then it remains to lift along the homomorphism $\pi$,
whose kernel is contained in
$\m_A^n$, and this can be done by the previous case.

Consider the cartesian diagram
$$
\sqdiagram{40}{\bigl(\pi_*X,A/(\m_A^n\cap\a)\bigr )&\mapright&
(Y,B)\cr
\mapdown&&\mapdown\cr
(X_{n-1},A_{n-1})&\mapright&(Y_{n-1},B_{n-1})\cr }
$$ (\ref{fiberproddef}). Because $(V_n, R_n)$ is $(n-1)$-versal
we get a lifting of the composition of
$$({\rm id},\rho)\circ(\psi,g)\colon (V_n,R_n)\to
(Y_{n-1},B_{n-1}),$$ where
$\rho\colon A\to A_{n-1}$ is the projection, to a homomorphism
$(V_n,R_n)\to (X_{n-1},A_{n-1})$; from the diagram above we get
a lifting
$(V_n,R_n)\to
\bigl(X,A/(\m_A^n\cap\a)\bigr)$, and we conclude that
$(V_n,R_n)$ is
$n$-versal.

So by taking as $R = \Lambda/I$, where $I = (f_1,\ldots,f_r)$,
where each $f_i$ is the limit of the $f_i^{(n)}$, we get a
deformation
$(V,R)$ which is $n$-versal for each $n$. Since every artinian
algebra has an order, this means that $(V,R)$ has the property of
\ref{versaldef} in the case that $A$ and $B$ are artinian.
\ref{weakcondvers} implies that $(V,R)$ is versal.

Now we only have to prove that the minimal number of generators
of the ideal $I = (f_1,\ldots,f_\ell)$ is $\ell =
\dim_\kappa\obs X_0$. For this we will produce a surjective
linear map $(I/\m_\Lambda I)^\vee\to \obs X_0$.

Consider the ideal $J_n = I\big/\bigl(I\cap
\m_\Lambda^{n+1}\bigr)$ in
$\Lambda_n =
\Lambda/\m_\Lambda^{n+1}$; we have $\Lambda_n/J_n = R_n$. Take
the induced surjective morphism
$$I/\m_\Lambda I\to J_n/\m_\Lambda J_n.$$

\proclaim Lemma. The surjective morphism $I/\m_\Lambda I\to
J_n/\m_\Lambda J_n$ is an isomorphism for $n\gg 0$.
\endproclaim

\proof This is equivalent to saying that $I\cap
\m_\Lambda^{n+1}\subseteq \m_\Lambda I$ for
$n\gg 0$, which follows from the Artin-Rees lemma.\endproof

Now consider the obvious surjective homomorphism of algebras
$\Lambda_n/\m_\Lambda J_n\to R_n$; the deformation
$(V_{n-1},R_{n-1})$ has an obstruction
$\omega\in (J_n/\m_\Lambda J_n)\otimes
\obs X_0$. Let us show that the associated linear map $u\colon
(J_n/\m_\Lambda J_n)^\vee\to
\obs X_0$ is surjective for $n\gg 0$.

In fact we can find a basis $\omega_1,\ldots,\omega_\ell$ of
$\obs X_0$ and for each
$i$ a small extension of artinian algebras $A'_i\to A_i$ with
kernel
$\a_i\simeq
\kappa$ and a deformation $(X_i,A_i)$ whose obstruction in
$\a_i\otimes \obs X_0\simeq \obs X_0$ is exactly $\omega_i$. Now
pick a homomorphism
$(\phi_i,f_i)\colon (V_n,R_n)\to (X_i,A_i)$ for some $n\gg 0$,
and lift $f_i\colon R_n\to A_i$ to a homomorphism $f'_i\colon
\Lambda_n/\m_\Lambda J_n\to A'_i$. Call $g_i\colon J_n/\m_\Lambda
J_n\to \a\simeq \kappa$ the restriction of $f'_i$; by the
functoriality of the obstruction class (\ref{absobsfunc}) we have
$$(g_i\otimes{\rm id})(\omega) = \omega_i,$$ which is equivalent
to $u(g_i) = \omega_i$. This proves the surjectivity of $u$, and
concludes the proof.\endproof

In the unobstructed case the miniversal deformations are easy to
characterize.

\proclaim Corollary. \call{charunobs} Assume $\obs X_0 = 0$. A
deformation $(X,A)$ is miniversal if and only if $A$ is a power
series algebra over
$\kappa$ and the \KS{} map
$K_X\colon
\tang A\to {\rm T}^1(X_0)$ is an isomorphism.
\endproclaim

\proof Let $(V,R)$ be a miniversal deformation constructed above;
then $R$ is a power series algebra. Let
$$(\phi,f)\colon (V',R')\to (V,R)$$ be a homomorphism; then
$(X,A)$ is miniversal if and only if $f$ is an isomorphism, and
it easy to check that this happens if and only if $A$ is a power
series algebra and the differential ${\rm d}f\colon \tang A \to
\tang R$ is an isomorphism. By \ref{ksmapfunc}
${\rm d}f\colon \tang A \to \tang R$ is an isomorphism if and
only if
$K_X\colon \tang A\to {\rm T}^1(X_0)$ is an isomorphism, and
this concludes the proof.\endproof

The following is one of the the simplest nontrivial examples of a
versal deformation space.

\proclaimr Example. Let $X_0\subseteq {\bf A}^n_\kappa$ a
hypersurface with isolated singularities, $n\ge 2$. Call
$$F\in \kappa[{\bf x}]\equaldef \kappa[x_1,\ldots,x_n]$$ a
generator of the ideal of $X$ in ${\bf A}^n_\kappa$, and set $A =
\kappa[{\bf x}]/(F)$. The ideal
$$ J = \left({\textstyle \partial F\over \textstyle \partial
x_1},\ldots, {\textstyle
\partial F\over \textstyle \partial x_1}\right)R\subseteq R
$$ generated by the images in $R$ of the partial derivatives of
$F$ is called the {\it Jacobian ideal\/} of
$X$. Look at the basic exact sequence
$$
\sequence{0& \mapright I_0/I_0^2&\mapright& \Omega_{{\bf
A}^n/\kappa}\otimes_{\kappa[{\bf x}]} A& \mapright&
\Omega_{X_0/\kappa}&\mapright&0 }\eqdef{succ}$$ where $I_0$ is
the ideal of $X_0$ in $\kappa[{\bf x}]$ (\ref{exseq}). The
$A$-module
$\Omega_{{\bf A}^n/\kappa}\rest{X_0}$ is free of rank $n$,
generated by
${\rm d}x_1,\ldots,{\rm d}x_n$, while $I_0/I_0^2$ is free of
rank 1, generated by the class of $F$. The matrix of the
differential $I_0/I_0^2\to \Omega_{{\bf A}^n/\kappa}$ is the
Jacobian matrix of $F$ restricted to $X_0$, so the image of its
adjoint $\left(\Omega_{{\bf
A}^n/\kappa}\rest{X_0}\right)^\vee\to
\left(I_0/I_0^2\right)^\vee$ is
$J\left(I_0/I_0^2\right)^\vee
\simeq J$. We have an exact sequence 
$$\sequence{ 0&\mapright& \left(\Omega_{X_0/\kappa}\right)^\vee&
\mapright&
\left(\Omega_{{\bf A}^n/\kappa}\otimes_{\kappa[{\bf x}]}
A\right)^\vee&\mapright&
\left(I_0/I_0^2\right)^\vee&\mapright^{\textstyle\partial}&
\ext^1_A \bigl(\Omega_{X_0/\kappa},A\bigr)&\mapright& 0 }$$
which gives an isomorphism $\ext^1_A
\bigl(\Omega_{X_0/\kappa},A\bigr)\simeq A/J$. This shows that
the ideal $J$ is the annihilator of $\ext^1_A
\bigl(\Omega_{X_0/\kappa},A\bigr)$ in $A$, so it is an invariant
of
$X_0$ and does not depend on the embedding $X_0\into {\bf
A}^n_\kappa$.

The dimension $\mu = \mu(X_0)$ of $\ext^1_A
\bigl(\Omega_{X_0/\kappa},A\bigr) = A/J$ as a
$\kappa$-vector space is called the {\it Milnor number\/} of
$X_0$. Choose polynomials
$G_1,\ldots,G_\mu\in \kappa[{\bf x}]$ whose images
$g_1,\ldots,g_\mu$ in $A/J$ form a basis, and consider the
subscheme $V\subseteq {\bf A}^n_{\kappa[[{\bf t}]]}$, where
${\bf t} = (t_1,\ldots,t_\mu)$, whose ideal is generated by
$$G = F+t_1G_1+\cdots+t_\mu G_\mu\in \kappa[[{\bf t}]][{\bf
x}];$$ we denote also by $V$ the formal deformation on $\Lambda =
\kappa[[{\bf t}]]$  induced by
$V$. I claim that $(V,R)$ is a miniversal deformation of $X_0$.

Observe that from the sequence \eqref{succ} we get that
$\ext^2_A(\Omega_{X_0/\kappa},A) = 0$, so $X_0$ is unobstructed.
By \ref{charunobs} it is enough to prove that the \KS{} map
$K_V\colon \tang R\to {\rm T}^1(X_0)$ is an isomorphism.

Consider the deformation $(V_1, \Lambda_1)$; we want to calculate
$k_V = e\bigl(V_1,X_0^{\Lambda_1}\bigr)$. According to
\ref{nuext} this is the image of
$$
\nu\bigl(V_1,X_0^{\Lambda_1}\bigr)\in
\bigl(\m_\Lambda/\m_{\Lambda}^2\bigr)\otimes_\kappa
(I_0/I_0^2)^\vee
$$ by the linear map
$${\rm id}\otimes\partial\colon
\bigl(\cot\Lambda\bigr)\otimes_\kappa
(I_0/I_0^2)^\vee\longrightarrow
\bigl(\cot\Lambda\bigr)\otimes_\kappa
\ext^1_A(\Omega_{X_0/\kappa},A).$$ Under the isomorphism of
$A$-modules
$A\simeq (I_0/I_0^2)^\vee$ the identity corresponds to the
element
$(I_0/I_0^2)^\vee =
\hom_A(I_0,A)$ which sends $F$ to 1, so to calculate the element
of
$$\bigl(\cot\Lambda\bigr)\otimes_\kappa (I_0/I_0^2)^\vee \simeq
\bigl(\cot\Lambda\bigr)\otimes_\kappa A$$ corresponding to
$\nu\bigl(V_1,X_0^{\Lambda_1}\bigr)$ it is enough to calculate
$\nu\bigl(V_1,X_0^{\Lambda_1}\bigr)(F)$. For this we follow the
definition of
$\nu\bigl(V_1,X_0^{\Lambda_1}\bigr)$ (\ref{embedded}). The image
of
$G\in
\Lambda[{\bf x}]$ in $\Lambda_1[{\bf x}]$ is an element of  the
ideal of $V_1$ in ${\bf A}^n_{\Lambda_1}$ mapping to $F$ in
$\kappa[{\bf x}]$, while the image of $F$ in 
$\Lambda_1[{\bf x}]$ is an element of the ideal of
$X_0^{\Lambda_1}$ mapping to $F$ in $\kappa[{\bf x}]$. The
difference
$$ G-F = \sum_{i=1}^\mu t_i G_i\in \m_{\Lambda_1}[{\bf x}] =
\bigl(\cot
\Lambda\bigr)\otimes \kappa[{\bf x}]
$$ maps to $\nu\bigl(V_1,X_0^{\Lambda_1}\bigr)(F)$ in $(\cot
\Lambda)\otimes A$; this means that the image $k_{V_1} =
e\bigl(V_1,X_0^{\Lambda_1}\bigr)$ of
$\nu\bigl(V_1,X_0^{\Lambda_1}\bigr)$ in
$\ext^1_A(\Omega_{X_0},A) = A/J$ is
$\sum_{i=1}^\mu t_i\otimes g_i$, where the $g_i$ are the classes
of the $G_i$ modulo $J$. Because by construction these form a
basis of $A/J$, it follows that the \KS{} map is an isomorphism,
as claimed.
\endproclaim

This example can be generalized to complete intersections of
positive dimension with isolated singularities in affine spaces.

\beginsection A technical lemma [technical]

Let $A$ be a noetherian commutative ring, $I_1$ and $I_2$ two
nilpotent ideals in $A$; set
$\widetilde A = A/(I_1\cap I_2)$, $A_1 = A/I_1$, $A_2 = A/I_2$,
$A_0 = A/(I_1+I_2)$, and call
$\pi_i\colon A_i\to A_0$ the two projections. Consider the
category
$\cal F$ of flat schemes of finite type over $\widetilde A$.
Define the objects of the category $\cal C$ to be quintuples
$(X_1,X_2,X_0,\alpha_1,\alpha_2)$, where
$X_i$ is a flat scheme of finite type over $A_i$ and
$\alpha_i\colon X_0\to X_i$ is a closed embedding of schemes
over $A$ inducing an isomorphism $X_i\rest{\spec {A_0}}\simeq
X_0$. The arrows from
$(X_1,X_2,X_0,\alpha_1,\alpha_2)$ to
$(Y_1,Y_2,Y_0,\beta_1,\beta_2)$ in
$\cal C$ are triples $(f_1,f_2,f_0)$ of morphisms $f_i\colon
X_i\to Y_i$ of schemes over $A$ such that $\beta_i f_0 = f_i
\alpha_i$ for $i = 1$,~2.

There is a functor $\Phi\colon {\cal F}\to \cal C$ which sends a
scheme
$X\to\spec\widetilde A$ into
$$ (X\rest{\spec{A_1}},X\rest{\spec{A_2}},
X\rest{\spec{A_0}},\alpha_1,\alpha_2)
$$ where $\alpha_i\colon X\rest{\spec{A_0}}\to
X\rest{\spec{A_i}}$ is the obvious embedding.

If $X_i = X\rest{\spec{A_i}}$ there are embeddings $\iota_1\colon
X_1\into X$ and
$\iota_2\colon X_2\into X$ with $\iota_1\circ \alpha_1 =
\iota_2\circ\alpha_2$, so from a morphism of $A$-schemes
$\phi\colon X\to Z$ we get two morphisms
$\psi_i = \phi\circ\iota_i$ with
$\psi_1\alpha_1 = \psi_2\alpha_2$.

\proclaim Lemma. \call{techlemma}The functor $\Phi$ is an
equivalence of categories.

Furthermore the construction above yields a bijective
correspondence between morphisms of
$A$-schemes $\phi\colon X\to Z$ and pairs of morphisms
$\psi_1\colon X_1\to Z$ and $\psi_2 = X_2\to Z$ with
$\psi_1\alpha_1 = \psi_2\alpha_2$.
\endproclaim

Call $\pi_i\colon A_i\to A_0$ and $\rho_i\colon \widetilde A\to
A_i$ the projections.

\proclaim Corollary. \call{fiberproddef}Suppose that $A$ is a
local artinian
$\kappa$-algebra with residue field $\kappa$. Let
$X^{(i)}$ be a deformation of
$X_0$ over
$A_i$ for $i = 0$, 1,~2, $(\phi_i,\pi_i)\colon
\bigl(X^{(i)},A_i\bigr)\to 
\bigl(X^{(0)},A_0\bigr)$ homomorphism of deformation. Then there
is a deformation $X$ over
$\widetilde A$ and homomorphisms $(\phi_i,\rho_i)\colon
\bigl(X,\widetilde A\,\bigr)\to
\bigl(X^{(i)},A_i\bigr)$ such that
$$
\sqdiagram{40}{\bigl(X,\widetilde
A\,\bigr)&\mapright^{(\psi_2,\rho_2)}&\bigl(X^{(2)},A_2\bigr)\cr
\mapdown\rt{(\psi_1,\rho_1)}&&\mapdown\rt{(\phi_2,\pi_2)}\cr
\bigl(X^{(1)},A_1\bigr)&\mapright^{(\phi_1,\pi_1)}
&\bigl(X^{(0)},A_0\bigr)\cr }
$$ is a cartesian diagram of deformations.
\endproclaim

\proof We have ${\pi_i}_*X^{(i)} = X^{(i)}\rest{\spec A_0}$, so
we get an object
$\bigl(X^{(1)},X^{(2)},X^{(0)},\alpha_1,\alpha_2\bigr)$ of $\cal
C$, where $\alpha_i$ is the composition of the isomorphism
$\phi_i^{-1}\colon X^{(0)}\simeq X^{(i)}\rest{\spec A_0}$ with
the embedding
$X^{(i)}\rest{\spec A_0}\into X^{(i)}$. If $X$ is a flat scheme
of finite type over
$\widetilde A$ such that $\Phi(X)\simeq
\bigl(X^{(1)},X^{(2)},X^{(0)},\alpha_1,\alpha_2\bigr)$ then we
get isomorphisms $\psi_i\colon {\rho_i}_*X = X\rest{\spec
A_i}\simeq X^{(i)}$. Then $\bigl(X,\widetilde A\,\bigr)$ is the
desired fiber product.\endproof

\proofof techlemma. We may assume $A = \widetilde A$. Let us
begin with some algebraic preliminaries.

Consider the category ${\cal F}_A$ of flat module over $A$, and
the category ${\cal F}_{A_\mini}$ whose objects are of quintuples
$$ M_\mini = (M_1,M_2,M_0,\alpha_1^{M_\mini},\alpha_2^{M_\mini})
$$ where $M_i$ is a flat module over $A_i$, and
$\alpha_i^{M_\mini}\colon M_i\to M_0$ is a homomorphism of
$A$-modules inducing an isomorphism
$M_i\otimes_{A_i}A_0 \simeq M_0$; we will call these objects
{\it flat $A_\mini$-modules}. The homomorphisms of
$A_\mini$-modules
$f_\mini\colon M_\mini\to N_\mini$ are triples
$f_\mini = (f_1,f_2,f_0)$ of homomorphisms of $A$-modules
$f_i\colon M_i\to N_i$ such that
$f_0\alpha_i^{M_\mini} = \alpha_i^{N_\mini}f_i$ for $i = 1$,~2.

The homomorphism $f_\mini$ is called {\it surjective\/} if each
$f_i$ is surjective.

There is functor
$U\colon {\cal F}_A\to {\cal F}_{A_\mini}$ which sends a flat
module
$M$ into
$$ M\otimes A_\mini = (M\otimes_A A_1,M\otimes_A A_2, M\otimes_A
A_0,\alpha_0^M,\alpha_1^M)
$$ where $\alpha_i\colon M_i\otimes_A A_i\to M_0\otimes_A A_0$
is induced by the projection
$\pi_i\colon A_i\to A_0$. A homomorphism of flat $A$-modules
$f\colon M\to N$ induces a homomorphism
$$ U(f) = (f\otimes{\rm id}_{A_1},f\otimes{\rm
id}_{A_2},f\otimes{\rm id}_{A_0})\colon M_\mini\longrightarrow
N_\mini,
$$ and this makes $U$ into a functor.

We want to show that $U$ is an equivalence of categories; let us
construct an inverse
$V\colon {\cal F}_{A_\mini}\to {\cal F}_A$. If $M_\mini$ is a
flat ${A_\mini}$module then we define $V(M_\mini)$ to be the
equalizer of the pair of morphisms
$(\alpha_1^{M_\mini},\alpha_2^{M_\mini})$, or, in other words,
the kernel of the homomorphism of $A$-modules
$\delta_{M_\mini}\colon M_1\times M_2\to M_0$ defined by
$\delta_{M_\mini}(x_1,x_2) =
\alpha_1^{M_\mini}(x_1)-\alpha_2^{M_\mini}(x_2)$. The obvious
homomorphism of rings $A\to A_1\times A_2$ makes the $A_1\times
A_2$-module $M_1\times M_2$ into an
$A$-module, and $V(M_\mini)$ is an $A$-submodule. If
$f_\mini\colon M_\mini\to N_\mini$ is a homomorphism of
$A_\mini$-modules then we define $V(f_\mini)$ to be the
restriction of
$f_1\times f_2\colon M_1\times M_2\to N_1\times N_2$, so that we
have a commutative diagram with exact rows
$$
\diagram{30}{0&\mapright& V(M_\mini)& \mapright& M_1\times M_2&
\mapright^{\delta_{M_\mini}}& M_0&\mapright&0\cr &&\mapdown\rt{
V(f_\mini)}&&\mapdown\rt{f_1\times f_2} &&
\mapdown\rt{f_0}\cr 0& \mapright& V(N_\mini)& \mapright&
N_1\times N_2&
\mapright^{\delta_{N_\mini}}& N_0&\mapright&0.\cr }
$$ This makes $V$ into a functor from ${\cal F}_{A_\mini}$ into
the category ${\cal M}_A$ of
$A$-modules.We need to show that $V(M_\mini)$ is flat over $A$
for any flat
$A_\mini$-module $M_\mini$, and to produce isomorphisms of
functors
$UV\simeq {\rm id}$ and
$VU\simeq {\rm id}$.

The chinese remainder theorem gives us an exact sequence of
$A$-modules
$$
\diagram{30}{0&\mapright& A&\mapright& A_1\times A_2&
\mapright^{\pi_1-\pi_2}& A_0&
\mapright& 0}
$$ which we can tensor with a flat $A$-module $M$ to get an exact
sequence
$$
\diagram{30}{0&\mapright& M&\mapright& (M\otimes_A
A_1)\times(M\otimes_A A_2)& \mapright& M\otimes_A A_0&
\mapright& 0.\cr}
$$ This gives a canonical isomorphism of $A$-modules $M\simeq
V(M\otimes A_\mini)$, which yields an isomorphism between the
functor $VU$ and the embedding of ${\cal F}_A$ into ${\cal M}_A$.

Now take a flat $A_\mini$-module $M_\mini$, and set $M =
V(M_\mini)$. The homomorphisms
$M\to M_1$ and $M\to M_2$ coming from the inclusion $M\subseteq
M_1\times M_2$, and the homomorphism $M\to M_0$ induced by
either of the two projections
$M_1\times M_2\to M_0$ induce homomorphisms of
$A_i$-modules $\phi_i^M\colon M\otimes_A A_i\to M_i$; we want to
show that the $\phi_i^M$ are isomorphisms, and that 
$M$ is flat. Once this is done we can restrict $V$ to a functor
$V\colon {\cal F}_{A_\mini}\to {\cal F}_A$, and the isomorphisms
$\phi_\mini\colon UV(M_\mini) = M\otimes_A A_\mini\to M_\mini$
will give an isomorphism of functors $UV\simeq {\rm id}$, as
claimed.

Choose a free $A$-module $F$ and a surjective homomorphism
$f_0\colon F\otimes A_0\to M_0$; this lifts to surjective
homomorphisms $f_1\colon F\otimes A_1\to M_1$ and $f_2\colon
F\otimes A_2\to M_2$, yielding a surjective homomorphism
$f_\mini\colon F_\mini = F\otimes_A A_\mini\to M_\mini$. Call
$K_i$ the kernel of
$f_i\colon F_i\to M_i$; because of the flatness of $M_i$ the
restriction
$\alpha_i^{F_\mini}\rest{K_i}\colon K_i\to K_0$ induces an
isomorphism
$K_i\otimes_A A_0
\simeq K_0$, so that
$$ K_\mini =
(K_1,K_2,K_0,\alpha_1^{F_\mini}\rest{K_1},\alpha_2^{F_\mini}\rest{K_2})
$$ is a flat $A_\mini$-module. Clearly $F = VU(F) = V(F_\mini)$,
and $\phi_i^F\colon F\otimes_A A_i\to F_i$ is the identity. Set
$K = V(K_\mini)$. We get a commutative diagram
$$
\diagram{30}{&&0&&0&&0\cr &&\mapdown&&\mapdown&&\mapdown\cr
0&\mapright& K& \mapright& K_1\times K_2&
\mapright^{\delta_{K_\mini}}& K_0&\mapright&0\cr
&&\mapdown&&\mapdown && \mapdown\cr 0&\mapright& F& \mapright&
F_1\times F_2&
\mapright^{\delta_{F_\mini}}& F_0&\mapright&0\cr
&&\mapdown&&\mapdown\rt{f_1\times f_2} && \mapdown\rt{f_0}\cr
0&\mapright& M& \mapright& M_1\times M_2&
\mapright^{\delta_{M_\mini}}& M_0&\mapright&0\cr
&&\mapdown&&\mapdown&&\mapdown\cr &&0&&0&&0\cr }
$$ with exact rows, whose last two columns are also exact; this
implies that the first column is also exact.

What follows is a familiar argument in commutative algebra. When
we we tensor it with
$A_i$ we get a right exact sequence, which is the top row of a
commutative diagram
$$
\diagram{30}{&&K\otimes A_i&\mapright& F\otimes A_i&\mapright&
M\otimes A_i& \mapright&0\cr
&&\mapdown\rt{\phi_i^K}&&\bivline\rt{\phi_i^F}&&
\mapdown\rt{\phi_i^M}\cr 0&\mapright&K_i& \mapright&
F_i&\mapright& M_i&\mapright&0\cr }
$$ with exact rows. This implies that $\phi_i^M$ is surjective;
since $M$ is an arbitrary module
$\phi_i^K$ will be surjective too, and by diagram chasing we get
that
$\phi_i^M$ is an isomorphism, as claimed.

Then $\phi_i^K$ will also be an isomorphism; this means that the
sequence 
$$
\sequence{0&\mapright&K\otimes A_i&\mapright& F\otimes
A_i&\mapright& M\otimes A_i&
\mapright&0\cr}
$$ is exact, so $\tor_1^A(A_i,M) = 0$. Since $M\otimes_A A_i =
M_i$, and
$M_i$ is flat over
$A_i$, the flatness of $M$ over $A$ follows immediately from the
local criterion of flatness (see [Matsumura], or the proof of
\ref{critflatness}).

Now we go from flat modules to flat algebras. This is
straightforward; call ${\cal F}^{\rm alg}_A$ the category of
flat $A$-algebras, and ${\cal F}^{\rm alg}_{A_\mini}$ the
category of flat $A_\mini$-algebras; that is, flat
$A_\mini$-modules $M_\mini$ such that every $M_i$ is endowed
with a structure of $A_i$-algebra, so that $\alpha_i^{M_\mini}$
is a homomorphism of
$A$-algebras. Then if $M$ is an $A$-algebra the $A_\mini$-module
inherits a structure of
$A_\mini$-algebra, so we get a functors $U^{\rm alg}\colon {\cal
F}^{\rm alg}_A\to {\cal F}^{\rm alg}_{A_\mini}$; conversely if
$M_\mini$ is an
$A_\mini$-algebra then $V(M_\mini)$ is an $A$-subalgebra of
$M_1\times M_2$, so $V(M_\mini)$ has a natural
$A$-algebra structure, and we get a functor $V^{\rm alg}\colon
{\cal F}^{\rm alg}_{A_\mini}\to {\cal F}^{\rm alg}_A$; these
together give an equivalence of categories.

Also from the construction we get immediately that if $M_\mini =
U(M)$ and $N$ is an
$A$-algebra then there is a bijective correspondence between
homomorphisms of algebras
$\phi\colon N\to M$ and pair of homomorphisms of $A$-algebras
$\psi_1\colon N\to M_1$ and
$\psi_2\colon N\to M_2$ such that $\alpha_1^{M_\mini}\psi_1 =
\alpha_2^{M_\mini}\psi_2$, obtained by composing $\phi$ with the
projections $M\to M_1$ and
$M\to M_2$.

If $M$ is a flat $A$-algebra of finite type then $M\otimes_A
A_i$ is of finite type over
$A_i$. On the other hand suppose that $M\otimes_A A_1$ is of
finite type over $A_1$, and choose a surjective homomorphism
$A_1[x_1,\ldots,x_n]\to M_1 = M/I_1$. By lifting the images of
the $x_i$ to $M$ we get a homomorphism $A[x_1,\ldots,x_n]\to M$,
which is easily seen to be surjective, due to the fact that
$I_1$ is nilpotent, so $M$ is also of finite type.

This proves \ref{techlemma} for affine schemes. The general case
follows from standard patching arguments, which we omit.\endproof

\references

M.\ Artin: (1)\enspace Lectures on Deformations of
Singularities. Tata Institute of Fundamental Research, Bombay
(1976).\quad (2)\enspace Versal deformations and algebraic
stacks. Invent{.} Math{.} {\bf 27} (1974), pp.\ 165--189.\\

P.\ Deligne, expos\'e  XVIII in: A{.} Grothendieck, with M.\
Artin and J{.} L{.} Verdier: Th\'eorie des topos et chomologie
\'etale de schemas, S\'eminaire de G\'eom\'etrie n. 4, vol{.}
III, Lecture Notes in Mathematics {\bf 305}. Springer-Verlag,
Berlin-New York, 1973.\\

R{.} Elkik: Alg\'ebrisation du module formel d'une singularit\'e
isol\'ee. In: Quelques probl\`emes de modules (S\'em{.} G\'eom{.}
Anal., \'Ecole Norm. Sup., Paris, 1971--1972), Ast\'erisque {\bf
16} (1974), Soc{.} Math{.} France, Paris, pp. 133--144.\\

L{.} Illusie: Complexe cotangent et d\'eformations I: Lecture
Notes in Mathematics {\bf 239}. Springer-Verlag, Berlin-New
York, 1971.
\enspace Complexe cotangent et d\'eformations II: Lecture Notes
in Mathematics {\bf 283}. Springer-Verlag, Berlin-New York,
1972.\\
 
J.\ Koll\'ar: Rational Curves on Algebraic Varieties. Ergebnisse
Der Mathematik Und Ihrer Grenzgebiete {\bf 32}, Springer-Verlag
(1995).\\

H.\ Matsumura: Commutative algebra. W{.} A{.} Benjamin, Inc{.},
New York (1970).\\

V.\ S.\ Retakh: Homotopical properties of 
categories of extensions, Russian Math. Surveys {\bf 41} (1986),
pp. 217--218.\\

M.\ Schlessinger: Functors of Artin rings. Trans{.} A.M.S{.}
{\bf 130} (1968).\\

K.\ H.\ Ulbrich, Group cohomology for Picard  categories. J. Alg.
{\bf 91} (1984).\\

\bye